\documentclass[12pt]{ociamthesis}
\usepackage{amsmath, amsthm, amsfonts, amscd, amssymb, epsfig, setspace}
\usepackage[cmtip,arrow]{xy}
\usepackage{pb-diagram,pb-xy}
\numberwithin{equation}{section}
\theoremstyle{plain}
\newtheorem{theorem}{Theorem}[section]
\newtheorem{lemma}[theorem]{Lemma}
\newtheorem{remark}[theorem]{Remark}
\newtheorem{corollary}[theorem]{Corollary}
\newtheorem{definition}[theorem]{Definition}
\newtheorem{proposition}[theorem]{Proposition}
\newtheorem{conjecture}[theorem]{Conjecture}

\newcommand{\myqed}{{\mbox{\rule{.1in}{.1in}}\par\vskip\baselineskip}}
\newcommand{\cO}{{\mathcal O}}
\newcommand{\cU}{{\mathcal U}}

\newcommand{\C}{{\mathbb C}}
\newcommand{\R}{{\mathbb R}}
\newcommand{\Z}{{\mathbb Z}}
\newcommand{\N}{{\mathbb{N}}}
\newcommand{\deebar}{\overline{\partial}}
\newcommand{\CP}{{\mathbb{P}^1}}
\newcommand{\QQ}{{\CP\times\CP\setminus\overline\Delta}}
\newcommand{\Q}{{\mathbf Q}}
\newcommand{\g}{{\mathfrak g}}
\newcommand{\cG}{{\mathcal G}}
\newcommand{\cA}{{\mathcal A}}
\newcommand{\T}{{\mathbf T}}
\newcommand{\Hyp}{{\bf H^3}}

\DeclareMathOperator{\supp}{supp}
\DeclareMathOperator{\End}{End}
\DeclareMathOperator{\Hom}{Hom}
\DeclareMathOperator{\rank}{rank}
\DeclareMathOperator{\Ker}{Ker}

\title{Differential Geometry of Monopole Moduli Spaces}
\author{Oliver Nash}
\college{Balliol College}
\degree{Doctor of Philosophy}
\degreedate{July 2006}

\begin{document}

\setcounter{secnumdepth}{3}
\setcounter{tocdepth}{3}

\maketitle
\begin{dedication}
   For my parents.
\end{dedication}

\begin{acknowledgements}
Throughout the course of my studies at Oxford, I have relied on the
help and support of many people. It is my pleasure to thank them here.
\par
Naturally I am most indebted to my supervisor, Nigel Hitchin.
It would be difficult to exaggerate how much he has done for me to
ensure that this thesis became a reality. Throughout the years
I worked with Nigel he demonstrated extraordinary
patience as well as exceptional generosity in terms of both time and 
ideas. I will always be grateful to Nigel with whom it has been a 
privilege to work.
\par
I am also very grateful to and for my wonderful family. Behind every
step of progress I made was the unwavering support of my parents
Charles and Edna as well as my sisters Elise and Mary and my brother
Nicholas.
\par
I could never have completed my studies at Oxford without the wonderful
friends who have always been beside me. I would like to thank them all,
especially those
from home in Dublin, my friends from the Mathematical Institute and
all my friends from my college, Balliol. Living in the community that is
Holywell Manor has been an amazing experience that I will, to say the
very least, find it difficult to leave behind.
\par
I would also like to thank my college advisors Keith Hannabuss and
Frances Kirwan for watching over me and regularly checking that all was
well. Balliol is lucky to have found two such remarkable people
in as many exceptional mathematicians.
\par
Finally I would like to thank the donors and sponsors of my scholarships:
The Scatcherd European Scholarship and the Foley-B\'ejar Scholarship. In
particular I would like to thank Martin Foley not only for his generous
financial aid but also for the genuine interest he consistently showed
in my progress.
\end{acknowledgements}

\begin{abstract}
   This thesis was motivated by a desire to understand the natural geometry
   of hyperbolic monopole moduli spaces. We take two approaches. Firstly
   we develop the twistor theory of singular hyperbolic monopoles and use
   it to study the geometry of their charge 1 moduli spaces. After this we
   introduce a new way to study the moduli spaces of both Euclidean and hyperbolic
   monopoles by applying Kodaira's deformation theory to the spectral curve.
   We obtain new results in both the Euclidean and hyperbolic cases. In
   particular we prove new cohomology vanishing theorems and find that the
   hyperbolic monopole moduli space appears to carry a new type of geometry
   whose complexification is similar to the complexification of hyperk\"ahler
   geometry but with different reality conditions.
\end{abstract}

\begin{romanpages}
   \tableofcontents
\end{romanpages}
\chapter{Overview and statement of results}
\section{A very short introduction to monopoles}
This thesis concerns non-Abelian magnetic monopoles, or just 
monopoles for short. Monopoles can be studied
on any oriented Riemannian 3-manifold but Euclidean space
and hyperbolic space are the two most important cases and
are the two which we shall study here. We
begin with a discussion of monopoles on Euclidean space, so
for now by a monopole we shall mean a Euclidean monopole.
\par
A monopole is a type of 3-dimensional soliton. More
precisely, a monopole is a pair $(A, \Phi)$ that satisfies
the Bogomolny equations
\begin{eqnarray*}
   F_A = *\nabla_A\Phi
\end{eqnarray*}
and satisfies certain boundary conditions which imply that
$F_A$ is $L^2$.
Here $A$ is a connection on a principal $G$-bundle $P$ over
$\R^3$ and $F_A$ is its curvature. $\Phi$ is a section of the
bundle of Lie algebras $adP$ over $\R^3$ that is associated to
$P$ by the adjoint action of $G$ on its Lie algebra, $\nabla_A$
is the covariant derivative operator induced on $adP$ and
$*$ is the Hodge $*$-operator on $\R^3$. For simplicity,
from now on we shall restrict ourselves to the case $G =
SU(2)$. Furthermore, we are only interested in solutions
$(A, \Phi)$ up to gauge equivalence; we identify
solutions that differ only by an automorphism of $P$.
\par
We have noted that the curvature of a monopole is $L^2$.
In fact it turns out that we have
\begin{eqnarray*}
   \int_{\R^3}\|F\|^2 = 4\pi k
\end{eqnarray*}
for an integer $k \ge 0$. This integer $k$ is called the charge
of the monopole. We exclude the trivial case $k=0$.
\par
By imposing spherical symmetry it is easy to write down
monopoles of charge $1$ but it is not obvious that there exist
solutions of higher charge. The first rigorous proof of existence
for monopoles of arbitrary charge was provided by Taubes who
showed how, roughly speaking,
a charge $k$ monopole may be constructed by gluing together
$k$ charge 1 monopoles. Furthermore, Taubes proved that the moduli
space (ie: the space of solutions up to gauge equivalence) of
monopoles of charge $k$ is a smooth manifold of dimension $4k - 1$.
\par
Although Taubes's results resolved the issue of existence, they
shed no light on how to solve the Euclidean Bogomolny equations.
Progress on this front was not long coming however. Hitchin,
using twistor theory, showed how a monopole corresponds to a
holomorphic vector bundle on $T\CP$ and using this introduced a
compact algebraic curve in $T\CP$, called the spectral curve, which
determines the monopole up to gauge equivalence. Around the same
time, Nahm, using a generalised Fourier transform (the Nahm
transform), showed how the Euclidean Bogomolny equations could be
reduced to the system (Nahm's equations)
\begin{eqnarray*}
   \frac{dT_i}{ds} = \epsilon_{ijk}[T_j, T_k]
\end{eqnarray*}
where $T_i$ is a $k\times k$ matrix valued function on
$(0,2)$ satisfying appropriate reality and boundary conditions.
Nahm also showed how to associate an algebraic curve to a
solution of these equations and Hitchin showed that this was the
spectral curve of the monopole. The spectral curve thus emerged as
a key feature in the theory of monopoles.
\par
After such successes solving the monopole equations, attention
turned to the study of their moduli spaces. We have already noted
that the moduli space of charge $k$ monopoles is a smooth manifold
of dimension $4k-1$. In fact there is a natural circle bundle
on the moduli space, called the gauged moduli space. Points in this
space should be thought of as a monopole together with an additional
phase factor. A major breakthrough was a theorem of Donaldson which
states that the charge $k$ gauged moduli space is naturally diffeomorphic
to the space of based degree $k$ rational maps $\CP \to \CP$.
\par
A natural parameterisation of the moduli space having been provided,
the next step was to understand the geometry of the moduli space. In
the course of his work on existence, Taubes had observed that the
gauged moduli space carries a natural quaternionic structure. The
natural question was whether this was integrable, ie: whether the
gauged moduli space was hyperk\"ahler. Atiyah and Hitchin answered
this question in the affirmative by showing how to view the gauged
moduli space as an infinite dimensional hyperk\"ahler quotient. They
also obtained explicit formulae for the metric on the centred charge
2 moduli space. Determining this metric was one of the great
achievements of Euclidean monopole theory.
\par
It was just as the study of Euclidean monopoles attained this
level of maturity that its hyperbolic younger brother was born.
Atiyah, observing that hyperbolic monopoles may be regarded as
$S^1$-invariant instantons on $S^4$, initiated their study. Because
of the link with instantons, some questions are much easier. For
example existence is trivial and an easy equivariant index
calculation shows that the moduli space of hyperbolic monopoles
is a smooth $4k-1$ dimensional submanifold of the space of
instantons.
\par
Hyperbolic monopoles have some key features in common with the
Euclidean variety. In particular they are also determined by a
spectral curve which in this case is a compact algebraic curve in
$\CP\times\CP$ and they too have a natural parameterisation in
terms of rational maps. However there are some important
differences. For example the Nahm transform reduces the hyperbolic
Bogomolny equations to a system of difference equations (whose
continuous limit is the Nahm equations above). A more striking
contrast was unveiled by Braam and Austin who showed that a
hyperbolic monopole is determined by its boundary value, the
$U(1)$ connection induced on the 2-sphere at infinity in hyperbolic
space. This is the exact opposite of what happens in the Euclidean
case where the boundary value depends only on the charge of the
monopole. Related to this is the fact that infinitesimal deformations
of hyperbolic monopoles are not $L^2$, ie: the natural $L^2$
\lq\lq metric\rq\rq on the moduli space diverges.
\par
In the Euclidean case, although it was a challenge to obtain the
aforementioned explicit formula for the centred charge 2 metric,
it was nevertheless clear from early on what \emph{type} of geometry
existed on the moduli spaces. By contrast, in the hyperbolic case
it remains unknown what type of geometry exists. This thesis was
motivated by a desire to determine this geometry.
\par
We attack this problem from two different points of view. Our first
approach was developed out of recognition of the need for more explicit
examples of hyperbolic monopole moduli spaces. Drawing on analogous
constructions in the Euclidean case, we introduce the class of 
\emph{singular} hyperbolic monopoles. We proceed to set up the relevant
twistor theory for their study and introduce the spectral data of a
singular hyperbolic monopole. After establishing the correspondence
between singular monopoles and their spectral data, we use it to identify
the charge 1 moduli spaces and study their geometry. The point of
introducing the singularities is that, because of their presence, even
the charge 1 moduli spaces carry interesting geometry while still
permitting explicit description. We find that these moduli spaces carry
a natural 2-sphere of conformally equivalent scalar-flat K\"ahler metrics.
\par
Our second approach exploits the fact that from the point of view
of the spectral curve, the theories of Euclidean and hyperbolic
monopoles are very similar. We show how to obtain the monopole moduli
spaces by applying Kodaira's deformation theory to the spectral curve in
an appropriate ambient space. Using this, we are able to recover the
natural hyperk\"ahler structure on the Euclidean monopole moduli spaces
from the point of view of the spectral curve. We go on to apply this
technique to study the geometry on the moduli spaces of hyperbolic
monopoles. We find that they carry a type of geometry whose
complexification is very similar to the complexification of
hyperk\"ahler geometry but which has different reality conditions.
It is, however, as yet unclear what type of real geometry this is the
complexification of.
\section{Summary of results}
{\bf Chapter \ref{SingHypChap}} contains our work on singular
hyperbolic monopoles. After providing the necessary definitions and
making some elementary observations, we embark on the task of setting
up the relevant twistor theory for their study. Using this we introduce
the spectral data of a singular monopole. This consists of a spectral
curve $S$ together with a divisor $D$ on $S$. $S$ is a compact algebraic
curve in the linear system $\cO(k, k)$ on $\QQ$ satisfying the
constraint (proposition \ref{sing_hyp_spec_crv_props} part (iii))
\begin{eqnarray*}
   L^{2m+k}(0,2l)|_S \simeq [D]
\end{eqnarray*}
where $m$ is the mass of the monopole and $k$, $l$
are charges. We prove that a singular monopole is determined up to
gauge equivalence by the pair $(S, D)$.
\par
After establishing this we go on to study the geometry of the charge
$1$ moduli spaces. Our results culminate in
\par\vskip0.5\baselineskip\noindent
{\bf Theorem \ref{Charge1ModSpaceGeom}.}
   {\sl Let $m > 0$, let $\{p_1,\ldots,p_n\} \subset \Hyp$ be $n$ distinct points in
   hyperbolic space and let $\{l_1,\ldots,l_n\} \subset \N$ be $n$ (strictly)
   positive integers. Let $M$ be the moduli space of gauged singular hyperbolic
   $SU(2)$ monopoles of non-Abelian charge 1 with Abelian charges $l_i$ at $p_i$.
   Then}
   \begin{itemize} {\sl
      \item
      $M$ carries a natural self-dual conformal structure
      \item
      For each point $u \in \partial\Hyp$ there is a volume form and complex
      structure $J^u$ on $M$ such that the metric determined in the conformal
      structure makes $M$ together with $J^u$ a scalar-flat K\"ahler manifold.}
   \end{itemize}
\par
In {\bf chapter \ref{DefSpecChap}} we introduce a new way to study the
moduli spaces of monopoles on both Euclidean and hyperbolic space. If
$S$ is the spectral curve of a monopole, the
idea is to use the triviality of $L^2|_S$ in the Euclidean case, or
$L^{2m+k}|_S$ in the hyperbolic case, to lift the spectral curve to a
curve $\hat S$ in $L^2 \setminus 0$ or $L^{2m+k}\setminus 0$. We then apply
Kodaira's deformation theory for compact complex submanifolds to $\hat S$
so that we obtain the complexified monopole moduli space as
a space of deformations. This means we have a model for the complexified
tangent space as $H^0(\hat S, \hat N)$ where $\hat N$ is the normal bundle of
$\hat S$ in its ambient space.
\par
Having set up the deformation theory we prove that, in the Euclidean case,
we have a natural decomposition of the complexified tangent space at $\hat S$:
\begin{eqnarray*}
   T_{\hat S}M_k\otimes\C \simeq H^0(\hat S, \hat N) \simeq H^0(\hat S,
   \hat N(-1))\otimes \C^2
\end{eqnarray*}
and that $H^0(\hat S, \hat N(-1))$ carries a natural quaternionic structure
and symplectic structure. This means that $H^0(\hat S, \hat N)$ is the
complexified tangent space to a hyperk\"ahler manifold and we show that
this is the usual hyperk\"ahler structure on the Euclidean monopole moduli
space. In the hyperbolic case we prove that there are \emph{two} corresponding
decompositions
\begin{eqnarray*}
   T_{\hat S}M_k \otimes\C \simeq H^0(\hat S, \hat N) &\simeq& H^0(\hat S,
   \hat N(-1, 0))\otimes \C^2\\
   T_{\hat S}M_k \otimes\C \simeq H^0(\hat S, \hat N) &\simeq& H^0(\hat S,
   \hat N(0, -1))\otimes \C^2
\end{eqnarray*}
Crucially, we identify a natural symplectic structure on both
$H^0(\hat S, \hat N(-1, 0))$ and $H^0(\hat S, \hat N(0, -1))$ but in place of
the quaternionic structure we obtain a conjugate linear isomorphism between
these two spaces.
\par
We must point out that for Kodaira's deformation theory to work,
it is necessary for a certain obstruction class to vanish and this is
usually established by proving the vanishing of a certain cohomology group.
We must prove such a vanishing theorem for our deformation theory to work
and much of the effort in this approach goes into proving this vanishing.
In the Euclidean case the vanishing theorem we need and prove is
\par\vskip0.5\baselineskip\noindent
{\bf Theorem \ref{Euc_vanishing_theorem}.} {\sl
   Let $\widetilde E$ be the holomorphic vector bundle on $\T$
   corresponding to an $SU(2)$ monopole of charge $k$ on $\R^3$. Let
   $S$ be the spectral curve of the monopole. Then
   \begin{eqnarray*}
      H^0(S, \widetilde EL(k-2)) = 0
   \end{eqnarray*}}
\indent
Although the statement of this theorem is similar to that of Hitchin's
crucial vanishing result $H^0(S, L^z(k-2)) = 0$ for $z\in (0,2)$, the proof
requires new ideas and is not merely a simple adaptation of known techniques.
A key observation is that there is a natural injection $H^0(S, \widetilde
EL(k-2)) \hookrightarrow H^0(\T, \cO(2k-2))$ which means that we may
think of our sections as polynomials.
\par
Similarly in the hyperbolic case we prove
\par\vskip0.5\baselineskip\noindent
{\bf Theorem \ref{hyp_vanishing}.} {\sl
   Let $\widetilde E$ be the holomorphic vector bundle on $\Q$
   corresponding to an $SU(2)$ monopole of charge $k$ on $\Hyp$. Let
   $S$ be the spectral curve of the monopole. Then
   \begin{eqnarray*}
      H^0(S, \widetilde EL^m(k-1, -1)) = 0
   \end{eqnarray*}}
\indent
Finally, we take a brief detour to answer a question that arose while
considering the analysis necessary to prove the vanishing theorem
mentioned above. Specifically, we prove that (for both Euclidean and
hyperbolic monopoles) the Penrose transform of the Higgs field can be
interpreted as a component of the Atiyah class of the holomorphic
vector bundle on twistor space corresponding to the monopole.
\par
We bring the thesis to a close with {\bf Chapter \ref{OpenIssueChap}},
a slightly more speculative chapter than its predecessors.
\par
In section \ref{HCInstQuotSect} we prove that the hypercomplex
structure (observed by Joyce and Teleman)
on the moduli spaces of instantons on a compact 4-dimensional
hypercomplex manifold may be obtained, at least formally, by applying
an infinite dimensional hypercomplex quotient construction introduced
by Joyce. This result is of interest independently of questions about
monopoles but we also outline some features that the geometry of
instantons on certain Hopf surfaces has in common with the
geometry of hyperbolic monopoles.
\par
In section \ref{HypMonop_KahMetric_sect}, we show how to obtain a
class of natural K\"ahler metrics on hyperbolic monopole moduli
spaces which are higher dimensional generalisations of the
metrics introduced in chapter \ref{SingHypChap}.
\par
In section \ref{TwistSpace_Section}, we outline how to construct
twistor spaces for the moduli spaces of hyperbolic monopoles. We find
that instead of a single twistor space together with a real structure,
we obtain two twistor spaces together with an anti-holomorphic
bijection between them. For this reason the reality condition on
the twistor lines is more complicated.

\chapter{Singular hyperbolic monopoles}\label{SingHypChap}
\section{Overview}
Singular monopoles were first studied by Kronheimer in \cite{Kronheimer}.
He showed how to extend the twistor theory of non-singular Euclidean
monopoles developed in \cite{MR649818} to allow for monopoles with a finite
number of prescribed singularities. Furthermore, he then used these
techniques to construct the twistor space of the charge $1$ moduli space
of singular monopoles on $\R^3$ and thus found that it carried a natural
hyperk\"ahler structure (indeed the moduli space is the smooth part of a
quotient of multi--Taub--NUT space).
\par
Here we offer a hyperbolic version of Kronheimer's work. That is, we show
how to extend the twistor theory of non-singular hyperbolic monopoles
developed in \cite{MR1399482} to allow for monopoles on hyperbolic space
with a finite number of prescribed singularities. Our main reason for
developing this theory is that, just as in the Euclidean case, this allows
us to construct the twistor space of the lowest dimensional moduli spaces
and hence study their geometry. As stated in the introduction, we find that
the moduli space of charge $1$ singular hyperbolic monopoles possesses a
natural $2$-sphere of scalar-flat K\"ahler metrics all within the same
conformal class. The $2$-sphere appears naturally as the boundary of
hyperbolic space.
\par
The motivation for this work is to understand the natural geometry of
hyperbolic monopole moduli spaces. Since Kronheimer's work made it possible
to see the natural hyperk\"ahler geometry of the charge $1$ singular
Euclidean monopole moduli spaces in a very explicit way it is natural to ask
the same question in the hyperbolic case.
\par
The geometric structure identified on the moduli space of charge 1 monopoles
perhaps deserves additional interest owing to the fact that a limiting case
of it has already been studied in some detail by LeBrun \cite{MR1114461}. As
we shall see, the spaces studied by LeBrun correspond to the zero mass limit
of our monopole moduli spaces.
\par
Finally we must point out that singular monopoles have recently been
studied by Kapustin and Witten in \cite{hep-th0604151} as part of
their work on the Geometric Langlands Programme. Although the context of
their work is very different to ours, \cite{hep-th0604151} nevertheless
provides a further justification for the study of singular monopoles.
\section{Definitions and elementary properties}\label{SngMonopDefsSect}
The monopoles we are interested in are solutions of the Bogomolny equations
with singularities at a fixed set of points in hyperbolic space. The definition
below gives the precise behaviour of the solutions at these points.
\begin{definition}\label{def_sngmpole}
   Let $\{p_1, \ldots ,p_n\} \subset \Hyp$ be $n$ distinct points in
   hyperbolic space. Let $U = \Hyp \setminus \{p_1, \ldots ,p_n\}$ and let
   $\pi : E \to U$ be a $C^\infty$ $SU(2)$ vector bundle on $U$. A singular
   hyperbolic $SU(2)$ monopole with singularities at $p_1, \ldots ,p_n$ is
   an $SU(2)$ connection
   \begin{eqnarray*}
      \nabla : \Omega^0(U, E) \to \Omega^1(U, E)
   \end{eqnarray*}
   and an $SU(2)$ endomorphism (the Higgs field) $\Phi \in \Omega^0(U,
   \End(E))$ such that:
   \begin{enumerate}
      \item
      $\Phi$ and $\nabla$ satisfy the Bogomolny equations:
      \begin{eqnarray*}
         \nabla \Phi = *F_\nabla
      \end{eqnarray*}
      (where $F_\nabla \in \Omega^2(U, \End(E))$ is the curvature of $\nabla$)
      \item
      $(\nabla, \Phi)$ satisfy the boundary conditions BC0, BC1, BC2 defined in
      \cite{MR1399482}. If we fix a point $O \in \Hyp$ then (for $SU(2)$ monopoles)
      these conditions can be described as follows. Let
      \begin{eqnarray*}
         \left[
         \begin{array}{cc}
            A & B\\
            -B^* & -A
         \end{array}
         \right]
      \end{eqnarray*}
      be the connection matrix of $\nabla$ in a gauge of unitary eigenvectors of
      $\Phi$. If $\rho$ is the hyperbolic distance from $O$ then we require
      \begin{itemize}
         \item
         $(\|\Phi\| - m)e^{2\rho}$ extends smoothly to $\partial\Hyp$ for some $m > 0$
         \item
         $A$ extends smoothly to $\partial\Hyp$
         \item
         $Be^{2m\rho}$ extends smoothly to $\partial\Hyp$
      \end{itemize}
      \item
      $\Phi$ has the following behaviour at singular points:
      \begin{eqnarray*}
         \lim_{\rho_i \to 0}\limits \left(\rho_i\|\Phi\|\right) \mbox{ exists and
         is (strictly) positive}\\
         d\left(\rho_i\|\Phi\|\right) \mbox{ is bounded in a neighbourhood of $p_i$}
      \end{eqnarray*}
      where $\rho_i$ is the hyperbolic distance from $p_i$.
   \end{enumerate}
\end{definition}
\begin{remark}
   We have chosen to define a monopole in the language of vector
   bundles. It may thus be useful (at least for the sake of fixing
   notation) to recall that by a $C^{\infty}$ $SU(2)$ vector bundle
   over $U$ is meant a rank $2$ complex vector bundle together
   with a symplectic form $\chi \in C^{\infty}(U, \wedge^2 E)$ and a
   quaternionic structure $j : E \to E$ (ie: $j$ is anti-linear on the fibres,
   covers the identity on $U$ and $j^2 = -1$) such that
   \begin{eqnarray*}
      \chi(jv,w) = \overline{\chi(v,jw)}
   \end{eqnarray*}
   and such that the Hermitian scalar product on $E$ defined by
   \begin{eqnarray*}
      (v,w) = i\chi(jv,w)
   \end{eqnarray*}
   is positive definite. A connection $\nabla$ on $E$ is an $SU(2)$ connection
   iff $\nabla\chi = \nabla j = 0$ and an endomorphism $\Phi$ of $E$ is an
   $SU(2)$ endomorphism iff $\lbrack\Phi, j\rbrack =0$ and $\Phi$ is skew
   adjoint with respect to $\chi$.
\end{remark}
\par
The boundary conditions (iii) of definition \ref{def_sngmpole} deserve
elaboration. To see where they come from, let $V : U \to \R$ be
\begin{eqnarray}\label{V_lambda_def}
   V = \lambda + \sum_{i=1}^n G_{p_i}
\end{eqnarray}
for some $\lambda \ge 0$ and where
\begin{eqnarray}\label{Hyp_Lapl_Green}
   G_p(x) = \frac{1}{e^{2\rho(p,x)} - 1}
\end{eqnarray}
is the Green's function for the hyperbolic Laplacian centred at $p$ and
normalised so that $\Delta G_p = -2\pi\delta_p$. Let $M$ be the principal
$U(1)$ bundle on $U$ with Chern class $\frac{1}{2\pi}[*dV]$. Let $\omega$ be a
connection on $M$ with curvature $\frac{1}{2\pi}*dV$ and define the
metric $g$ on $M$ by
\begin{eqnarray*}
   g = Vg_\Hyp + V^{-1}\omega\otimes\omega
\end{eqnarray*}
We give $M$ the orientation defined by $vol_{\Hyp}\wedge\omega$.
As shown in \cite{MR1114461}, we may add in a fixed point $\hat p_i$ of the
$S^1$ action on $M$ over each $p_i \in \Hyp$ to obtain a smooth manifold
$\hat M = M \cup \{\hat p_1, \ldots, \hat p_n\}$ and the metric extends smoothly
to $\hat M$.  Now if $(\nabla, \Phi)$
is a connection and Higgs field on $U$ then (suppressing the notation for pull
backs) we define a connection $\hat\nabla$ on $M$ according to the correspondence
\begin{eqnarray}\label{Bog_Inst_corresp}
   (\nabla, \Phi) \mapsto \hat\nabla = \nabla - V^{-1}\Phi\otimes\omega
\end{eqnarray}
Using the formulae
\begin{eqnarray*}
   F_{\hat\nabla} &=& F_\nabla - V^{-1}\Phi\otimes d\omega + V^{-2}\Phi\otimes
   dV\wedge\omega - V^{-1}\nabla\Phi\wedge\omega
\end{eqnarray*}
and
\begin{eqnarray*}
   d\omega &=& *dV\\
   \hat *\alpha &=& V^{-1}(*\alpha)\wedge\omega \quad\mbox{for $\alpha \in \wedge^2 T^*U$}\\
   \hat *(\alpha\wedge\omega) &=& V*\alpha \quad\mbox{for $\alpha \in T^*U$}
\end{eqnarray*}
(where $\hat *$ denotes the Hodge $*$-operator on $M$ and $*$ is the Hodge
$*$-operator on $U$)
it follows that $(\nabla, \Phi)$ satisfy the Bogomolny equations
on $U$ iff $\hat\nabla$ satisfies the anti-self-dual Yang--Mills equations on
$M$. We thus have a correspondence between solutions of the Bogomolny equations
on $U$ and $S^1$-invariant solutions of the anti-self-dual Yang--Mills equations
on $M$. The promised elucidation of the aforementioned boundary conditions can
now be stated as
\begin{lemma}\label{sngmpole_inst_link}
   In the above notation, $(\nabla, \Phi)$ satisfy the boundary conditions
   (iii) of definition \ref{def_sngmpole} iff the corresponding solution of
   the anti-self-dual Yang-Mills equations on $M$ extends to a solution on $\hat M$.
\end{lemma}
\noindent {\bf Proof} The proof is completely analogous to the corresponding
result in \cite{Kronheimer}. \myqed
Now if we have an $S^1$ invariant instanton on a bundle $E \to \hat M$, the
fibre of $E$ over the point $\hat p_i \in \hat M$ lying above a singular point
$p_i \in \Hyp$ will carry a representation of $S^1$. Since the $S^1$ action is
compatible with the $SU(2)$ structure of $E$, this action must have weights
$(l_i, -l_i)$ for some integer $l_i \ge 0$. The question arises of identifying this
integer $l_i$ in terms of the corresponding solution of the Bogomolny equations on $U$.
In fact
\begin{eqnarray}\label{Ab_charge_formula}
   l_i = 2\lim_{p\to p_i}\left(\rho(p,p_i)\|\Phi(p)\|\right)
\end{eqnarray}
To see why, fix a trivialisation of $E$ in a neighbourhood of $\hat p_i$. Let $A$
be the corresponding matrix of $1$-forms and let $A_0 = A(X)$ where $X$ is the
vector field on $\hat M$ generated by the $S^1$ action. Now for each $p \in \Hyp$
near $p_i$ choose a gauge transformation $g : M_p \to SU(2)$ on the corresponding
$S^1$ orbit that takes our fixed trivialisation to an $S^1$-invariant one. In view
of \eqref{Bog_Inst_corresp} we thus have
\begin{eqnarray*}
   -V^{-1}\Phi = g^{-1}A_0g + g^{-1}X(g)
\end{eqnarray*}
Now as $p\to p_i$
\begin{eqnarray*}
   \|g^{-1}A_0g\| = \|A_0\| \to 0
\end{eqnarray*}
since $X$ vanishes at $p_i$. Furthermore
\begin{eqnarray*}
   \|g^{-1}X(g)\| \to l_i
\end{eqnarray*}
since $g$ is approaching the $S^1$ representation with weights $(l_i, -l_i)$. Equation
\eqref{Ab_charge_formula} now follows upon noting that $\lim_{p\to p_i}\limits V^{-1}\|\Phi\|
= 2 \lim_{p\to p_i}\limits \rho_i\|\Phi\|$ since $\lim_{p\to p_i}\limits 2\rho_i V = 1$.
\begin{definition}
   In the above notation, we define the Abelian charge $l_i$ of the monopole at $p_i$
   by equation \eqref{Ab_charge_formula}. We also define the total Abelian charge
   $l$ of the monopole as $l = \sum_{i=1}^{n}\limits l_i$.
\end{definition}
\begin{definition}
   Let $O \in \Hyp$. From condition (ii) of definition \ref{def_sngmpole} the limit
   \begin{eqnarray*}
      m = \lim_{\rho(p,O)\to \infty}\|\phi(p)\| \in \R
   \end{eqnarray*}
   exists and is (strictly) positive. We define $m$ to be the mass of the monopole.
\end{definition}
Fix a point $O \in \Hyp$. Since $\|\Phi(p)\| \to m > 0$ as $\rho(p, O) \to \infty$
we can choose a sphere $S$ in $\Hyp$ centred at $O$ large enough that the singular
points and zeros of $\Phi$ all lie inside $S$. On such a sphere, the bundle $E$
splits as a direct sum of eigenbundles of $\Phi$
\begin{eqnarray*}
   E|_S = M^+\oplus M^-
\end{eqnarray*}
where $M^\pm$ is the bundle corresponding to the eigenvalue $\pm i\|\Phi\|$. (Note
that these bundles are interchanged by $j$.)
\begin{definition}
   In the above notation and using the natural orientation of $S$, we define the
   total charge $N$ of the monopole by
   \begin{eqnarray*}
      N = c_1(M^+)[S]
   \end{eqnarray*}
\end{definition}
\begin{definition}
   We define the non-abelian charge $k$ of a monopole to be $k = N + l$.
\end{definition}
It is important to address the issue of existence of singular
hyperbolic monopoles. As we shall see, the key is a result of LeBrun in
\cite{MR1114461}.
\par
We have seen that given a harmonic function $V$ on $U$ as in \eqref{V_lambda_def} we
obtain the Riemannian manifold $\hat M$ of lemma \ref{sngmpole_inst_link}.
The function $V$ depends on a choice of $\lambda \ge 0$ and LeBrun \cite{MR1114461}
shows that for $\lambda = 1$, $\hat M$ has an $S^1$-equivariant conformal
compactification $M^c$ obtained by adding a $2$-sphere of fixed points of
the $S^1$ action on $\hat M$ and gluing along the boundary of $\Hyp$
(which $\hat M$ fibres over). Furthermore, after reversing the orientation,
$M^c$ is diffeomorphic to $n\mathbb{CP}^2 = \mathbb{CP}^2 \# \cdots \# \mathbb{CP}^2$
and the conformal class (which is self-dual) contains a metric of positive scalar
curvature. For $n=0$ (ie: no singularities) this construction is of course the usual
observation that round $S^4$ is an $S^1$-equivariant conformal compactification of
$\Hyp\times S^1$ which was used very successfully by Atiyah in \cite{MR893593} to
study monopoles on $\Hyp$. For $n=1$ we obtain $\mathbb{CP}^2$ with the usual
Fubini-Study conformal structure.
\par
Now in \cite{MR1399482}, it is noted that an $S^1$-invariant instanton
on $S^4$ corresponds to a solution of the Bogomolny equations on $\Hyp$
that satisfies the boundary conditions (ii) of definition \ref{def_sngmpole}.
Similarly and in view of lemma \ref{sngmpole_inst_link}, an $S^1$-invariant
anti-self-dual instanton on $M^c$ corresponds to a solution of the
Bogomolny equations satisfying conditions (ii) and (iii) of \ref{def_sngmpole}.
This is the same as a self-dual instanton on $n\mathbb{CP}^2$ (we don't have
to be careful whether we consider anti-self-dual or self-dual instantons on
$S^4$ since it carries an orientation reversing diffeomorphism).
Existence of our singular monopoles then follows from the existence of
$S^1$-invariant self-dual instantons on the self-dual manifolds $n\mathbb{CP}^2$.
Indeed a careful equivariant index calculation can be used to calculate
the dimension of the moduli space, $4k - 1$.
\par
In fact, in view of Buchdahl's construction \cite{MR857374} of instantons on
$\mathbb{CP}^2$, it should even be possible to obtain explicit formulae for
singular hyperbolic monopoles just as the same is possible by applying an
$S^1$-invariant version of the ADHM construction for instantons on $S^4$ to
obtain formulae for non-singular hyperbolic monopoles. 
\section{The Hitchin--Ward correspondence}
The Hitchin--Ward transform is the fundamental theorem that tells us how
to interpret solutions of the Bogomolny equations on twistor space. In this
section we address the question of what happens to the data on twistor space
when the solutions of the Bogomolny equations have singularities as prescribed
in definition \ref{def_sngmpole}. Before stating the theorem, we find it
convenient to introduce some terminology and make some elementary observations
about hyperbolic space.
\begin{definition}
   Given an oriented geodesic $\gamma$ in hyperbolic space and a point $O \in \Hyp$
   there exists a unique parameterisation of $\gamma$ such that $\gamma$ is
   parameterised by arc length and $\gamma(0)$ is the closest point to $O$
   on $\gamma$. We call this the parameterisation of $\gamma$ determined by
   $O \in \Hyp$.
\end{definition}
\begin{lemma}\label{safe_radius}
   Let $\{p_1, \ldots ,p_n\} \subset \Hyp, O \in \Hyp$. Let $R > 0$
   be large enough that $\{p_1, \ldots ,p_n\} \subset B(O, R)$, let $\gamma
   : \R \to \Hyp$ be a geodesic with the parameterisation determined by $O \in
   \Hyp$ and let $|t| \ge R$. Then $\gamma(t) \notin B(O, R)$.
\end{lemma}
\noindent{\bf Proof} Even in hyperbolic space, the hypotenuse of a right angled
triangle is the longest side. \myqed
We shall denote the twistor space (ie: the set of oriented geodesics) of $\Hyp$ by
$\Q$. If $x \in \Hyp$, we shall denote the corresponding twistor line (the set of
all geodesics passing through $x$) in $\Q$ by $P_x$. Finally, if we fix a point
$O \in \Hyp$, then we can identify
\begin{eqnarray*}
\Q\simeq \QQ
\end{eqnarray*}
where (the so called anti-diagonal) is
\begin{eqnarray*}
   \overline\Delta = (p, \tau (p))
\end{eqnarray*}
and $\tau$ is the usual (anti-podal) real structure on $\CP$. The diagonal $\Delta
\subset \QQ$ appears as the twistor line of the chosen point $O \in \Hyp$. 
\par
Now consider the natural double fibration
\begin{eqnarray*}
   \begin{diagram}
      \node[2]{S\Hyp} \arrow{sw,t}{\nu} \arrow{se,t}{\mu}\\
      \node{\Hyp} \node[2]{\Q}
   \end{diagram}
\end{eqnarray*}
where $S\Hyp \subset T\Hyp$ is the unit tangent bundle of $\Hyp$.
Using the hyperbolic metric, we have $T\Hyp \simeq T^{*}\Hyp$ and so pulling
back the natural 1-form on $T^{*}\Hyp$, $T\Hyp$ and hence $S\Hyp$ carries a
natural 1-form
\begin{eqnarray}\label{theta_hat_def}
   \hat\theta \in \Omega^1(S\Hyp)
\end{eqnarray}
Thus if $\hat f : \Q \to S\Hyp$ is a section of $\mu$ we obtain a 1-form
$\theta = \hat f^{*}\hat\theta$ on $\Q$.
\begin{remark}\label{theta01_def}
   A point $O\in\Hyp$ determines a section of $\mu$, namely $\gamma \mapsto
   \dot\gamma(0)$ in the parameterisation of $\gamma$ determined by $O$. Thus,
   by the above $\QQ$ carries a natural 1-form $\theta$.
   As we shall see later, it is really the $(0,1)$ component of $\theta$
   that interests us. Using the coordinates, $(\lbrack z,1\rbrack,\lbrack
   w,1\rbrack)$ on (an open set of) $\QQ$, the explicit formula for the $(0,1)$
   component of $\theta$ is:
   \begin{eqnarray}\label{theta_formula}
      \theta^{0,1} = (z-w)\left(\frac{d\bar z}{(1+z\bar z)(1+\bar zw)} +
      \frac{d\bar w}{(1+w\bar w)(1+z\bar w)}\right)
   \end{eqnarray}
   Note that $\deebar \theta^{0,1} = 0$ and that $\theta$ vanishes on
   the twistor line $\Delta$ (where $z=w$) in ${\QQ}$ corresponding to those geodesics passing
   through $O$, and has a singularity along $\overline\Delta$ (where $z = -1/\overline w$).
\end{remark}
\begin{lemma}\label{lemma_theta_hat}
   Let $\hat f$ be a section of $\mu$ above and let $f = \nu\circ\hat f : \Q \to
   \Hyp$. Let $\theta = \hat f^{*}\hat\theta$ and let $\omega \in \wedge^2 T^{*}
   \Hyp$. Then
   \begin{eqnarray*}
      (f^{*}\omega)^{0,2} + i\theta^{0,1}\wedge (f^{*}(*\omega))^{0,1} = 0
   \end{eqnarray*}
   where $*$ is the Hodge star on $\Hyp$.
\end{lemma}
\noindent{\bf Proof} This follows from results in \cite{MR1399482}. (In particular
see equation (3.5) of \cite{MR1399482}). \myqed
\begin{remark}
   The above results also hold with $\R^3$ in place of $\Hyp$. In the case of
   $\R^3$ if we use the usual coordinates $(\eta,\zeta)\mapsto \eta\frac{\partial}
   {\partial\zeta}$ on its twistor space $\T \simeq TS^2$, then the formula for
   its natural 1-form is:
   \begin{eqnarray}\label{Euc_theta01_def}
      \theta_{\T}^{0,1} = \frac{2\eta}{(1+\zeta\bar\zeta)^2}d\bar\zeta
   \end{eqnarray}
   This 1-form is also the 1-form obtained by using the usual round metric
   on $S^2$ to obtain $TS^2 \simeq T^{*}S^2$ and pulling back the
   natural 1-form on $T^{*}S^2$.
\end{remark}
\begin{definition}
   Let $\{p_1, \ldots ,p_n\} \subset \Hyp$ be $n$ distinct points in hyperbolic
   space and let $x \in \Hyp$. Suppose that there exists a geodesic $\gamma \in
   P_x$ and $p_i, p_j \in \gamma$ such that $x$ separates $p_i$ and $p_j$ on
   $\gamma$. Then we say $x$ is geodesically trapped by $\{p_1, \ldots ,p_n\}$.
\end{definition}
We are finally ready to state the Hitchin--Ward correspondence. The proof used
here owes most to the proof of the corresponding result in \cite{MR1399482}.
\begin{theorem}\label{hw_bogo}
   Let $\{p_1, \ldots ,p_n\} \subset \Hyp$ be $n$ distinct points in hyperbolic
   space. Let $U = \Hyp \setminus \{p_1, \ldots ,p_n\}$. Let $P = P_1 \cup \cdots
   \cup P_n \subset \Q$ (where $P_i = P_{p_i}$). Then to each solution of the
   $SU(2)$ Bogomolny equations on $U$, there corresponds a pair of rank 2
   holomorphic vector bundles $(E^+, E^-)$ on $\Q$ together with an isomorphism
   \begin{eqnarray*}
      h : E^+|_{\Q \setminus P} \to E^-|_{\Q \setminus P}
   \end{eqnarray*}
   of holomorphic vector bundles such that:
   \begin{enumerate}
      \item\label{hw_cond1}
      If $x \in U$ is a point that is not geodesically trapped by $\{p_1, \ldots 
      ,p_n\}$, then there exists a partition of $P_x \cap P$ into two disjoint
      sets: $Q^{+}_x, Q^{-}_x$ such that $\widetilde E^x$ is naturally isomorphic
      to the trivial vector bundle with fibre $E_x$, where $\widetilde E^x$ is
      the vector bundle over $P_x$ obtained by gluing $E^+|_{P_x \setminus Q^{+}_x}$
      and $E^-|_{P_x \setminus Q^{-}_x}$ together over $P_x \setminus P$ using $h$.
      \item\label{hw_cond2}
      $E^{\pm}$ carry holomorphic symplectic structures compatible with $h$
      \item\label{hw_cond3}
      There exists an anti-holomorphic (anti-linear) map
      \begin{eqnarray*}
         \tilde j : E^+\to E^-
      \end{eqnarray*}
      (covering $\sigma : \Q \to \Q$) such that
      \begin{eqnarray*}
         (h^{-1} \tilde j)^2 = -1
      \end{eqnarray*}
      over $\Q\setminus P$.
   \end{enumerate}
      Furthermore the holomorphic data determines the solution of the Bogomolny
      equations.
\end{theorem}
\noindent{\bf Proof} The essence of the theorem is by now standard. We sketch
   a proof that emphasises the differences that arise because of the singular
   points $p_i$.
   \par
   Thus, fix a point $O \in \Hyp$ and let $R > 0$ be large enough that
   $\{p_1, \ldots ,p_n\} \subset B(O, R)$. Define
   \begin{eqnarray*}
      f : \Q \to \Hyp
   \end{eqnarray*}
   by
   \begin{eqnarray*}
      \gamma \mapsto \gamma(R)
   \end{eqnarray*}
   where $\gamma$ is given the parameterisation determined by $O \in \Hyp$.
   Note that in view of lemma \ref{safe_radius}, we have $f(\Q) \subset U
   \subset \Hyp$ and so it makes sense to define\footnote{Although $E^+$
   depends on $R$, different values of $R$ give naturally isomorphic bundles.}
   \begin{eqnarray*}
      E^+ = f^{*}E
   \end{eqnarray*}
   and
   \begin{eqnarray*}
      \deebar : \Omega^0(\Q, E^+) \to \Omega^{0,1}(\Q, E^+)
   \end{eqnarray*}
   by
   \begin{eqnarray*}
      \deebar s = ((f^*\nabla) s - i(f^*\Phi)(s)\otimes\theta)^{0,1}
   \end{eqnarray*}
   where
   $\theta = \hat f^*\hat\theta$, $\hat\theta$ is the 1-form of equation
   \eqref{theta_hat_def} and $\hat f$ is the map:
   \begin{eqnarray*}
      \hat f : \Q &\to& S\Hyp\\
      \gamma &\mapsto& \dot\gamma(R)
   \end{eqnarray*}
   We thus have $\deebar^2$ = $F_{\hat\nabla}^{0,2}$ where $F_{\hat\nabla}$ is
   the curvature of the connection $\hat\nabla = f^*\nabla - if^*\Phi\otimes
   \theta$. But
   \begin{eqnarray*}
      F_{\hat\nabla} = f^*F_\nabla + i\theta\wedge f^*(\nabla\Phi) - if^*\Phi
      \otimes d\theta
   \end{eqnarray*}
   Using the Bogomolny equations $F_\nabla = *\nabla\Phi$ and the fact that
   $\deebar\theta^{0,1} = 0$ we find
   \begin{eqnarray*}
      \deebar^2 = (f^* F_\nabla)^{0,2} + i\theta^{0,1}\wedge (f^*(*F_\nabla))^{0,1}
   \end{eqnarray*}
   In view of lemma \ref{lemma_theta_hat} we thus have $\deebar^2 = 0$ and so
   we have a holomorphic structure on $E^+$. Define the holomorphic bundle
   $E^-$ using the same construction as for $E^+$ but with $f\circ\sigma$ in place
   of $f$. Note that we thus have $E^- = \sigma^* E^+$ as complex (but obviously
   not holomorphic) vector bundles.
   \par
   To define
   \begin{eqnarray*}
      h : E^+|_{\Q \setminus P} \to E^-|_{\Q \setminus P}
   \end{eqnarray*}
   note that if $\gamma \in \Q$ then $\sigma(\gamma)$ is just $\gamma$
   parameterised in the opposite direction, ie:
   \begin{eqnarray*}
      \sigma(\gamma) : t \mapsto \gamma(-t)
   \end{eqnarray*}
   Thus
   \begin{eqnarray*}
      E^-_{\gamma} \simeq E^+_{\sigma(\gamma)} \simeq E_{\gamma(-R)}
   \end{eqnarray*}
   and so to define $h$ we must define an isomorphism:
   \begin{eqnarray*}
      h_{\gamma} : E_{\gamma(R)} \simeq E_{\gamma(-R)}
   \end{eqnarray*}
   for all $\gamma \in \Q\setminus P$. Thus fix $\gamma \in \Q \setminus P$ and
   let $v \in E_{\gamma(R)}$. Note that $\gamma \subset U$ and so let $s$ be the
   unique section of $E$ along $\gamma$ such that
   \begin{eqnarray*}
      s(\gamma(R)) = v
   \end{eqnarray*}
   and
   \begin{eqnarray*}
      (\nabla_{\dot \gamma} - i\Phi)s = 0
   \end{eqnarray*}
   We define
   \begin{eqnarray*}
      h_{\gamma}(v) = s(\gamma(-R))
   \end{eqnarray*}
   It is straightforward to check that $h$ is indeed a holomorphic bijection.
   \par
   To see that condition \ref{hw_cond1} holds, let $x \in U$ be a point that
   is not geodesically trapped by $\{p_1, \ldots, p_n\}$. Let
   \begin{eqnarray*}
      Q^{\pm}_{x} = \{\gamma \in P_x\quad|\quad\mbox{there exists $p_i \in \gamma$
      separating $x$ and $\gamma(\pm R)$ on $\gamma$\}}
   \end{eqnarray*}
   Note that $Q^{+}_{x} \cup Q^{-}_{x} = P \cap P_x$ and $Q^{+}_{x} \cap Q^{-}_{x}
   = \emptyset$ since $x$ is not geodesically trapped.
   We define
   \begin{eqnarray*}
      \psi^{+} : E^+|_{P_x\setminus Q^{+}_{x}} \to (P_x\setminus Q^{+}_{x})
      \times E_x
   \end{eqnarray*}
   as follows. Let $\gamma \in P_x\setminus Q^{+}_{x}$ and let $\gamma^{+}$
   be the closed segment of $\gamma$ joining $x$ and $\gamma(R)$. Let $v \in
   E^+_{\gamma} \simeq E_{\gamma(R)}$ and let $s$ be the unique section of
   $E$ over $\gamma^{+} \subset U$ such that
   \begin{eqnarray*}
      s(\gamma(R)) = v
   \end{eqnarray*}
   and
   \begin{eqnarray*}
      (\nabla_{\dot \gamma} - i\Phi)s = 0
   \end{eqnarray*}
   Define
   \begin{eqnarray*}
      \psi^{+}(v) = (\gamma, s(x)) \in (P_x\setminus Q^{+}_{x}) \times E_x
   \end{eqnarray*}
   Similarly define
   \begin{eqnarray*}
      \psi^{-} : E^-|_{P_x\setminus Q^{-}_{x}} \to (P_x\setminus Q^{-}_{x})
      \times E_x
   \end{eqnarray*}
    by using a section $s$ over the closed segment $\gamma^{-}
   \subset U$ of $\gamma$ joining $x$ and $\gamma(-R)$ such that
   \begin{eqnarray*}
      s(\gamma(-R)) = v
   \end{eqnarray*}
   (since $E^-_{\gamma} \simeq E_{\gamma(-R)}$). Then define
   \begin{eqnarray*}
      \psi : \widetilde E^{x} \to P_x \times E_x
   \end{eqnarray*}
   by
   \begin{eqnarray*}
      [v] \mapsto \left\{
      \begin{array}{ll}
         \psi^{+}(v) & \mbox{if $v \in E^+_{P_x \setminus Q^{+}_{x}}$}\\
         \psi^{-}(v) & \mbox{if $v \in E^-_{P_x \setminus Q^{-}_{x}}$}
      \end{array} \right.
   \end{eqnarray*}
   It is straightforward to verify that $\psi$ is well defined and is the required
   trivialisation.
   \par
   To see that condition \ref{hw_cond2} holds, note that $E$ carries a
   symplectic structure
   \begin{eqnarray*}
      \chi \in C^{\infty}(U, \wedge^2 E^{*})
   \end{eqnarray*}
   Using $f$, we pull this back to a symplectic structure
   \begin{eqnarray*}
      \tilde\chi = f^{*}\chi \in C^{\infty}(\Q, \wedge^2 {E^+}^{*})
   \end{eqnarray*}
   Similarly, $E^-$ carries a symplectic structure. It is straightforward
   to verify that these are holomorphic and compatible with $h$.
   \par
   Finally for condition \ref{hw_cond3} let $j : E \to E$ be the
   quaternionic structure carried by $E$ and $\hat\sigma : E^- \to E^+$
   be the bijection induced by $\sigma : \Q \to \Q$. We define
   \begin{eqnarray*}
      \tilde j : E^+ \to E^-
   \end{eqnarray*}
   by
   \begin{eqnarray*}
      \tilde j = \hat \sigma^{-1} \circ (f^{*}j)
   \end{eqnarray*}
   Again it is straightforward to verify that $\tilde j$ has the required
   properties.
   \par
   We shall omit the proof that the holomorphic data determines the solution
   to the Bogomolny equations. \myqed
\begin{remark}
   Theorem \ref{hw_bogo} is of course well known \cite{MR893593}, \cite{MR1399482}
   in the case $n=0$ (ie: no
   singularities). In this case, the proof we have offered holds if we take
   $R = 0$. Thus $h$ is an isomorphism on all $\Q$ and so we really obtain
   a holomorphic vector bundle on $\Q$ rather than a triple $(E^+, E^-, h)$
   as in the singular case. Furthermore, it is clear that the proof could
   be generalised to groups other than $SU(2)$. In particular we may consider
   $U(1)$ monopoles with no singularities on $\Hyp$. In this
   case the trivial unit mass $U(1)$ monopole on $\Hyp$ yields a holomorphic
   line bundle L over $\Q$. Note that $L$ will carry a canonical anti-holomorphic
   (anti-linear) bijection:
   \begin{equation}
      L \to L^{*}
   \end{equation}
   covering $\sigma : \Q \to \Q$ since we started with a $U(1)$ bundle. Now it
   is well known \cite{MR893593}, \cite{MR1399482} that $L\simeq \cO(1, -1)$
   however it is worth mentioning that we may recover this result with minimal
   effort given the approach we have taken.
\end{remark}
\begin{proposition}\label{def_L}
Let
   \begin{eqnarray*}
      \pi : L \to \Q
   \end{eqnarray*}
   be the holomorphic line bundle corresponding to the trivial unit mass
   $U(1)$ monopole on $\Hyp$. Then
   \begin{eqnarray*}
      L \simeq \cO(1,-1)
   \end{eqnarray*}
\end{proposition}
\noindent{\bf Proof}
   From the recipe of theorem \ref{hw_bogo}, the form defining the
   $\deebar$-operator of $L$ is $\theta^{0,1}$. Using the formula of equation
   \eqref{theta_formula} we can now read off that $L\simeq\cO(1,-1)$ as
   required. \myqed
\begin{remark}
   Consider again the 1-form $\theta$ of definition \ref{theta01_def}.
   Since $\deebar\theta^{0,1} = 0$ we have a cohomology class
   \begin{eqnarray*}
      \lbrack \theta^{0,1}\rbrack \in H_{\deebar}^{0,1}(\Q)
      \simeq H^1(\Q,\cO)
   \end{eqnarray*}
   If $exp : H^1(\Q, \cO) \to H^1(\Q, \cO^*)$ is the usual exponential
   map, then $\theta^{0,1}$ defines the line bundle
   $exp(\lbrack\theta^{0,1}\rbrack)$ and this is if course our line
   bundle $L$.
   \par
   Since $L$ is in the image of the map $exp$, it is trivial as a $C^{\infty}$
   complex line bundle. This is the reason that we can (and will, see theorem
   \ref{hw_E_ext}) raise it to non-integral powers. Indeed $H^1(\Q,\cO)$ is a
   complex vector space and so $L$ can be raised to any complex power.
   \par
   We also note that the above method for finding the line bundle $L$ corresponding
   to the trivial $U(1)$ monopole can, of course, also be applied in the Euclidean
   case. Using the formula (\ref{Euc_theta01_def}) we can thus recover Hitchin's
   line bundle $L$, cf \cite{MR649818}.
\end{remark}
\par
Now that we have the Hitchin--Ward correspondence for solutions of the
Bogomolny equations over $\Hyp\setminus \{p_1, \ldots, p_n\}$, the next
step is to work out what consequences the boundary conditions in definition
\ref{def_sngmpole} have for the holomorphic data $(E^+, E^-, h)$. We first
deal with conditions (iii) of definition \ref{def_sngmpole}. Clearly the behaviour
of $\Phi$ near the singularities $p_i$ will be reflected in the behaviour of
$h$ near $P = P_1\cup\cdots\cup P_n$, where $P_i = P_{p_i}$ (the twistor line
that is the set of geodesics passing through $p_i$). Indeed
let $\tilde p_i \in H^0(\Q, \cO(1,1))$ be a section with divisor $P_i$ and let
\begin{eqnarray*}
   \tilde p = \prod_{i=1}^{n} \tilde p_i^{l_i} \in H^0(\Q, \cO(l, l))
\end{eqnarray*}
for some $\{l_1, \ldots ,l_n\} \subset \N$ and $l = \sum l_i$.
Note that we can regard
\begin{eqnarray*}
   h \in H^0(\Q\setminus P, \Hom(E^+, E^-))
\end{eqnarray*}
and that the singular set of $h$ is exactly the same as the zero set of $\tilde p$.
With this notation in place, we can state
\begin{theorem}\label{hw_h_sings}
   Let $(E^+, E^-, h)$ be the holomorphic data corresponding to a
   solution to the Bogomolny equations as in theorem \ref{hw_bogo}.
   Let $I \subset \Q$ be the set of geodesics in $\Hyp$ that pass through
   at least two of the singular points $p_i \in \Hyp$. Then the solution
   of the Bogomolny equations satisfies conditions (iii) of definition
   \ref{def_sngmpole} and has Abelian charges $l_1, \ldots ,l_n$
   if and only if the section
   \begin{eqnarray*}
      \tilde ph \in H^0(\Q\setminus P, \Hom(E^+, E^-)(l, l))
   \end{eqnarray*}
   extends across $P$ to a holomorphic section on all of $\Q$ and is non-vanishing
   on $\Q \setminus I$.
\end{theorem}
\noindent {\bf Proof}
   We show that a solution of the Bogomolny equations satisfying conditions (iii)
   of definition \ref{def_sngmpole} has the required property and omit the
   proof of the converse since we do not require it.
   \par
   Now since $I$ is a discrete set of points and $\Q$ has complex dimension 2,
   if $\tilde ph$ extends as a holomorphic section to $\Q\setminus I$ then, by
   Hartog's theorem, the isolated singularities at the points of $I$ are removable.
   \par
   Let $P_i \subset \Q$ be the twistor line in $\Q$ corresponding to a singularity
   $p_i \in \Hyp$ and let $x \in P_i \setminus I$. Let $U_x \subset \Q\setminus
   \bigcup_{j\ne i}\limits P_j$ be an open neighbourhood of $x$ in $\Q\setminus
   \bigcup_{j\ne i}\limits P_j$.
   \par
   Consider $\tilde p_j$, $1 \le j \ne i \le n$. This is non-vanishing on $U_x$
   and so $\tilde ph$ has a removable singularity along $U_x \cap P_i$ iff $\tilde
   p_i^{l_i}h$ does and $\tilde ph$ is non-vanishing on $U_x \cap P_i$ iff $\tilde
   p_i^{l_i}h$ is. This means that we can deal with each singularity separately. 
   We just need to prove the result for a single singularity. Identifying $\Q
   \simeq \QQ$, we may take this
   singularity to be at the point $O$ whose twistor line is the diagonal $\Delta
   \subset \QQ$. We let $l$ be the Abelian charge.
   \par
   Now, there is a holomorphic trivialisation of $\cO(1,1)$ over the open set $V$ of
   $\QQ$ with coordinates $([z,1],[w,1])$ such that $\tilde O \in
   H^0(\Q, \cO(1,1))$ is trivialised as the function
   \begin{eqnarray*}
      (z,w) \mapsto z-w
   \end{eqnarray*}
   Thus if $h$ has matrix
   \begin{eqnarray*}
      \left[
      \begin{array}{cc}
         a_1 & a_2\\
         a_3 & a_4
      \end{array}\right]
   \end{eqnarray*}
   relative to local trivialisations of $E^{\pm}$, then we need to
   investigate the behaviour of the functions
   \begin{eqnarray*}
      (z,w) \mapsto (z-w)^l a_i(z,w)
   \end{eqnarray*}
   as $z \to w$. It is sufficient to do this for each fixed value $w=w_0$. 
   Furthermore, we may without loss of generality assume $w_0 = 0$ since we
   may always choose our coordinates to arrange for this. From
   now on we thus work on a fixed slice $w=0$.
   \par
   Now it will follow from our work below that the functions $a_i$ (defined on
   a neighbourhood of $0$ in $\C^*$) cannot have essential singularities at $z = 0$.
   At worst they have poles. Thus there exist unique integers $m,n\in\Z$ such that
   \begin{eqnarray*}
      z^n\left[
      \begin{array}{l}
         a_4(z)\\
         -a_3(z)
      \end{array}\right]
   \end{eqnarray*}
   and
   \begin{eqnarray*}
      z^m\left[
      \begin{array}{l}
         -a_2(z)\\
         a_1(z)
      \end{array}\right]
   \end{eqnarray*}
   have removable singularities at $0$ are are non-vanishing in a neighbourhood of
   $0$. Using these to define a new trivialisation of $E^+$ we find that $h$ has
   matrix
   \begin{eqnarray}\label{h_matrix}
      \left[
      \begin{array}{cc}
         z^n & 0\\
         0   & z^m
      \end{array}\right]
   \end{eqnarray}
   Since $h$ is compatible with the symplectic structures on $E^+$ and $E^-$, it must
   have a regular determinant and so we must have $n=-m$. Without loss of generality,
   we may assume $n > 0$. Evidently we will be done if we can show that $n=l$ and,
   in view of \eqref{h_matrix}, to do this it is enough to show that there exist
   $C_1, C_2 > 0$ and $1 > \epsilon_1, \epsilon_2 > 0$ such that for all non-zero $z$
   with $|z|$ small enough we have 
   \begin{eqnarray*}
      C_1|z|^{\epsilon_1} < |z|^l\|H(z)\| < C_2|z|^{-\epsilon_2}
   \end{eqnarray*}
   where $H$ is the matrix of $h$ with respect
   to local trivialisations of $E^\pm$ and we are using the $l^1$ norm on matrices.
   \par
   We will show this by choosing appropriate local trivialisations as detailed below.
   Thus for any $\delta > 0$ define
   \begin{eqnarray*}
      \Gamma(\delta) = \{([z,1],[0,1]) \in \Q \quad |\quad 0 < |z| < \delta\} \subset
      \Q\setminus \Delta
   \end{eqnarray*}
   and
   \begin{eqnarray*}
      \Omega(\delta) = \{\gamma(t)\in\Hyp\quad |\quad \gamma \in \Gamma(\delta)
      \mbox{ and } |t| < \delta\} \subset \Hyp\setminus\{O\}
   \end{eqnarray*}
   where we have given $\gamma$ the parameterisation determined by $O \in \Hyp$.
   If $\delta$ is small enough, then $\Phi$ is non-zero on $\Omega(\delta)$. Thus,
   for $j=0,1$, let $e_j$ be a unitary eigensection of $E|_{\Omega(\delta)}$ with
   eigenvalue $(-1)^ji\|\Phi\|$. Also let
   \begin{eqnarray*}
      \left[\begin{array}{cc}
         A    & B\\
         -B^* & -A
      \end{array}\right]
   \end{eqnarray*}
   be the matrix of the monopole connection $\nabla$ with respect to the local
   trivialisation defined by $e_0, e_1$. Note that $B$ is bounded in a neighbourhood
   of $O$ since
   \begin{eqnarray*}
      \left\| \nabla\left(\frac{\Phi}{\|\Phi\|}\right)\right\| = 2\|B\|
   \end{eqnarray*}
   and it follows easily from conditions (iii) of definition \ref{def_sngmpole}
   that $\nabla(\Phi/\|\Phi\|)$ is bounded in a neighbourhood of $O$.
   \par
   Now if $s$ is a section of $E$ over $\{\gamma(t)\quad|\quad|t| < \delta\}$
   for some $\gamma \in \Gamma(\delta)$ then in terms of the trivialisation of $E$
   determined by $e_0,e_1$ the equation $(\nabla_{\dot\gamma} - i\Phi)s = 0$ becomes
   \begin{eqnarray}\label{scatt_eqn_removsingprf}
      \frac{ds}{dt} = \left[
      \begin{array}{cc}
         \|\Phi\| + A(\dot\gamma) & -B^*(\dot\gamma)\\
         B(\dot\gamma)            & -\|\Phi\| - A(\dot\gamma)
      \end{array}
      \right]s
   \end{eqnarray}
   Let $H(z, t)$ be the matrix solution of this equation such that
   $H(z, -\delta) = I$. Then $H(z) = H(z, \delta)$ is the matrix of $h$ in the
   local trivialisations of $E^\pm$ determined by $e_0, e_1$.
   \par
   To proceed with the required analysis of \eqref{scatt_eqn_removsingprf} define
   \begin{eqnarray*}
      G(z,t) = \exp\left(-\int_{-\delta}^t (\|\Phi\| + A(\dot\gamma))ds\right)H(z,t)
   \end{eqnarray*}
   Then $H$ solves \eqref{scatt_eqn_removsingprf} iff $G$ solves
   \begin{eqnarray}\label{scatt_eqn2_removsingprf}
      \frac{dG}{dt} = \left[
      \begin{array}{cc}
         0             & -B^*(\dot\gamma)\\
         B(\dot\gamma) & -2(\|\Phi\| + A(\dot\gamma))
      \end{array}
      \right]G
   \end{eqnarray}
   We estimate the behaviour of solutions of this by regarding the off diagonal terms
   (for which we have a bound) as a perturbation of the diagonal terms and transforming
   to an integral equation. Thus let $U$ be the fundamental solution of the diagonal
   equation, ie:
   \begin{eqnarray*}
      U(t) = \left[
      \begin{array}{cc}
         1 & 0\\
         0 & \exp\left(-2\int_{-\delta}^t (\|\Phi\| + A(\dot\gamma))ds\right)
      \end{array}
      \right]
   \end{eqnarray*}
   Note that $1 \le \|U(t)\| \le 2$. Define the integral operator $T$ by
   \begin{eqnarray*}
      (TG)(t) = U(t)\int_{-\delta}^t U(s)^{-1}\left[
      \begin{array}{cc}
         0             & -B^*(\dot\gamma)\\
         B(\dot\gamma) & 0
      \end{array}
      \right]G(s)ds
   \end{eqnarray*}
   Since $B$ is bounded, let $M$ be a constant such that $\|B\| < M$ on $\Omega(\delta)$.
   Then for $t \in [-\delta, \delta]$ we have
   \begin{eqnarray*}
      \|(TG)(t)\| &\le& \sup_{[-\delta,\delta]}\|G\|M\int_{-\delta}^t\left\|\left[
      \begin{array}{cc}
         1 & 0\\
         0 & \exp\left(-2\int_{s}^t (\|\Phi\| + A(\dot\gamma))du\right)
      \end{array}
      \right]\right\|ds\\
      &\le& 4M\delta \sup_{[-\delta,\delta]}\|G\|
   \end{eqnarray*}
   Thus by taking $\delta > 0$ small enough, we can have $\|T\| < \epsilon$ for any
   $\epsilon > 0$, where we are using the sup-norm for $T$. Now $G$ solves
   \eqref{scatt_eqn2_removsingprf} iff it solves
   \begin{eqnarray*}
      G = U + TG
   \end{eqnarray*}
   and so using the sup-norm for everything, we have
   \begin{eqnarray*}
      \|G - U\| &=& \|(T + T^2 + \cdots)U\|\\
               &\le& (\|T\| + \|T\|^2 + \cdots)\|U\|\\
                &=& \frac{2}{\|T\|^{-1} - 1}\|U\|
   \end{eqnarray*}
   Thus provided we take $\delta > 0$ small enough we have
   \begin{eqnarray*}
      \|G - U\| < \epsilon
   \end{eqnarray*}
   for any $\epsilon > 0$ and thus
   \begin{eqnarray*}
      \|U\| - \epsilon < &\|G\|& < \|U\| + \epsilon\\
      \Rightarrow 1 - \epsilon < &\|G\|& < 2 + \epsilon
   \end{eqnarray*}
   and so finally
   \begin{eqnarray}\label{blah_feck}
      (1 - \epsilon)\exp\left(\int_{-\delta}^\delta \|\Phi\|dt\right)
      < \|H(z)\| < (2 + \epsilon)\exp\left(\int_{-\delta}^\delta \|\Phi\|dt\right)
   \end{eqnarray}
   \par
   It only remains to deal with the behaviour of $\exp\left(\int_{-\delta}^\delta 
   \|\Phi\|dt\right)$ as $z \to 0$. To do this, note that from conditions (iii) of
   definition \ref{def_sngmpole} for $\delta > 0$ sufficiently small and
   $0 < |z| < \delta$, $|t| < \delta$ we have
   \begin{eqnarray*}
      l - \frac{1}{2} < 2\rho(\gamma(t), O)\|\Phi\| < l + \frac{1}{2}
   \end{eqnarray*}
   Next by Pythagoras's theorem for hyperbolic space (applied to the triangle
   with vertices $O, \gamma(t), \gamma(0)$ which has a right angle at $\gamma(0)$)
   we have
   \begin{eqnarray*}
      \cosh\rho(\gamma(t), O) = \rho(\gamma(0), O)\cosh t
   \end{eqnarray*}
   and furthermore\footnote{This a special case of the general formula $\cosh(\rho)
   = \frac{\sqrt{(1+|z|^2)(1+|w|^2)}}{|1+z\overline w|}$.}
   \begin{eqnarray*}
      \cosh\rho(\gamma(0), O) = \sqrt{1 + |z|^2}
   \end{eqnarray*}
   so that we have
   \begin{eqnarray*}
      l - \frac{1}{2} < 2\cosh^{-1}\left(\sqrt{1 + |z|^2}\cosh t\right)\|\Phi\|
      < l + \frac{1}{2}
   \end{eqnarray*}
   and so
   \begin{eqnarray*}
      \frac{l - 1/2}{2\sqrt{t^2 + |z|^2}} <
      \frac{\cosh^{-1}\left(\sqrt{1 + |z|^2}\cosh t\right)}{\sqrt{t^2 + |z|^2}}\|\Phi\|
      < \frac{l + 1/2}{2\sqrt{t^2 + |z|^2}}
   \end{eqnarray*}
   Now since the function
   \begin{eqnarray*}
      f(s, t) = \frac{\cosh^{-1}\left(\sqrt{1 + s^2} \cosh t\right)}{\sqrt{t^2 + s^2}}
   \end{eqnarray*}
   which is a priori defined on $\R^2\setminus \{0\}$ in fact extends as a continuous
   function on $\R^2$ with value $1$ at $0$ we have
   \begin{eqnarray*}
      \frac{1}{1 + \epsilon'}\frac{l - 1/2}{2\sqrt{t^2 + |z|^2}} < \|\Phi\| <
      \frac{1}{1 - \epsilon'}\frac{l + 1/2}{2\sqrt{t^2 + |z|^2}}
   \end{eqnarray*}
   for any $\epsilon' > 0$, provided $\delta > 0$ is sufficiently small and $|t|,
   |z| < \delta$. Integrating, we thus have
   \begin{eqnarray*}
      \frac{l - 1/2}{1 + \epsilon'}\sinh^{-1}(\delta/|z|) < \int_{-\delta}^\delta
      \|\Phi\|dt < \frac{l + 1/2}{1 - \epsilon'}\sinh^{-1}(\delta/|z|)
   \end{eqnarray*}
   Taking exponentials and remembering that $\exp(\sinh^{-1}x) = x + \sqrt{x^2 + 1}
   \sim 2x$ as $x \to \infty$, we thus have
   \begin{eqnarray*}
      C'_1 |z|^{-l + \epsilon_1} < \exp\left(\int_{-\delta}^\delta \|\Phi\| dt\right)
      < C'_2 |z|^{-l - \epsilon_2}
   \end{eqnarray*}
   for constants $C'_1, C'_2 > 0$ and $1 > \epsilon_1, \epsilon_2 > 0$. Combining this
   with \eqref{blah_feck} we thus have
   \begin{eqnarray*}
      C_1 |z|^{-l + \epsilon_1} < \|H(z)\| < C_2 |z|^{-l + \epsilon_2}
   \end{eqnarray*}
   for appropriate constants, as required. \myqed
Note that the above theorem also applies to $h^{-1}$. From now on, when we refer
to $\tilde ph$ and $\tilde ph^{-1}$ it shall be understood that they are regarded as
being defined over all of $\Q$.
\par
Having dealt with conditions (iii) of definition \ref{def_sngmpole} we need only note
that conditions (ii) have the same effect on $E^{\pm}$ as in the non-singular
case. Indeed, suppose that $(E^+, E^-, h)$ is the holomorphic data corresponding to
a solution of the Bogomolny equations as in theorem \ref{hw_bogo}. Suppose also that
this solution satisfies conditions (ii) of definition \ref{def_sngmpole}. Define
$L^+ \subset E^+$ as follows. Let
\begin{eqnarray*}
   v \in E^+_{\gamma} \simeq E_{\gamma(R)}
\end{eqnarray*}
and let $s$ be the unique section of $E$ over $\left\{\gamma(t)\quad|\quad t\ge R
\right\}$ such that $s(\gamma(R)) = v$ and $(\nabla_{\dot\gamma} - i\Phi)s = 0$. We
define $v \in L^+$ iff
\begin{eqnarray*}
   s(\gamma(t)) \to 0 \mbox{ as } t \to \infty
\end{eqnarray*}
We then have
\begin{theorem}\label{hw_E_ext}
   In the above notation, $L^+$ is a holomorphic line sub-bundle of $E^+$. Furthermore
   \begin{eqnarray*}
      L^+ \simeq L^m(0, -N)
   \end{eqnarray*}
   where $L$ is the line bundle defined in proposition \ref{def_L} and $N$ is the
   total charge of the monopole. Also, since $E^+$ has a holomorphic symplectic
   structure, we have
   \begin{eqnarray*}
      E^+ / L^+ \simeq (L^+)^{*}
   \end{eqnarray*}
   and so we can express $E^+$ as an extension:
   \begin{eqnarray}\label{LE_ext}
      0 \to L^+ \to E^+ \to (L^+)^{*} \to 0
   \end{eqnarray}
\end{theorem}
\noindent{\bf Proof}
   This can be proved by modifying a proof in the non-singular case.
   Since it is only the asymptotic behaviour of the sections of $E$
   which matters, the singularities do not cause any complication (note
   that $h$ does not even enter the statement of the theorem). For a proof
   of the result in the non-singular case see \cite{MR1399482}. \myqed
Recall now that we have the map $\tilde j : E^+ \to E^-$.
We thus have a bundle
\begin{eqnarray*}
   L^- = \tilde j L^+ \subset E^-
\end{eqnarray*}
Using the facts that $\nabla j = 0$, $j$ is anti-linear and that $j$
covers $\sigma$, it is easy to identify $L^-$ in the same way as $L^+$. Let
\begin{eqnarray*}
   v \in E^-_{\gamma} \simeq E_{\gamma(-R)}
\end{eqnarray*}
and let $s$ be the unique section of $E$ over $\{\gamma(t)\quad|\quad
t\le -R\}$ such that $s(\gamma(-R)) = v$ and $(\nabla_{\dot\gamma} - i\Phi)s = 0$.
Then $v \in L^-$ iff
\begin{eqnarray*}
   s(\gamma(t)) \to 0 \mbox{ as } t \to -\infty
\end{eqnarray*}
We thus have a corresponding expression of the bundle $E^-$ as an extension:
\begin{eqnarray}\label{LE-_ext}
   0 \to L^- \to E^- \to (L^-)^* \to 0
\end{eqnarray}
Furthermore a choice of isomorphism $L^+ \simeq L^m(0,-N)$ induces an isomorphism
$L^- \simeq L^{-m}(-N,0)$.
\par
In the case of non-singular monopoles,
the situation is a little simpler since $h$ is a global isomorphism
and so we can work with $h^{-1}L^- \subset E^+$.
\section{The spectral data}
We are now in a position to combine the results of the previous section
and identify the spectral data which determine a singular hyperbolic $SU(2)$
monopole.
\begin{definition}\label{spec_curve}
   Let $(E^+, E^-, h)$ be the holomorphic data corresponding to a singular
   hyperbolic $SU(2)$ monopole as in theorem \ref{hw_bogo}. Define the map
   $\psi$ as follows:
   \begin{eqnarray}\label{psi_def}
      \psi : L^+ \to E^+ \to E^-(l, l) \to (L^-)^{*}(l, l)
   \end{eqnarray}
   where the second arrow is the map $\tilde ph$ of theorem \ref{hw_h_sings}
   and the last arrow is formed by tensoring the projection of the exact sequence
   \eqref{LE-_ext} with the identity map on $\cO(l, l)$. Note that using
   an isomorphism $L^+ \simeq L^m(0, -N)$ as in theorem \ref{hw_E_ext}
   we can regard
   \begin{eqnarray*}
      \psi \in H^0(\Q, \cO(k,k))
   \end{eqnarray*}
   where $k$ is the non-Abelian charge of the monopole. We define the
   spectral curve $S$ of the monopole to be the divisor of $\psi$.
\end{definition}
\begin{remark}
   In the case of non-singular monopoles, the geometric interpretation of the
   spectral curve is clear. Since $S$ is a subset of twistor space which
   is the set of all oriented geodesics in $\Hyp$ and $S$ is preserved by
   $\sigma$, $S$ really defines a set of  unoriented lines in $\Hyp$.
   These are known as the spectral lines of the monopole. Chasing
   through the definitions one finds that a line $\gamma$ in $\Hyp$ is a
   spectral line iff there exists a non-zero section $s$ of $E$ along $\gamma$
   such that $(\nabla_{\dot\gamma} - i\phi)s = 0$ and $s(\gamma(t)) \to 0$ as $t
   \to \pm\infty$.
   \par
   For our singular monopoles, the situation is more complicated. The
   spectral curve $S$ of the singular monopole still defines a set of
   spectral lines in $\Hyp$ but it is not as easy to identify them
   geometrically. For a line $\gamma$ which does not pass through a singular
   point, the rule for deciding if $\gamma$ is a spectral line is the same as
   for a non-singular monopole. However if $\gamma \in S \cap P$ this no longer
   makes sense. Indeed if we return to the definition of $S$ using $\tilde
   ph$ we see that to define $S$ we had to note that $\tilde ph$ had a
   removable singularity along $P$. Obtaining the value of a holomorphic
   function at a removable singularity requires a limiting process and this
   means that to decide if a line $\gamma$ passing through a singular point
   is a spectral line we will need to examine the behaviour of sections of
   $E$ along geodesics in a neighbourhood of $\gamma$ in $\Hyp$.
\end{remark}
For non-singular (Euclidean or hyperbolic) monopoles, the spectral
curve determines the monopole. As we shall see, this is almost true for singular
monopoles. Except for the special case (which we shall not consider) when
$S \cap P$ is not finite (ie when $P_i$ is a connected component of $S$ for some
$i$), exactly one additional piece of spectral data is needed to identify a
singular monopole.
\par
Let $S$ be the spectral curve of a singular hyperbolic $SU(2)$ monopole for
which $S \cap P$ is finite. Note that by definition of $S$ and the exactness of
the sequence
\begin{eqnarray*}
   0 \to L^-(l, l) \to E^-(l, l) \to (L^-)^{*}(l, l) \to 0
\end{eqnarray*}
it follows that the image of $\tilde ph$ restricted to $L^+|_{S}$ is in fact
contained in $L^-(l, l)|_{S}$. We thus have a map $\xi^- : L^+|_{S} \to
L^-(l, l)|_{S}$. Note that we can regard
\begin{eqnarray}\label{xi_def}
   \xi^- \in H^0(S, (L^+)^{*}L^-(l, l))
\end{eqnarray}
Similarly, considering $\tilde ph^{-1}$ restricted to $L^-|_{S}$ we have
\begin{eqnarray}\label{xi_bar_def}
   \xi^+ \in H^0(S, (L^-)^{*}L^+(l, l))
\end{eqnarray}
Clearly these two satisfy $\xi^-\xi^+ = \tilde p^2|_{S}$. Let $D$ be the divisor of
$\xi^+$. Using the fact that $(h^{-1}\tilde j)^2 = -1$ outside $P$ we find that
the total map:
\begin{eqnarray*}
   L^+|_{S} \overset{\tilde j}\to L^-|_{S} \overset{\xi^+}\to
   L^+(l, l)|_{S} \overset{\tilde j}\to L^-(l, l)|_{S}
\end{eqnarray*}
is just the map:
\begin{eqnarray*}
-\xi^- : L^+|_{S} \to L^-(l, l)|_{S}
\end{eqnarray*}
It thus follows (since $\tilde j$ covers $\sigma$) that the divisor of $\xi^-$ is
$\sigma(D)$.
\begin{definition}\label{spec_div_def}
   Using the above notation, we call $D$ the spectral divisor of the monopole.
   Also since $\sigma(S \cap P) = S \cap P$ we see that $S \cap P$ really
   defines a set of unoriented lines in $\Hyp$. We call these lines the
   singular spectral lines of the monopole. We shall refer to $(S, D)$ as
   the spectral data of a monopole.
\end{definition}
Note that the support $|D|$ of $D$ is contained in $S \cap P$
since $\tilde ph^{-1}$ is an isomorphism outside $P$. Also note that since
$\xi^-\xi^+ = \tilde p^2|_{S}$ we have
\begin{eqnarray*}
   |D| \cup \sigma(|D|) = S \cap P
\end{eqnarray*}
Since $\sigma$ has no fixed points
\begin{eqnarray*}
   |D| \cap \sigma(|D|) = \emptyset
\end{eqnarray*}
and so $|D|$ defines a partition of the set of singular spectral lines
into disjoint conjugate subsets. Recalling that $\sigma(\gamma)$ is just
$\gamma$ parameterised in the opposite direction this means that $|D|$ really
defines an orientation for each singular spectral line. Finally note that since
\begin{eqnarray*}
   D + \sigma(D) = (\tilde p^2|_{S})
\end{eqnarray*}
it follows that that $|D|$ determines $D$ (provided we know $\tilde p^2$,
ie: the locations of the singularities).
As we shall see the spectral data $(S, D)$ determines the monopole and so
the spectral data for a singular hyperbolic $SU(2)$ monopole for which
$S \cap P$ is finite may be regarded as the set of spectral lines in $\Hyp$
together with an orientation for each singular spectral line.
\par
We gather together a few important properties of the spectral data for
a singular monopole.
\begin{proposition}\label{sing_hyp_spec_crv_props}
   Let $(S, D)$ be the spectral data of a singular hyperbolic
   $SU(2)$ monopole of non-Abelian charge $k$ for which $S \cap P$ is finite. Then
   \begin{enumerate}
      \item
      $S$ is compact.
      \item
      $S$ is real (ie: preserved by $\sigma$).
      \item
      $L^{2m+k}(0, 2l)|_{S} \simeq \lbrack D\rbrack$
      \item
      If $S$ is non-singular then it has genus $(k-1)^2$.
   \end{enumerate}
\end{proposition}
\noindent{\bf Proof}
   \begin{enumerate}
      \item
      Let $\gamma$ be a geodesic in $\Hyp$ that does not pass through any
      of the singular points, $p_i$. We have already noted that the condition
      for $\gamma$ to be a spectral line is that there exists a non-zero
      solution $s$ to $(\nabla_{\dot\gamma} - i\Phi)s = 0$ along $\gamma$ such
      that $s(t) \to 0$ as $t \to \pm\infty$. Thus if we fix a point $O \in \Hyp$,
      the argument used in \cite{MR649818} to prove compactness of the spectral
      curve of a non-singular Euclidean monopole shows that there exists $M > 0$
      such that if $\gamma$ is a spectral line then it meets the ball $B(O, M)
      \subset \Hyp$ of radius $M$ centred about $O$. Thus if we choose $M$
      sufficiently large that we also have $\{p_1, \ldots, p_n\} \subset B(O, M)$
      we find that all spectral lines (ie: including the singular spectral lines
      passing through the points $p_i$) meet $B(O, M)$.
      \par
      Now consider the ball model of $\Hyp$ and use this to identify
      $\Q \subset \CP\times\CP$. Give the conformal boundary of the ball the usual
      round metric of $\CP\simeq S^2$ and let $d$ be the geodesic distance
      function. In view of the above, there exists $N_M >0$ such that
      \begin{eqnarray*}
         S \subset \left\{(p,q)\in \CP\times\CP \quad |\quad d(p,q) \le N_M\right\}
      \end{eqnarray*}
      and so $S$ being a closed subset of a compact space is compact.
      \item
      This follows since $\tilde j : L^+ \to L^-$ covers $\sigma$.
      \item
      An isomorphism $L^+ \simeq L^m(0, -N)$ as in theorem \ref{hw_E_ext} induces
      an isomorphism
      \begin{eqnarray*}
         (L^-)^*L^+(l,l) \simeq L^{2m+k}(0,2l)
      \end{eqnarray*}
      Since the spectral divisor $D$ is the divisor of a section of $(L^-)^*
      L^+(l,l)|_S$ the result follows.
      \item
      This follows from the adjunction formula since $S$ is a divisor of $\cO(k,k)$.
   \end{enumerate}
   \myqed
The most important of the properties listed in the above proposition is
(iii) since it is the only property which differs from the non-singular
case. It replaces the condition that $L^{2m+k}|_S$ must be trivial which is what
holds in the non-singular case. As we shall see in the next section, it is exactly
condition (iii) that means that the spectral curves of $k=1$ singular hyperbolic
monopoles lift to twistor lines in appropriate twistor spaces.
\par
We wish to show that the spectral data determine the monopole. To do this we
shall need the following lemma.
\begin{lemma}\label{lemma_ext_cobound}
   Let $S$ be the spectral curve of a singular hyperbolic $SU(2)$ monopole
   for which $S \cap P$ is finite and let $\xi^+$ be the section of
   \eqref{xi_bar_def}. Let $a \in H^1(\Q, (L^+)^2)$ be the class
   representing the extension \eqref{LE_ext} and let
   \begin{eqnarray}\label{cobound}
      \delta : H^0(S, (L^-)^{*} L^+(l, l)) \to H^1(\Q, (L^+)^2)
   \end{eqnarray}
   be the connecting homomorphism associated to the short exact sequence
   of sheaves
   \begin{eqnarray}\label{SES_ext_cobound}
     0 \to \cO_{\Q}((L^+)^2) \to \cO_{\Q}((L^-)^{*}L^+(l, l)) \to \cO_S((L^-)^{*}
     L^+(l, l)) \to 0
   \end{eqnarray}
   Then $\delta\xi^+ = a$.
\end{lemma}
\noindent {\bf Proof} By abuse of notation let $a \in \Omega^{0,1}(\Q, (L^+)^2)$
be a Dolbeault representative for the extension class of \eqref{LE_ext}. By
exactness of the long exact sequence of cohomology groups associated to
\eqref{SES_ext_cobound} $[a]$ is in the image of $\delta$ if and only if
\begin{eqnarray}\label{img_delta_cond}
   \psi a = \deebar b
\end{eqnarray}
for some $b \in \Omega^0(\Q, (L^-)^*L^+(l,l))$ where $\psi$ is the section
of equation \eqref{psi_def}. Furthermore, in this case
\begin{eqnarray*}
   b|_S \in H^0(S, (L^-)^*L^+(l,l))
\end{eqnarray*}
is the class mapped to $[a]$ under $\delta$. Now if we fix a smooth splitting
$E^+ = L^+ \oplus (L^+)^*$ of \eqref{LE_ext} on $\Q$ then the $\deebar$-operator
of $E^+$ is
\begin{eqnarray*}
   \deebar = \left[
   \begin{array}{cc}
      \deebar_{L^+} & a\\
      0             & \deebar_{(L^+)^*}
   \end{array}\right]
\end{eqnarray*}
and from this we can see that \eqref{img_delta_cond} holds if and only if
we have a meromorphic splitting of \eqref{LE_ext} with a pole along $S$.
To see this, note that in our fixed smooth splitting, to define
the meromorphic splitting we only need to define the map $(L^+)^* \to L^+$
over $\Q\setminus S$ and this is
\begin{eqnarray*}
   b/\psi : (L^+)^* \to L^+
\end{eqnarray*}
To prove the lemma, we thus need a meromorphic
splitting $\alpha : (L^+)^* \to E^+$ of \eqref{LE_ext} with a pole along $S$ such
that if $\alpha' : (L^+)^* \to L^+$ is the induced map (using our fixed smooth
splitting) then $(\psi\alpha')|_S = \xi^+$. However by definition of $\psi$
this is the same as requiring that $\xi^+$ is given by the map
\begin{eqnarray}\label{xi_split_condition}
   (L^+)^* \overset{\alpha}\to E^+ \overset{\tilde ph}\to E^-(l,l) \to
   (L^-)^*(l,l)
\end{eqnarray}
(restricted to $S$). We thus need only exhibit a meromorphic splitting of
\eqref{LE_ext} that satisfies \eqref{xi_split_condition}. We thus define
$\alpha$ to be
\begin{eqnarray*}
   (L^+)^* \overset{\psi^{-1}}\to L^-(-l,-l) \overset{\tilde ph^{-1}}\to E^+
\end{eqnarray*}
It is trivial to check that $\alpha$ is indeed a splitting of \eqref{LE_ext} and
satisfies condition \eqref{xi_split_condition}. \myqed
\begin{theorem}\label{spec_data_inv}
Let $(S, D)$ be the spectral data of a singular hyperbolic $SU(2)$ monopole
of mass $m$ and non-Abelian charge $k$ with Abelian charge $l_i$ at the singular
points $p_i$, $i=1,\ldots,n$ for which $S \cap P$ is finite.
This data determines the monopole.
\end{theorem}
\noindent{\bf Proof} We must show that we can recover the data $(E^+, E^-, h)$
   of theorem \ref{hw_bogo} from the spectral data $(S, D)$.
   \par
   Thus let $N = k - l$. Let $s$ be a section of $L^{2m+k}(0, 2l)$ on $S$ with
   divisor $D$ and such that $ss^* = 1$. Let $a' = \delta s \in H^1(\Q, L^{2m}
   (0, -2N))$ where
   \begin{eqnarray*}
      \delta : H^0(S, L^{2m+k}(0, 2l)) \to H^1(\Q, L^{2m}(0, -2N))
   \end{eqnarray*}
   is the connecting homomorphism associated to the short exact sequence of
   sheaves
   \begin{eqnarray*}
      0 \to \cO_{\Q}(L^{2m}(0, -2N)) \to \cO_{\Q}(L^{2m+k}(0, 2l)) \to
      \cO_{S}(L^{2m+k}(0, 2l)) \to 0
   \end{eqnarray*}
   Let ${E'}^+$ be the bundle defined (up to equivalence) as an extension of
   $L^m(0, -N)$ by $L^{-m}(0, N)$ using $a'$ as extension class.
   \par 
   Choose the unique isomorphism $L^+ \simeq L^m(0,-N)$ such that the
   image of $s$ under the induced isomorphism
   \begin{eqnarray*}
      H^0(S, L^{2m+k}(0, 2l)) \simeq H^0(S, (L^+)^{*}L^-(l, l))
   \end{eqnarray*}
   is $\xi^+$ of \eqref{xi_bar_def}. Then by lemma \ref{lemma_ext_cobound}
   the image of $a'$ under the induced isomorphism
   \begin{eqnarray*}
      H^1(\Q, L^{2m}(0, -2N)) \simeq H^1(\Q, (L^+)^2)
   \end{eqnarray*}
   is the extension class of \eqref{LE_ext}. It follows that
   we must have ${E'}^+ \simeq E^+$. Similarly, we can recover $E^-$.
   \par
   It remains only to show that $h$ is determined by $(S, D)$. Now we saw
   in the course of the proof of lemma \ref{lemma_ext_cobound} that the
   sequence \eqref{LE_ext} has a natural splitting on $\Q\setminus S$,
   ie: on $\Q\setminus S$, we naturally have
   \begin{eqnarray*}
      E^+ \simeq L^+ \oplus (L^+)^*
   \end{eqnarray*}
   and similarly for $E^-$. A quick review of the relevant definitions reveals
   that the matrix of $h$ with respect to these splittings is (cf
   \eqref{h_Euc_matrix})
   \begin{eqnarray*}
      \left[
      \begin{array}{ll}
         0             & -\tilde p/\psi\\
         \psi/\tilde p & 0
      \end{array}\right]
   \end{eqnarray*}
   and so we see that $h$ is determined by $\psi$ (which is determined by $S$)
   and $\tilde p$ (which is determined by the points $p_i$) as required. \myqed
\begin{remark}\label{twst_gauging_rmk}
   In the course of the above theorem we made a choice of section $s$ of
   $L^{2m+k}(0, 2l)$ on $S$ with divisor $D$ and such that $ss^* = 1$. This choice
   is unique up to a factor of $U(1)$. Thus the space of all spectral data $(S, D)$
   together with a choice of such a section $s$ is naturally a $U(1)$ bundle over
   the space of spectral data $(S, D)$, ie: over the monopole moduli space. This
   natural $U(1)$ bundle on the monopole moduli space is of course the gauged
   monopole moduli space and is the space whose geometry we are interested in.
\end{remark}
\section{The charge $1$ moduli space}
We are interested in the natural geometric structure on the moduli space
of hyperbolic monopoles\footnote{When we talk of the monopole moduli space,
we shall always mean the gauged moduli space as discussed in remark
\ref{twst_gauging_rmk}.}. The moduli space for $k=1$ (non-singular)
monopoles on $\Hyp$ is simply $\Hyp\times S^1$ but the $k=1$
moduli spaces for singular monopoles are more interesting and have
essentially already been studied in some detail by LeBrun \cite{MR1114461}\footnote{
In the notation of \eqref{def_V}, LeBrun studied the case with $\lambda = 1$ and
$l_1 = \cdots = l_n = 1$. As we shall see for monopoles of mass $m > 0$, we
have $\lambda = 1 + 2m$.}. The goal is to gain some insight into the natural
geometric structure of the moduli spaces of hyperbolic monopoles in general by studying
this simpler ($k=1$) case.
\par
The easiest way to identify the moduli space of $k=1$ singular $SU(2)$ monopoles is
to identify the twistor space. The spectral curves $S$ have genus $0$ for $k=1$ and
the spectral data $(S, D)$ of a monopole can be naturally identified with a twistor line.
We are most interested in the local differential geometric structure of the monopole
moduli space and so it is enough to work with generic spectral curves. Furthermore, the
twistor space naturally carries some additional structure which defines the geometry of
the moduli space which is what we're really interested in. In fact, this twistor space
arises naturally from the Hitchin--Ward transform for a certain singular $U(1)$ monopole.
\subsection{Hitchin--Ward transform for singular $U(1)$ monopoles}\label{HWSingU1Sect}
In this subsection, we will review the construction of certain twistor spaces introduced
by LeBrun in \cite{MR1114461} and show how it can be viewed from the point of view of
the Hitchin--Ward transform for singular $U(1)$ monopoles. In the next subsection, we
will make use of the observation that a twistor line in this space neatly encodes the
spectral data of a charge $1$ singular hyperbolic $SU(2)$ monopole to study the
geometry of the moduli space of such monopoles.
\par
Let $\{p_1,\ldots ,p_n\} \subset \Hyp$ be $n$ distinct points in hyperbolic
space, let $\{l_1,\ldots ,l_n\} \subset \N$ be $n$ (strictly) positive
integers and let $\lambda \in [0,\infty)$ be a non-negative real number.
Let $G_i = G_{p_i}$ be the Green's function for the hyperbolic Laplacian introduced
in equation \eqref{Hyp_Lapl_Green}.
A solution to the $U(1)$ Bogomolny equations is just a solution to the
Laplace equation so that
\begin{eqnarray}\label{def_V}
   V = \lambda + \sum l_i G_i
\end{eqnarray}
defines a $U(1)$ monopole on $U = \Hyp\setminus \{p_1,\ldots ,p_n\}$. We can
thus apply the $U(1)$ version of theorem \ref{hw_bogo} and so obtain the triple
$(K^+, K^-, h)$ where $K^\pm$ are holomorphic {\it line} bundles on $\Q$ and
\begin{eqnarray*}
   h : K^+|_{\Q\setminus P} \to K^-|_{\Q\setminus P}
\end{eqnarray*}
is a holomorphic isomorphism. We can encode this data in a singular 3 dimensional
complex space $\widetilde Z$ as
\begin{eqnarray*}
   \widetilde Z = \left\{ (x,y) \in \left(K^+\oplus (K^-)^*\right)|_{Q\setminus P}
   \quad | \quad h(x)y = 1\right\}
\end{eqnarray*}
\par
Now we have already seen that if $V = 1$ then $K^+ = K^- = L = \cO(1, -1)$ with
$h$ the identity map. It follows (for example from the linearity of the $U(1)$
Penrose transform) that if $V = \lambda$ then $K^+ = K^- = L^\lambda$ with $h$
the identity map. Next, it follows from the work of LeBrun in \cite{MR1114461}
that $V = G_i$ corresponds to $K^+ = \cO(0, 1)$, $K^- = \cO(-1, 0)$. Thus,
in  this case $h \in H^0(\Q\setminus P, \cO(-1, -1))$ and if $\tilde p_i \in
H^0(\Q, \cO(1, 1))$ is a section corresponding to $p_i \in \Hyp$ then by the $U(1)$
version of theorem \ref{hw_h_sings} $\tilde p_i h$ is a non-vanishing holomorphic
function on $\Q$. It is thus constant and we may take this constant to be 1.
We thus have $h = \tilde p_i^{-1}$. Combining these observations we find that
for $V$ given by equation \eqref{def_V} above we have
\begin{eqnarray*}
   K^+ &\simeq& L^\lambda(0, l)\\
   K^- &\simeq& L^\lambda(-l, 0)\\
   h &=& \tilde p^{-1} \in H^0(\Q\setminus P, \cO(-l, -l))
\end{eqnarray*}
where $\tilde p = \prod_{i=1}^{n} \tilde p_i^{l_i} \in H^0(\Q, \cO(l, l))$ and
$l = \sum_{i=1}^n l_i$ as in theorem \ref{hw_h_sings}. We can thus explicitly
identify $\widetilde Z$ as
\begin{eqnarray}\label{twst_space_def}
   \widetilde Z \simeq \left\{(x,y) \in \left(L^\lambda(0, l)\oplus L^{-\lambda}
   (l, 0)\right)|_{\Q\setminus P}\quad |\quad xy = \tilde p(u)\right\}
\end{eqnarray}
(where $(x,y)$ is in the fibre of $L^\lambda(0, l)\oplus L^{-\lambda}(l, 0)$ over
$u \in \Q$.) 
\par
Since the space $\widetilde Z$ encodes the data $(K^+, K^-, h)$ which determines our
singular $U(1)$ monopole, the natural question is how to recover the monopole
directly from $\widetilde Z$. This question has essentially already been
answered by LeBrun \cite{MR1114461} (using ideas of Hitchin \cite{MR520463}.)
A natural desingularisation $Z$ of $\widetilde Z$ is the twistor space for the
4 dimensional real manifold that is the total space of the principal $U(1)$
bundle $M$ of the monopole we started with, endowed with a Gibbons--Hawking
type conformal structure (introduced by LeBrun \cite{MR1114461}). Since we
are following the setup of \cite{MR1114461}, we shall summarise the details
in the language which will be most useful to us, but without providing proofs.
\par
Note also that the fact that we can view $\widetilde Z$ as arising from
our version of the Hitchin--Ward transform for singular $U(1)$ monopoles, though
satisfying, is not relevant to the rest of our work here. Equation
\eqref{twst_space_def} essentially appears in \cite{MR1114461} and may be taken
as the starting point of our work here.
\par
In order to connect the space $\widetilde Z$ with our $SU(2)$ monopole moduli
space, we will need to understand the real structure and the twistor lines.
Of course it is LeBrun's desingularisation $Z$ of $\widetilde Z$
that is the twistor space, however the real structure and twistor lines are
first defined on $\widetilde Z$ and then lifted to $Z$ so we shall work on
$\widetilde Z$. This is also consistent with the approach taken by Hitchin in
\cite{MR520463}.
\par
We define the real structure
\begin{eqnarray*}
   \hat \tau : \widetilde Z \to \widetilde Z
\end{eqnarray*}
by restricting the anti-holomorphic (anti-linear) map $L^\lambda(l, 0) \to
L^{-\lambda}(0, l)$ induced by $\sigma$ on $\Q$. This restricts to $\widetilde
Z$ since $\tilde p$ is real.
\par
To identify the twistor lines in $\widetilde Z$, it is convenient to use
the following elementary lemma.
\begin{lemma}\label{p_factor_lemma}
   In the above notation, let $P_q \subset \Q$ be the twistor line of a
   point $q \in U$. Then we can consistently choose $x, y \in H^0(\CP, \cO(l))$
   unique up to a factor of $U(1)$ such that $x = y^*$ and
   \begin{eqnarray}
      (\pi_1^*x)(\pi_2^*y) = \tilde p
   \end{eqnarray}
   on $P_q$.
\end{lemma}
\noindent {\bf Proof} Use $q$ to identify $\Q \simeq \QQ$ so that $P_q$
corresponds to the diagonal $\Delta$. Let $\zeta$ be the natural coordinate
on $\Delta \simeq \CP$, then restricted to $\Delta$, $\tilde p$ is given by
\begin{eqnarray*}
   \tilde p = \prod (a_i\zeta^2 + 2b_i\zeta - \overline a_i)^{l_i}
\end{eqnarray*}
where $a_i \in \C$ and $b_i \in \R$. We now simply follow the recipe given
in \cite{MR520463}. The discriminant $4(b_i^2 + |a_i|^2)$ of the $i^{\rm th}$ 
quadratic is positive and so we can without ambiguity define $\Delta_i
= \sqrt{b_i^2 + |a_i|^2}$ to be the positive square root and
the roots $\alpha_i$ and $\beta_i$ by
\begin{eqnarray*}
   \alpha_i &=& \frac{-b_i + \Delta_i}{a_i}\\
   \beta_i &=& \frac{-b_i - \Delta_i}{a_i}
\end{eqnarray*}
We then define
\begin{eqnarray*}
   x &=& A\prod(\zeta - \alpha_i)^{l_i}\\
   y &=& B\prod(\zeta - \beta_i)^{l_i}
\end{eqnarray*}
where $A$ and $B$ satisfy $AB = \prod a_i^{l_i}$. $x$, $y$ now satisfy
the required conditions if and only if
\begin{eqnarray*}
   |A|^2 = \prod(b_i - \Delta_i)^{l_i}
\end{eqnarray*}
We thus have the required factoring of $\tilde p$ and the indeterminacy
in the argument of $A$ corresponds to the extra $U(1)$ factor. \myqed
Now we also have a natural holomorphic projection
$\widetilde Z \to \Q\setminus P$, compatible with real structures. Composing
this with the projection $\pi_1$ of $\Q \subset \CP\times\CP$ onto the first
factor we thus have a holomorphic projection
\begin{eqnarray*}
   \pi : \widetilde Z \to \CP
\end{eqnarray*}
(Note that if we instead used the projection $\pi_2$ onto the second factor
in $\CP\times\CP$, we would just obtain the map $\tau\circ\pi\circ\hat\tau$
so that provided we remember the real structures on $\widetilde Z$ and $\CP$,
there is no new information.)
\par
We can now exhibit the family of twistor lines in $\widetilde Z$.
Each twistor line will be the image of a holomorphic
section of $\pi$. As we shall see, this has geometrical significance for $M$.
Indeed if $\pi$ preserved the real structures on $\widetilde Z$ and $\CP$, then
(see \cite{MR1631124}) $M$ would have a hypercomplex structure. This is, of course,
not the case here. We define these sections as follows. Pick a point $q \in U =
\Hyp \setminus \{p_1, \ldots, p_n\}$. Let $P_q$ be the corresponding twistor line
in $\Q$. Note that
\begin{eqnarray*}
   \pi_i : P_q \to \CP
\end{eqnarray*}
is a holomorphic bijection for $i=1,2$. (A geometric reason for this is that
$P_q$ is the set of all geodesics through the point $q \in \Hyp$ so if we know
the start or end point of the geodesic on sphere at infinity in $\Hyp$, then we
know the geodesic.) Now choose $x,y \in H^0(\CP, \cO(l))$ as in lemma
\ref{p_factor_lemma} and let
\begin{eqnarray*}
   s : P_q \to L
\end{eqnarray*}
be a holomorphic trivialisation of $L|_{P_q}$ such that $ss^* = 1$. Evidently $s$
is unique up to a factor of $U(1)$. $s$ defines trivialisations $s^\lambda$ of
$L^\lambda$ and $s^{-\lambda}$ of $L^{-\lambda}$ over $P_q$. The pair of sections
\begin{eqnarray*}
   (s^\lambda \pi_2^* x, s^{-\lambda} \pi_1^* y) \in H^0(P_q, L^\lambda(0, l)
   \oplus L^{-\lambda}(l, 0))
\end{eqnarray*}
define a real lifting of the twistor line $P_q$ in $\Q$ to the required twistor
line in $\widetilde Z$.
\par
With this explicit description of the twistor lines of $\widetilde Z$, we are
ready to identify them with the spectral data of charge 1 singular hyperbolic
$SU(2)$ monopoles. This means that $Z$ is the twistor space of the moduli space
$M$ of these monopoles. Furthermore, LeBrun \cite{MR1114461} explicitly identified
the space which $Z$ is the twistor space of and so by studying its natural geometry
we are studying the natural geometry of the monopole moduli space. We summarise
here the geometry of $M$.
\par
We saw in section \ref{SngMonopDefsSect} that a harmonic function $V$ as in 
equation \eqref{def_V} defines a Riemannian manifold. The manifold $M$ is the total
space of the $U(1)$ bundle on $U$ with Chern class $\frac{1}{2\pi}[*dV]$ and
the metric (which is anti-self-dual with respect to the orientation defined by
$dx\wedge dy\wedge dz\wedge\omega$) is
\begin{eqnarray}\label{Gibb_Hawk_conf}
   \hat V \hat h + \hat V^{-1}\omega\otimes \omega
\end{eqnarray}
where $\hat h$ is the pull back of the hyperbolic metric $h$ to $M$, $\hat V$
is the pull back of $V$ to $M$ and $\omega$ is a connection on $M$ with curvature
$*dV$. According to \cite{MR1114461}, the twistor space of $M$ with the conformal
structure defined by this metric is $Z$ (the natural desingularisation of $\widetilde
Z$).
\par
We noted above that the fact that the twistor lines were the images of holomorphic
sections of the natural holomorphic projection $\pi : \widetilde Z \to \CP$ would
have geometric significance for $M$. In fact if we fix a point $u \in \CP$ then the
fibre $\Sigma = \pi^{-1}(\{u\})$ is a divisor corresponding to a
complex structure on $M$. Furthermore, as shown in \cite{MR1114461}, the divisor
$\Sigma + \overline\Sigma$ represents the line bundle $K^{-1/2}$ where $K$ is the
canonical bundle of $Z$. Since
\cite{MR1050087}, holomorphic sections $H^0(Z, K^{-1/2})$ correspond to
scalar-flat K\"ahler metrics\footnote{Strictly speaking an element of
$H^0(Z, K^{-1/2})$ gives a pair of complex structures $J$, $-J$ on $M$ and we
cannot tell them apart. So the K\"ahler metric is well defined but the complex
structure is only defined up to conjugation. This is not an issue for us however
since we have $\Sigma$.}
in the conformal class on $M$, $\Sigma$ defines a K\"ahler metric on
$M$ up to scale. We are fortunate that we may appeal to \cite{MR1114461} for an
explicit description of these metrics. The details are as follows.
\par
Choose coordinates for $\Hyp$ so that it is represented in the upper half-space
model
\begin{eqnarray*}
   \Hyp &\simeq& \{ (x,y,z)\in \R^3\quad |\quad z > 0\}\\
   h &=& \frac{dx^2 + dy^2 + dz^2}{z^2}
\end{eqnarray*}
Suppressing the notation the pull back of data from $U$ to $M$,
note that $dx, dy, dz, \omega$ trivialise $T^*M$ so
that we can define an almost complex structure $J$ on $M$ by
\begin{eqnarray*}
   Jdx &=& dy\\
   Jdz &=& z V^{-1}\omega
\end{eqnarray*}
Noting that
\begin{eqnarray*}
   J(dx + idy) &=& -i(dx + idy)\\
   J(z^{-1}Vdz + i\omega) &=& -i(z^{-1}V dz + i\omega)\\
   d(dx + idy) &=& 0\\
   d(z^{-1}V dz + i\omega) &=& (dx + idy)\wedge\left(\frac{1}{z}
   (\frac{\partial V}{\partial x} - i\frac{\partial V}{\partial y})
   dz + \frac{i}{z}\frac{\partial V}{\partial z}dy\right)
\end{eqnarray*}
we see that $J$ is integrable. Now define the metric
\begin{eqnarray*}
   g = z^2(Vh + V^{-1}\omega\otimes\omega)
\end{eqnarray*}
on $M$ (which is in the conformal class defined by \eqref{Gibb_Hawk_conf}).
We claim that $(M, J, g)$ is a scalar-flat K\"ahler manifold. If we let
$\Omega = g(J\cdot, \cdot)$ be the associated 2-form then a quick
calculation reveals
\begin{eqnarray}\label{Kahler_form}
   \Omega = -(Vdx\wedge dy + zdz\wedge \omega)
\end{eqnarray}
from which it follows easily (using the fact that $d*dV = 0$) that $d\Omega = 0$.
Thus $g$ is a K\"ahler metric and we need only note that an anti-self-dual K\"ahler
manifold is scalar-flat (see for example \cite{MR707181}).
\par
Finally note that to define this metric on $M$, we
represented $\Hyp$ in the upper half-space model. This singles out a point on the
conformal 2-sphere at infinity and so we in fact have an entire $S^2$ of such
metrics. This corresponds to the choice of point in $\CP$ giving the divisor
$\Sigma$ above. To see this more clearly, note that the K\"ahler form
\eqref{Kahler_form} can be written as
\begin{eqnarray*}
   \Omega = -\frac{1}{2}(V*d(z^2) + d(z^2)\wedge\omega)
\end{eqnarray*}
where $*$ is the Hodge $*$-operator on $\Hyp$. Thus to define the K\"ahler
structure, we only need the function $z$. Furthermore, $z$ is just a horospherical
height function on $\Hyp$ and as we shall see, any such function will define a
K\"ahler structure. Let us briefly recall some elementary facts about horospherical
height functions on $\Hyp$.
\par
A horospherical height function is the exponential of a Busemann function.
To define a Busemann function on a complete Riemannian manifold with
distance function $\rho$, we choose a geodesic $\gamma$ parameterised by
arc-length and define the associated Busemann function $b_\gamma$ by
\begin{eqnarray}\label{Busemann_def}
   b_\gamma(x) = \lim_{t \to \infty}\left(t - \rho(x, \gamma(t))\right)
\end{eqnarray}
Note that $b_\gamma(\gamma(t)) = t$ and that if we change the parameterisation of
$\gamma$ by $t \mapsto \gamma(t + a)$ for some constant $a$ then $b_\gamma$ is
replaced with $b_\gamma - a$. Let us consider Busemann functions on $\Hyp$. We
claim that the level sets of $b_\gamma$ are the horospheres passing through
$\gamma(\infty) \in \partial\Hyp$. To see this, we may choose
coordinates (x, y, z) so that $\Hyp$ is represented in the upper half-space 
model and our geodesic $\gamma$ appears as
\begin{eqnarray*}
   \gamma : t \mapsto (0, 0, e^t)
\end{eqnarray*}
The horospheres tangential to $\gamma(\infty)$ are given in these
coordinates by $z=const$. Using the formula
\begin{eqnarray*}
   \rho((x,y,z), \gamma(t)) = \cosh^{-1}\left(\frac{x^2 + y^2 + z^2 + e^{2t}}{2ze^t}\right)
\end{eqnarray*}
it is easy to check from \eqref{Busemann_def} that
\begin{eqnarray*}
   e^{b_\gamma((x,y,z))} = z
\end{eqnarray*}
so that $z$ is the associated horospherical height function as required.
Since, if we fix a point $p$ on the boundary, $\Hyp$ is foliated by the horospheres
tangential to $\partial\Hyp$ at $p$, we see that to define a Busemann function on $\Hyp$, we
just need to define a value on each horosphere. Choosing a geodesic such that
$\gamma(\infty) = p$ we define the value on each horosphere using $b_\gamma(\gamma(t)) = t$.
The only remaining degree of freedom is which horosphere contains $\gamma(0)$. Thus
$p$ determines the Busemann function up to the addition of a constant.
\par
We have thus seen that a horospherical height function is determined up to a positive
scale factor by a point on $\partial\Hyp$ and that any such function can be represented
as the function $z$ in upper half-space coordinates. If we fix a point $O \in \Hyp$, we
thus have a natural $2$-sphere of horospherical height functions: for each point on 
$\partial\Hyp$ we take the associated horospherical height function $q$ such that $O$ lies
in the level set $q=1$. In this way, a point in $\Hyp$ determines a natural $2$-sphere of
scalar-flat K\"ahler metrics on $M$.
\subsection{The $SU(2)$ moduli space}
In the previous subsection, we reviewed a construction of LeBrun \cite{MR1114461}
(and also noted how it fits with our singular $U(1)$ Hitchin--Ward transform). The
essential point was the geometry of the space $M$ and its twistor space. We now
wish to show that $M$ is in fact the moduli space of $k=1$ singular $SU(2)$ monopoles.
We will do this by showing how the spectral curve $S$ of such a monopole naturally
lifts to a twistor line in $\widetilde Z$ and that this twistor line also encodes
the spectral divisor of the monopole.
\par
Thus consider a singular $SU(2)$ monopole with non-Abelian charge $k=1$ and Abelian
charges $l_i$ at the points $p_i \in \Hyp$, $i=1,\ldots, n$. Let
\begin{eqnarray*}
   V = 1 + 2m + 2\sum l_i G_i
\end{eqnarray*}
So comparing with \eqref{def_V}, we have $\lambda = 1 + 2m$ and $2l_i$ in place of
$l_i$. As we saw in the previous subsection, associated with $V$ is a twistor space
\begin{eqnarray*}
   \widetilde Z \simeq \left\{(x,y) \in \left(L^{1+2m}(0, 2l)\oplus L^{-1-2m}
   (2l, 0)\right)|_{\Q\setminus P}\quad |\quad xy = \tilde p\right\} 
\end{eqnarray*}
The key observation is that a twistor line of $\widetilde Z$ exactly encodes the
spectral data $(S, D)$ for a charge $1$ singular hyperbolic monopole. We need only
observe that $D$ is the divisor of a section
\begin{eqnarray*}
   \xi^+ \in H^0(S, (L^-)^{*}L^+(l, l))
\end{eqnarray*}
and using the isomorphisms $L^+ \simeq L^m(0, -N)$ and $L^- \simeq L^{-m}(-N, 0)$
we see that we can regard
\begin{eqnarray*}
   \xi^+ \in H^0(S, L^{1+2m}(0,2l))
\end{eqnarray*}
Similarly we have
\begin{eqnarray*}
   \xi^- \in H^0(S, L^{-1-2m}(2l, 0))
\end{eqnarray*}
Thus $\xi^\pm$ provide a natural lifting of the spectral curve (which is a
rational curve) to the twistor space $\widetilde Z$. The sections $\xi^\pm$
are determined up to $U(1)$ factor by the spectral divisor $D$. This $U(1)$
factor, of course, corresponds to the gauging of the monopole. This is how
a monopole corresponds to a twistor line, and hence a point in $M$. Conversely,
by construction, a twistor line in $\widetilde Z$ projects down to the spectral
curve of a singular hyperbolic monopole and this curve obviously comes with
a section of $L^{1+2m}(0, 2l)$ which gives the spectral divisor. To summarise,
we thus have
\begin{theorem}\label{Charge1ModSpaceGeom}
   Let $m > 0$, let $\{p_1,\ldots,p_n\} \subset \Hyp$ be $n$ distinct points in
   hyperbolic space and let $\{l_1,\ldots,l_n\} \subset \N$ be $n$ (strictly)
   positive integers. Let $M$ be the moduli space of gauged singular hyperbolic
   $SU(2)$ monopoles of non-Abelian charge 1 with Abelian charges $l_i$ at $p_i$.
   Then
   \begin{itemize}
      \item
      $M$ carries a natural self-dual conformal structure
      \item
      For each point $u \in \partial\Hyp$ there is a volume form and complex
      structure $J^u$ on $M$ such that the metric determined in the conformal
      structure makes $M$ together with $J^u$ a scalar-flat K\"ahler manifold.
   \end{itemize}
\end{theorem}

\chapter{Deforming the spectral curve }\label{DefSpecChap}
\section{Overview}\label{DefSpecOverview}
Ever since \cite{MR709461} it has been known that the $(4k-1)$-dimensional
moduli space of charge $k$ $SU(2)$ Euclidean monopoles can be naturally
identified with the space of spectral curves, ie: those compact algebraic
curves $S$ in the linear system $|\cO(2k)|$ on $\T \simeq T\CP$ such that
\begin{itemize}
   \item
   $S$ is real
   \item
   $S$ has no multiple components
   \item
   $L^2$ is trivial on $S$ and $L(k-1)$ is real
   \item
   $H^0(S,L^z(k-2))=0$ for $z\in(0,2)$
\end{itemize}
Furthermore, this identification lifts naturally to yield an identification
of the $4k$-dimensional moduli space of \emph{gauged} Euclidean monopoles
and the space of pairs consisting of a spectral curve $S$ as above together
with a choice of trivialisation of $L^2$ (satisfying a reality condition)
over $S$.
\par
One of the most interesting features of the moduli space of gauged
Euclidean monopoles is the natural hyperk\"ahler structure it carries. This
structure has previously been understood from at least three points of view:
the existence of the natural $L^2$ metric on the moduli space, the fact that
the Bogomolny equations can be regarded as moment map equations for an
infinite dimensional hyperk\"ahler quotient or from a construction of the
twistor space of the moduli space. These facts are well laid out in
\cite{MR934202}. However the question of understanding the hyperk\"ahler
structure from the point of view of the spectral curve has not been previously
addressed.
\par
This is an interesting question because the spectral curve approaches to Euclidean
and hyperbolic monopoles (see \cite{MR649818} and \cite{MR1399482} as well as
\cite{MR893593}) have much in common. By understanding how the hyperk\"ahler
structure on the space of spectral curves arises in the Euclidean case and by
adapting the same techniques to the hyperbolic case, we are able to learn about the
geometry of the moduli spaces of (gauged) hyperbolic monopoles. This is important
because the natural geometric structure on the moduli space of hyperbolic monopoles
remains a mystery to date (indeed the natural $L^2$ \lq\lq metric\rq\rq diverges,
see \cite{MR1010167}).
\par
Consider the Euclidean case. Of the conditions listed overleaf that characterise
a curve in $|\cO(2k)|$ as a spectral curve, the most important is the condition
that $L^2|_S$ be trivial. Indeed as pointed out in \cite{MR649818}, a naive
parameter count indicates that the space of such curves has dimension $4k$.
Our idea is to view the $L^2|_S$ triviality condition more geometrically by
using a trivialisation (satisfying a reality condition) to lift $S$ to a real curve
$\hat S$ in $L^2\setminus 0$. We show that it is possible to apply Kodaira's
deformation theory to $\hat S$ (at least when it is smooth) and thus obtain the
space of curves in $|\cO(2k)|$ satisfying the $L^2|_S$ triviality condition as a
deformation space. This means that we have a model for the tangent space of such
curves and so for the monopole moduli space from which we are able to identify the
hyperk\"ahler structure.
\par
It is worth pointing out that in the charge 1 case, the spectral curves have genus 0
and lift to twistor lines in $L^2\setminus 0$ which is of course the twistor space
of $\R^3\times S^1$ (the charge 1 gauged moduli space). The equivalent statement in
the hyperbolic case also holds and our work in chapter \ref{SingHypChap} generalised
this to \emph{singular} charge 1 hyperbolic monopoles (whose spectral curves still
have genus 0). Our approach here is thus a generalisation in a different direction,
a generalisation that involves deformations of higher genus curves. 
\par
As we will see, although questions about the geometry of the hyperbolic monopole
moduli space remain, this approach does yield significant insights. We find that
the hyperbolic monopole moduli space appears to carry a new type of geometry whose
complexification is very similar to the complexification of hyperk\"ahler geometry
but with different reality conditions.
\par
Finally it must be pointed out that much of the work in this approach goes into
proving a cohomology vanishing theorem. For example, in the Euclidean case
we prove
\begin{eqnarray*}
   H^0(S, \widetilde EL(k-2)) = 0
\end{eqnarray*}
where $\widetilde E$ is the holomorphic bundle on $\T$ corresponding to the
monopole. A corollary of this result is that we get a proof that the space of smooth
curves in $|\cO(2k)|$ satisfying the $L^2|_S$ triviality condition has dimension $4k$
and so we obtain the dimension of the monopole moduli space without having to appeal
to the analytical results of Taubes \cite{MR723549}.

\section{The Euclidean case}
\subsection{Recovering the hyperk\"ahler structure}
Let $S \subset \T$ be the spectral curve of a charge $k$ Euclidean monopole. For
simplicity we shall assume $S$ is non-singular for the rest of this chapter.
Let $\phi$ be a non-vanishing holomorphic section of $L^2$ on $S$ such that
\begin{eqnarray*}
   \phi^*\phi = -1
\end{eqnarray*}
Note that $\phi$ is unique up to a factor of $U(1)$. Let
\begin{eqnarray}\label{spec_crv_lift}
   \hat S \subset L^2\setminus 0
\end{eqnarray}
be the image
of $\phi$ and let $\hat N$ be the normal bundle of $\hat S$ in $L^2 \setminus 0$.
As we said in section \ref{DefSpecOverview} we are going to consider the space
$M_k^\C$ of deformations of $\hat S$ in $L^2\setminus 0$. Since $L^2\setminus 0$
carries a real structure, the same is true of $M_k^\C$. We will find that $M_k^\C$
has dimension $4k$ and that in the neighbourhood of a genuine spectral curve, the
real points can be naturally identified with an open set in the moduli space $M_k$
of monopoles. We thus have a natural isomorphism
\begin{eqnarray*}
   T_{\hat S}M_k^\C \simeq T_{\hat S} M_k\otimes_{\R}\C
\end{eqnarray*}
Now to get the deformation theory to work, we appeal to a well known result of
Kodaira \cite{MR0133841} which states that if
\begin{eqnarray}\label{def_obs_cond}
   H^1(\hat S, \hat N) = 0
\end{eqnarray}
then we have a well behaved\footnote{By \lq\lq well behaved\rq\rq we mean that
$M_k^\C$ is a complete maximal family of deformations. See \cite{MR0133841}.}
space $M_k^\C$ of deformations of $\hat S$ in $L^2\setminus 0$ and furthermore there
is a natural isomorphism
\begin{eqnarray*}
   T_{\hat S}M_k^\C \simeq H^0(\hat S, \hat N)
\end{eqnarray*}
We thus have a model for the complexified tangent space to the moduli space
of monopoles
\begin{eqnarray}\label{Euc_tang_spc_model}
   T_{\hat S}M_k\otimes_{\R}\C \simeq H^0(\hat S, \hat N)
\end{eqnarray}
Evidently, we must address the issue of \eqref{def_obs_cond}. We shall find
that it does indeed hold. As a first step to establishing this result, we
need to identify $\hat N$. To this end we have some lemmas.
\par
\begin{lemma}\label{NormExtLemma}
   Let $X$ be a complex manifold, $S \subset X$ a complex submanifold and
   $\pi : V \to X$ a holomorphic vector bundle on $X$. Suppose that we have
   a section $\phi \in H^0(S, V)$ with image $\hat S \subset V$. Let $N$ be
   the normal bundle of $S$ in $X$ and $\hat N$ be the normal bundle of
   $\hat S$ in $V$. Then (regarding $\hat N$ as a bundle on $S$ using $\phi$)
   $\hat N$ naturally fits into the short exact sequence of vector bundles
   on $S$
   \begin{eqnarray}\label{NormExtSeq}
      0 \to V|_S \to \hat N \to N \to 0
   \end{eqnarray}
\end{lemma}
\noindent {\bf Proof}
The three natural exact sequences
\begin{eqnarray*}
   &0& \to \pi^*V \to T V \to \pi^* T X \to 0\\
   &0& \to T S \to T X|_S \to N \to 0\\
   &0& \to T\hat S \to T V|_{\hat S} \to \hat N \to 0
\end{eqnarray*}
fit together into the commutative diagram
\begin{eqnarray*}
   \begin{CD}
      &&&& 0 && 0\\
      &&&& @VVV @VVV\\
      &&&& T\hat S @>\simeq >> T S\\
      &&&& @VVV @VVV\\
      0 @>>> \pi^*V|_{\hat S} @>>> T V|_{\hat S} @>>> \pi^*T X|_S @>>> 0\\
      &&&& @VVV @VVV\\
      &&&& \hat N && N\\
      &&&& @VVV @VVV\\
      &&&& 0 && 0
   \end{CD}
\end{eqnarray*}
A quick diagram chase reveals that there exists a unique map $\hat N \to N$
for which the diagram commutes and furthermore that the resulting sequence
\begin{eqnarray*}
   0 \to V|_{S} \to \hat N \to N \to 0
\end{eqnarray*}
is exact. \myqed
\begin{lemma}\label{NormIdLemma}
   In the notation of lemma \ref{NormExtLemma}, suppose that $V = L$ is a
   line bundle and that $S$ has codimension $1$.
   Let $\alpha \in H^1(S, LN^*)$ be the extension class of the sequence
   (\ref{NormExtSeq}) and let $\delta : H^0(S, L) \to H^1(X, LN^*)$
   be the connecting homomorphism associated to the following short exact
   sequence of sheaves on $X$
   \begin{eqnarray}\label{XSSheafSeq}
      0 \to \cO_X(LN^*) \to \cO_X(L) \to \cO_S(L) \to 0
   \end{eqnarray}
   Then $(\delta \phi)|_{S} = \alpha$.
\end{lemma}
\noindent {\bf Proof}
The extension class $\alpha$ is defined as $\alpha =
\delta^{'} 1$ where $\delta^{'} : H^0(S, \cO) \to H^1(S, LN^*)$ is the
connecting homomorphism associated to the short exact sequence obtained by
tensoring (\ref{NormExtSeq}) with $N^*$. We shall calculate an explicit \v{C}ech
cocycle representative for $\delta\phi$ and for $\delta^{'}1$ with respect to a
Leray cover $\cU = \{U_i\}$ for $X$ and show that they define the same
class in $H^1(S, LN^*)$.
\par
Thus let $\cU$ be a sufficiently fine Leray cover for $X$ and let the
functions $\psi_i : U_i \to \C$ cut out $S_i = S\cap U_i$ for each open set $U_i
\in \cU$. To obtain an explicit \v{C}ech representative for $\delta\phi$ we must
consider the commutative diagram of \v{C}ech cochain groups associated to the
sequence
(\ref{XSSheafSeq})
\begin{eqnarray*}
   \begin{CD}
      0 @>>> \underset{i}\prod H^0(U_i, LN^*) @>>> \underset{i}\prod H^0(U_i, L)
      @>>> \underset{i}\prod H^0(S_i, L) @>>> 0\\
      && @VVV @VVV @VVV\\
      0 @>>> \underset{i,j}\prod H^0(U_{ij}, LN^*) @>>> \underset{i,j}\prod
      H^0(U_{ij}, L) @>>> \underset{i,j}\prod H^0(S_{ij}, L) @>>> 0\\
   \end{CD}
\end{eqnarray*}
where $U_{ij} = U_i \cap U_j$ and $S_{ij} = S_i \cap S_j$. Now fix trivialisations
of $L$ on each open set $U_i$ and let $l_{ij} : U_{ij} \to \C^*$ be the transition
functions. Let $\phi_i : S_i \to \C$ represent $\phi|_{S_i}$ with respect to these
trivialisations and let $\hat\phi_i : U_i \to \C$ be a holomorphic extension of
$\phi_i$ to $U_i$. Then, chasing through the above diagram, we see that a \v{C}ech
representative for $\delta\phi$ is
\begin{eqnarray*}
(\delta\phi)_{ij} = \frac{\hat\phi_i - l_{ij}\hat\phi_j}{\psi_i}
\end{eqnarray*}
Note that we have used the chosen trivialisation of $L$ on $U_i$ and the
trivialisation of $N$ on $S_i$ determined by $\psi_i$ so that $(\delta\phi)_{ij}$
takes values in $\C$ rather than $LN^*$.
\par
With this in hand, we turn to the calculation of a \v{C}ech cocycle for $\delta^{'}1$.
The commutative diagram of \v{C}ech cochain groups associated to the exact sequence
with $\delta^{'}$ as connecting homomorphism is
\begin{eqnarray*}
   \begin{CD}
      0 @>>> \underset{i}\prod H^0(S_i, LN^*) @>>> \underset{i}\prod
      H^0(S_i, \hat NN^*) @>b>> \underset{i}\prod H^0(S_i, \cO) @>>> 0\\
      && @VVV \partial @VVV @VVV\\
      0 @>>> \underset{i,j}\prod H^0(S_{ij}, LN^*) @>a>> \underset{i,j}\prod
      H^0(S_{ij}, \hat NN^*) @>>> \underset{i,j}\prod H^0(S_{ij}, \cO) @>>> 0\\
   \end{CD}
\end{eqnarray*}
Now note that the trivialisation of $L$ on $U_i$ together with the extension
$\hat\phi_i$ of $\phi_i$ together determine an isomorphism
\begin{eqnarray*}
   \hat N|_{S_i} &\simeq& N|_{S_i} \oplus \cO_{S_i}\\
   \left[(v,w )\right] &\mapsto& (\left[ v\right], w - v(\hat\phi_i))
\end{eqnarray*}
where $(v,w) \in TU_i \oplus \C \simeq T(U_i\times \C) \simeq T(L|_{U_i})$.
Also the explicit trivialisation of $N|_{S_i}$ determined by $\psi_i$ is
\begin{eqnarray*}
   N|_{S_i} &\simeq& \cO_{S_i}\\
   \left[ v\right] &\mapsto& v(\psi_i)
\end{eqnarray*}
where $v \in TU_i$. With these isomorphisms in place, each bundle appearing
in our latest commutative diagram is trivialised.
\par
We wish to calculate $\delta^{'}1$, thus take the cocycle $\{1\}\in \underset{i}
\prod H^0(S_i, \cO)$. Using our isomorphisms $f_i: H^0(S_i, \hat NN^*) \simeq
H^0(S_i, \cO)^2$ a cochain that is mapped to this under $b$ is $\{(1, 0)\}$. We
must apply the \v{C}ech coboundary map to this. A careful examination of the
definition of $f_i$ reveals
\begin{eqnarray*}
   \partial \{(1, 0)\} = \left\{\left(1 - \frac{\psi_j}{\psi_i}v_j(\psi_i),
   \frac{\psi_j}{\psi_i}v_j(\hat\phi_i - l_{ij}\hat\phi_j)\right)\right\}
\end{eqnarray*}
where $v_j \in H^0(S_{ij}, TU_{ij})$ satisfies $v_j(\psi_j) = 1$. The
cocycle mapped to this under $a$ represents $\delta^{'}1$. In fact we can read off
\begin{eqnarray*}
   (\delta^{'}1)_{ij} &=& \frac{\psi_j}{\psi_i}v_j(\hat\phi_i - l_{ij}\hat\phi_j)\\
   &=& \frac{\psi_j}{\psi_i}v_j\left(\frac{\hat\phi_i - l_{ij}\hat\phi_j}
   {\psi_j}\psi_j\right)\\
   &=& \frac{\psi_j}{\psi_i}\frac{\hat\phi_i - l_{ij}\hat\phi_j}{\psi_j}v_j (\psi_j) +
   \frac{\psi_j}{\psi_i}v_j\left(\frac{\hat\phi_i - l_{ij}\hat\phi_j}{\psi_j}\right)
   \psi_j\\
   &=& \frac{\hat\phi_i - l_{ij}\hat\phi_j}{\psi_i}
\end{eqnarray*}
since $v_j(\psi_j) = 1$ and $\psi_j$ vanishes on $S_{ij}$.
\par
Thus the same cocycle represents both $\delta\phi$ and $\delta^{'}1$. \myqed
\begin{corollary}\label{Euc_norm_bnd_id}
   Let $\hat N$ be the normal bundle of $\hat S \subset L^2\setminus 0$ as in
   \eqref{spec_crv_lift}. Let $\widetilde E$ be the holomorphic vector
   bundle on $\T$ corresponding to the monopole. Then (identifying $\hat S$ and
   $S$)
   \begin{eqnarray*}
      \hat N \simeq \widetilde EL(k)|_{S}
   \end{eqnarray*}
\end{corollary}
\noindent {\bf Proof} This is an immediate consequence of previous lemma together
with Hitchin's construction of $\widetilde E$ from $S$ (see \cite{MR649818}). \myqed
\begin{remark}
   Note that since the space of deformations of spectral curves has a natural map
   to the space of deformations of holomorphic bundles on $\T$ there should
   be a natural map
   \begin{eqnarray*}
      H^0(S, \hat N) \to H^1(\T, \End(\widetilde E))
   \end{eqnarray*}
   and in view of corollary \ref{Euc_norm_bnd_id} this means that there should be
   a natural map
   \begin{eqnarray*}
      H^0(S, \widetilde EL(k)) \to H^1(\T, \End(\widetilde E))
   \end{eqnarray*}
   It is instructive to check that this is indeed the case.
   Thus note that since the spectral curve is a divisor of $\cO(2k)$ we have the
   short exact sequence of sheaves on $\T$
   \begin{eqnarray*}
      0 \to \cO_T(\widetilde EL(-k)) \to \cO_T(\widetilde EL(k)) \to \cO_S
      (\widetilde EL(k)) \to 0
   \end{eqnarray*}
   The connecting homomorphism of the induced long exact sequence of cohomology
   is a map
   \begin{eqnarray*}
      H^0(S, \widetilde EL(k)) \to H^1(\T, \widetilde EL(-k))
   \end{eqnarray*}
   and since $L(-k)$ is a subbundle of $\widetilde E$ we also have a map
   \begin{eqnarray*}
      H^1(\T, \widetilde EL(-k)) \to H^1(\T, \End(\widetilde E))
   \end{eqnarray*}
   (where we have used the fact that $\widetilde E^* \simeq \widetilde E$). We
   thus have a natural map
   \begin{eqnarray*}
      H^0(S, \widetilde EL(k)) \to H^1(\T, \End(\widetilde E))
   \end{eqnarray*}
   as required.
\end{remark}
\par
Having identified the normal bundle of $\hat S$, the next step is to establish the
vanishing of $H^1(S, \hat N)$ required by Kodaira's theorem \cite{MR0133841}.
In fact we shall prove a stronger cohomology vanishing theorem which will also
be of use elsewhere. First however, it is convenient for us to introduce a
definition and a technical lemma that will be used in the proof of the vanishing
theorem.
\begin{definition}\label{special_gauge}
Fix a charge $k$ $SU(2)$ monopole on $\R^3$ and let $r$ be the radial distance from
a point $O \in \R^3$. We define a special gauge to be a
unitary gauge which diagonalises the Higgs field and such that the connection
matrix takes the form
\begin{eqnarray*}
   \left[\begin{array}{cc}
   A    & B\\
   -B^* & -A
   \end{array}\right]
\end{eqnarray*}
where $B = O(r^{-2})$.
\end{definition}
It follows easily from the boundary conditions given in \cite{MR649818} that special
gauges exist for any monopole.
\begin{lemma}\label{vanishing_lemma}
   Consider an $SU(2)$ monopole on $\R^3$ with underlying complex vector bundle $E$.
   Fix a point $O \in \R^3$ and let $\gamma$ be an oriented line in
   $\R^3$ that is not a spectral line of the monopole. Let $x$ be the point on $\gamma$
   closest to $O \in \R^3$. Let $L^\pm_\gamma \subset E_x$
   be the subspaces of initial conditions to the scattering equation $\nabla_{\dot\gamma}
   - i\Phi = 0$ along $\gamma$ that give solutions decaying at the positive and
   negative ends of $\gamma$. Using the decomposition $E_x = L^+_\gamma \oplus L^-_\gamma$
   define the endomorphism $M_\gamma : E_x \to E_x$ by
   \begin{eqnarray*}
      M_\gamma = \left[
      \begin{array}{cc}
         1 & 0\\
         0 & -1
      \end{array}\right]
   \end{eqnarray*}
         Then $\|M_\gamma\|$ is bounded for all $\gamma$ sufficiently far from $O$ by a
         constant independent of $\gamma$.
\end{lemma}
\noindent {\bf Proof}
Let $s_0$ be a solution to the scattering equation decaying at the
positive end of $\gamma$ and $s'_0$ be a solution decaying at the negative end. Let $v
\in E_x$ and let $v = v_+ + v_-$ where $v_\pm \in L^\pm_\gamma$. Then
\begin{eqnarray*}
   v_+ &=& \frac{\langle v, s'_0(x)\rangle}{\langle s_0, s'_0\rangle}s_0(x)\\
   v_- &=& -\frac{\langle v, s_0(x)\rangle}{\langle s_0, s'_0\rangle}s'_0(x)
\end{eqnarray*}
Now $M_\gamma v = v_+ - v_-$. Thus
to bound $\|M_\gamma \|$ it is enough to bound $\|v_\pm\|$ for all $v$ such that
$\|v\| = 1$. In this case, 
\begin{eqnarray*}
   \|v_\pm\| \le \frac{\|s_0(x)\|\|s'_0(x)\|}{|\langle s_0, s'_0\rangle |}
\end{eqnarray*}
Let $e_0, e_1$ be a special gauge for $E$ such that $e_0$ lies in the negative
eigenspace of $i\Phi$. We claim that for all $\epsilon > 0$,
there exists $R > 0$ and $s_0$, $s'_0$ such that if $\|x\| \ge R$ then
$\|s_0(x) - e_0(x)\| < \epsilon$ and $\|s'_0(x) - e_1(x)\| < \epsilon$. In this case,
a quick calculation shows we have
\begin{eqnarray*}
   \frac{\|s_0(x)\|\|s'_0(x)\|}{|\langle s_0, s'_0\rangle |} < \frac{(1 + \epsilon)^2}
   {1 - 2\epsilon - \epsilon^2}
\end{eqnarray*}
and so we have the required bound for $\|v_\pm\|$. It remains to prove our claim.
As we shall see, it is enough to prove it only for $s_0$ since the case for $s'_0$
follows by a similar argument.
\par
Thus let $s_0 = y_0 e_0 + y_1 e_1$ and let $t$ be the arc-length parameter for
$\gamma$ such that $x$ corresponds to $t=0$. Then the equation $(\nabla_{\dot\gamma}
- i\Phi)s_0 = 0$ becomes
\begin{eqnarray*}
   \frac{d}{dt}\left[
   \begin{array}{c}
      y_0\\
      y_1
   \end{array}\right]
   + \left[
   \begin{array}{cc}
      A(\dot\gamma) + \|\Phi\| & -B^*(\dot\gamma)\\
      B(\dot\gamma)            & -A(\dot\gamma) - \|\Phi\|
   \end{array}\right]
   \left[
   \begin{array}{c}
      y_0\\
      y_1
   \end{array}\right] = 0
\end{eqnarray*}
where
\begin{eqnarray*}
   \left[
   \begin{array}{cc}
   A    & B\\
   -B^* & -A
   \end{array}\right]
\end{eqnarray*}
is the connection matrix in this gauge.
Using $B = O(r^{-2})$ and $\|\Phi\| = 1 - k/2r + O(r^{-2})$ our equation becomes
\begin{eqnarray}\label{scatt_eqn_A}
   \frac{d}{dt}\left[
   \begin{array}{c}
      y_0\\
      y_1
   \end{array}\right]
   = \left(\Lambda_\gamma(t) + C_\gamma(t)\right)\left[
   \begin{array}{c}
      y_0\\
      y_1
   \end{array}\right]
\end{eqnarray}
where
\begin{eqnarray*}
   \Lambda_\gamma(t) = \left[
   \begin{array}{cc}
      -1 + \frac{k}{2\sqrt{t^2 + \|x\|^2}} - A(\dot\gamma) & 0\\
      0        & 1 - \frac{k}{2\sqrt{t^2 + \|x\|^2}} + A(\dot\gamma)
   \end{array}\right]
\end{eqnarray*}
and
\begin{eqnarray*}
   \|C_\gamma(t)\| = O((t^2 + \|x\|^2)^{-1})
\end{eqnarray*}
Thus not only is $\|C_\gamma(t)\|$ integrable but
$\int_0^\infty \|C_\gamma(t)\|dt \to 0$ as $\|x\| \to \infty$. Now define
\begin{eqnarray*}
   \left[
   \begin{array}{c}
      z_0(t)\\
      z_1(t)
   \end{array}\right]
   = e^{h(t)}\left[
   \begin{array}{c}
      y_0(t)\\
      y_1(t)
   \end{array}\right]
\end{eqnarray*}
where
\begin{eqnarray*}
h(t) = \int_0^t\left(1 - \frac{k}{2\sqrt{s^2 + \|x\|^2}} + A(\dot\gamma)\right)ds
\end{eqnarray*}
Then $\left[\begin{array}{c}y_0\\y_1\end{array}\right]$ solves \eqref{scatt_eqn_A}
if and only if $\left[\begin{array}{c}z_0\\z_1\end{array}\right]$ solves
\begin{eqnarray}\label{scatt_eqn_B}
   \frac{d}{dt}\left[
   \begin{array}{c}
      z_0\\
      z_1
   \end{array}\right]
   = \left(\Lambda^-_\gamma(t) + C_\gamma(t)\right)\left[
   \begin{array}{c}
      z_0\\
      z_1
   \end{array}\right]
\end{eqnarray}
where
\begin{eqnarray*}
   \Lambda^-_\gamma(t) = \left[
   \begin{array}{cc}
      0 & 0\\
      0 & 2\left(1 - \frac{k}{2\sqrt{t^2 + \|x\|^2}} + A(\dot\gamma)\right)
   \end{array}\right]
\end{eqnarray*}
Define the integral operator $T$ operating on bounded $\C^2$-valued
functions on $[0,\infty)$ by
\begin{eqnarray*}
   \left(T\left[
   \begin{array}{c}
      z_0\\
      z_1
   \end{array}\right]\right)(t)
   = - \int_t^\infty K(s,t)C_\gamma(s)\left[\begin{array}{c}
      z_0\\
      z_1
   \end{array}\right](s)ds
\end{eqnarray*}
where
\begin{eqnarray*}
   K(s, t) =
   \left[\begin{array}{cc}
      1 & 0\\
      0 & e^{2(h(t) - h(s))}
   \end{array}\right]
\end{eqnarray*}
Note that provided $\|x\| > k/2$ we have $|e^{2(h(t) - h(s))}| \le 1$ for all
$s \ge t \ge 0$ since $|e^{2(h(t) - h(s))}| = e^{2\Re(h(t) - h(s))}$ and
\begin{eqnarray*}
\Re(h(t)) = \int_0^t\left(1 - \frac{k}{2\sqrt{s^2 + \|x\|^2}}\right)ds
\end{eqnarray*}
is an increasing function if $\|x\| > k/2$. This shows that $K$ is uniformly
bounded (by a constant independent of $x$) and hence that $T$ is well defined
(ie: the integral converges).
\par
Now $\left[\begin{array}{c}z_0\\z_1\end{array}\right]$ solves \eqref{scatt_eqn_B}
and produces a solution $s_0$ with the right decay if it solves the integral equation
\begin{eqnarray}\label{z_int_eqn}
   \left[
   \begin{array}{c}
      z_0\\
      z_1
   \end{array}\right]
   = \left[
   \begin{array}{c}
      1\\
      0
   \end{array}\right]
   + T \left[
   \begin{array}{c}
      z_0\\
      z_1
   \end{array}\right]
\end{eqnarray}
But in view of our bound for $K$ and the fact that $\int_0^\infty\|C_\gamma(t)\|dt
\to 0$ as $\|x\| \to \infty$ we have $\|T\|_{L^\infty} \to 0$
as $\|x\| \to \infty$. Thus for sufficiently large $\|x\|$
\begin{eqnarray*}
   \left[
   \begin{array}{c}
      z_0\\
      z_1
   \end{array}\right]
   - \left[
   \begin{array}{c}
      1\\
      0
   \end{array}\right]
   = (T + T^2 + \cdots)\left[
   \begin{array}{c}
      1\\
      0
   \end{array}\right]
\end{eqnarray*}
and so we have
\begin{eqnarray*}
   \left\|\left[
   \begin{array}{c}
      z_0\\
      z_1
   \end{array}\right]
   - \left[
   \begin{array}{c}
      1\\
      0
   \end{array}\right]\right\|_{L^\infty}
   \le \frac{1}{\|T\|^{-1}_{L^\infty} - 1} \to 0 \qquad \mbox{as $\|x\| \to \infty$}
\end{eqnarray*}
In particular this gives
\begin{eqnarray*}
   \left\|\left[
   \begin{array}{c}
      y_0(0)\\
      y_1(0)
   \end{array}\right]
   - \left[\begin{array}{c}
      1\\
      0\end{array}\right]\right\|
   < \epsilon
\end{eqnarray*}
for all $\|x\|$ sufficiently large and so $\|s_0(x) - e_0(x)\| < \epsilon$
as required.
\myqed
\begin{remark}
   Note that in the course of the proof of lemma \ref{vanishing_lemma}
   we showed that there exists a compact subset of $\R^3$ such that for geodesics
   which do not meet this set there is an upper bound for
   \begin{eqnarray*}
      \frac{\|s_0(x)\|\|s'_0(x)\|}{|\langle s_0, s'_0\rangle |}
   \end{eqnarray*}
   In particular this shows that $\langle s_0, s'_0\rangle$ can only vanish for
   geodesics that meet this compact subset of $\R^3$, ie: all spectral
   lines meet this compact subset. This shows that the spectral curve is compact.
   Our result may thus be regarded as a slight sharpening of the result that
   the spectral curve is compact.
\end{remark}
We are ready to prove the required vanishing theorem.
\begin{theorem}\label{Euc_vanishing_theorem}
   Let $\widetilde E$ be the holomorphic vector bundle on $\T$
   corresponding to an $SU(2)$ monopole of charge $k$ on $\R^3$. Let
   $S$ be the spectral curve of the monopole. Then
   \begin{eqnarray*}
      H^0(S, \widetilde EL(k-2)) = 0
   \end{eqnarray*}
\end{theorem}
\noindent {\bf Proof} We shall adapt the methods used in \cite{MR709461}
to prove the crucial vanishing theorem $H^0(S, L^z(k-2)) = 0$ for $z \in
(0,2)$ (see also \cite{MR929140}). Thus our method is to show that we have an
injection
\begin{eqnarray*}
H^0(S, \widetilde EL(k-2)) \hookrightarrow H^1(\T, S^2\widetilde E(-2))
\end{eqnarray*}
(where $S^2\widetilde E$ denotes the symmetric square of $\widetilde E$)
and apply the Penrose transform to the resulting class to get a solution
$\phi$ of a covariant Laplace equation on $\R^3$. We will establish that
$\phi$ is decaying and so by applying a maximum principle to the function
$\|\phi\|^2$ deduce that it must in fact vanish.
\par
We begin by noting that $\widetilde EL(k-2)$ has no non-zero sections
on $\T$. This follows since we have a short exact sequence
\begin{eqnarray}\label{E_ext1}
   0 \to \cO(-2) \to \widetilde EL(k-2) \to L^2(2k-2) \to 0
\end{eqnarray}
which is obtained from the first expression of $\widetilde E$ as an extension
\begin{eqnarray*}
   0 \to L^*(-k) \to \widetilde E \to L(k) \to 0
\end{eqnarray*}
Now in \cite{MR709461} Hitchin proves $H^0(\T, L^m(p)) = 0$ for any $p
\in \Z$ and any $m\ne 0$. Thus in particular $H^0(\T, L^2(2k-2)) = 0$.
Since also $H^0(\T, \cO(-2)) = 0$ it follows from the long exact cohomology
sequence associated to \eqref{E_ext1} that we must also have
$H^0(\T, \widetilde EL(k-2)) = 0$ as claimed.
\par
Next note that since $S$ is a divisor of $\cO(2k)$ on $\T$ we have a short
exact sequence of sheaves on $\T$
\begin{eqnarray*}
   0 \to \cO_{\T}(\widetilde EL(-k-2)) \to \cO_{\T}(\widetilde EL(k-2)) \to
   \cO_S(\widetilde EL(k-2)) \to 0
\end{eqnarray*}
and since $H^0(\T, \widetilde EL(k-2)) = 0$ we have an injection
\begin{eqnarray*}
   \delta : H^0(S, \widetilde EL(k-2)) \hookrightarrow H^1(\T, \widetilde EL(-k-2))
\end{eqnarray*}
\par
However we also have an exact sequence
\begin{eqnarray}\label{E_ext2}
   0 \to \widetilde EL(-k-2) \to S^2\widetilde E(-2) \to L^{-2}(2k-2) \to 0
\end{eqnarray}
obtained from the second expression of $\widetilde E$ as an extension
\begin{eqnarray*}
   0 \to L(-k) \to \widetilde E \to L^*(k) \to 0
\end{eqnarray*}
Since $H^0(\T, L^{-2}(2k-2)) = 0$ the long exact cohomology sequence associated
to \eqref{E_ext2} supplies us with another injection
\begin{eqnarray*}
   i : H^1(\T, \widetilde EL(-k-2)) \hookrightarrow H^1(\T, S^2 \widetilde E(-2))
\end{eqnarray*}
We thus have our injection $i\delta : H^0(S, \widetilde EL(k-2)) \hookrightarrow 
H^1(\T, S^2\widetilde E(-2))$. As in \cite{MR709461}, we will represent a class
in the image of $i\delta$ in Dolbeault cohomology by differential forms with
specific decay and support properties.
\par
Specifically, we claim that if $s \in H^0(S, \widetilde EL(k-2))$, then we can
find a Dolbeault representative $\theta^+$ for $\delta(s)$ such that
\begin{enumerate}
   \item
   $\theta^+$ has compact support
   \item
   $\theta^+ \in \Omega^{0,1}(\T, (L^+)^2(-2)) \subset \Omega^{0,1}(\T, \widetilde
   EL^+(-2)) = \Omega^{0,1}(\T, \widetilde EL(-k-2))$
\end{enumerate}
where $L^+ \subset \widetilde E$ is the line subbundle $\simeq L(-k)$ of
$\widetilde E$ introduced in \cite{MR649818} (whose fibre consists of solutions
of the scattering equation $\nabla - i\Phi = 0$ decaying at the positive end of a
geodesic).
\par
The crucial step in establishing this claim is to note that in view of the exact
sequence
\begin{eqnarray*}
   0 \to L^2(-2) \to \widetilde EL(k-2) \to \cO(2k-2) \to 0
\end{eqnarray*}
we have an injection $H^0(S, \widetilde EL(k-2)) \hookrightarrow H^0(S,
\cO(2k-2))$. Furthermore the restriction map $H^0(\T, \cO(2k-2)) \to H^0(S, \cO(2k-2))$
is an isomorphism so $s$ defines an element
\begin{eqnarray*}
   t \in H^0(\T, \cO(2k-2))
\end{eqnarray*}
We will use this polynomial $t$ when constructing our Dolbeault representative $\theta^+$
to ensure that it has the right properties.
\par
We now recall the explicit description of $\widetilde E$ used in \cite{MR709461}.
Thus define
\begin{eqnarray*}
   E^+ = L(-k)\oplus L^*(k)
\end{eqnarray*}
and define a $\deebar$-operator on $E^+$ by
\begin{eqnarray*}
   \deebar = \left[
   \begin{array}{cc}
      \deebar & \deebar\alpha/\psi\\
      0 & \deebar
   \end{array}\right]
\end{eqnarray*}
where $\alpha$ is a $C^\infty$ section
of $L^2$ on $\T$ supported in a compact neighbourhood of $S$ that restricts
to a holomorphic trivialisation of $L^2$ on $S$ (satisfying $\alpha\alpha^* = -1$
on $S$) and $\psi \in H^0(\T, \cO(2k))$ is a section defining the spectral curve
$S$. With this $\deebar$-operator on $E^+$ we have $E^+ \simeq \widetilde E$ and we
can explicitly see the subbundle $L^+ \subset \widetilde E$ as $L(-k) \subset E^+$.
\par
Now let $\{V_i\}$ be an open cover of $S$ by sufficiently small open balls
in $\T$ and such that the cover lies inside a compact subset $K$ of $\T$. Extend
this to an open cover of $K$ by adding in sufficiently small open balls $\{W_i\}$
in $\T$ such that the new cover lies inside a compact subset $K^{'}$ of $\T$.
Finally, extend this to an open cover $\cU$ of $\T$ in such a way that no open
set in $\cU\setminus (\{V_i\}\cup \{W_i\})$ meets $K$.
\par
Note \cite{MR709461} that we may use the cover $\{V_i\}$ to compute $\alpha$
and so we may assume $\supp(\deebar\alpha/\psi) \subseteq K$.
\par
Using $\widetilde E \simeq E^+$, over each open set $U\in \cU$ we have
\begin{eqnarray*}
   s|_{S\cap U} = \left[
   \begin{array}{l}
      s_1^U\\
      s_2^U
   \end{array}
   \right]
\end{eqnarray*}
with $s_1^U \in \Omega^0(U\cap S, L^2(-2))$ and $s_2^U \in \Omega^0(U\cap S,
\cO(2k-2))$. For each $U$, we choose a holomorphic extension
\begin{eqnarray*}
   \sigma = \left[
   \begin{array}{l}
      \sigma_1^U\\
      \sigma_2^U
   \end{array}
   \right]
\end{eqnarray*}
of $s|_{S\cap U}$ to all of $U$ such that $\sigma_2^U = t|_{U}$ and such that
for $U \in \cU\setminus (\{V_i\}\cup \{W_i\})$ we choose $\sigma_1^U = 0$. We
can do this because on such $U$ we have $E^+L(k-2) = L^2(-2)\oplus \cO(2k-2)$
holomorphically since $\supp(\deebar\alpha/\psi) \subset K$ and $U\cap K = \emptyset$.
\par
Now if $\{\phi_U\}$ is a partition of unity subordinate to $\cU$, a Dolbeault
representative for $\delta(s)$ is
\begin{eqnarray*}
   \theta^+ &=& \frac{1}{\psi}\deebar\sum_{U\in\cU}\limits\phi_U\left[
   \begin{array}{l}
      \sigma_1^U\\
      \sigma_2^U
   \end{array}
   \right]
   = \frac{1}{\psi}\deebar\left[
   \begin{array}{l}
      \sigma\\
      t
   \end{array}
   \right]
\end{eqnarray*}
where $\sigma = \sum_{U\in\cU}\limits \phi_U\sigma_1^U \in \Omega^0(T, L^2(-2))$.
Thus since $\deebar t = 0$ we have
\begin{eqnarray*}
   \theta^+ = 
   \frac{1}{\psi}\left[
   \begin{array}{c}
      \deebar\sigma + \frac{\deebar\alpha}{\psi} t\\
      0
   \end{array}
   \right]
\end{eqnarray*}
Evidently $\theta^+$ is supported in $K^{'}$ and takes values in $L^2(-2k-2) \simeq
(L^+)^2(-2)$ and so our claim is true.
\par
Similarly the same claim is true with $L^-$ in place of $L^+$ and so we also
have $\theta^- \in \Omega^{0,1}(\T, (L^-)^2(-2))$ with compact support. Using
$(L^{\pm})^2 \hookrightarrow S^2 \widetilde E$ both $\theta^+$ and
$\theta^-$ represent $i\delta(s)$. We thus have $\theta^+ - \theta^- =
\deebar\gamma$ for some $\gamma \in \Omega^0(\T, S^2 \widetilde E(-2))$.
In fact if $E^- = L^*(-k)\oplus L(k)$ with $\deebar$-operator
\begin{eqnarray*}
   \deebar = \left[
   \begin{array}{cc}
      \deebar & \deebar\alpha^*/\psi\\
      0 & \deebar
   \end{array}\right]
\end{eqnarray*}
then using
\begin{eqnarray}\label{h_Euc_matrix}
   h = \left[\begin{array}{cc}
      -\alpha^* & -\frac{\alpha\alpha^* + 1}{\psi}\\
      \psi & \alpha
   \end{array}\right]
   : E^+ \simeq E^- 
\end{eqnarray}
as in \cite{MR709461}, we can choose
\begin{eqnarray*}
   \theta^- = \frac{1}{\psi}\left[
   \begin{array}{c}
   \deebar\left(\frac{\alpha\alpha^* + 1}{\psi}\alpha^* t + (\alpha^*)^2
   \sigma\right) + \frac{\deebar\alpha^*}{\psi}t\\
   0
   \end{array}\right]
\end{eqnarray*}
Now using the smooth isomorphism $S^2 E^- \simeq L^{-2}(-2k)
\oplus \cO \oplus L^2(2k)$ we can get an explicit formula for $\gamma \in
\Omega^0(\T, S^2 E^-(-2))$. Indeed we have
\begin{eqnarray*}
   && S^2 h\left[
   \begin{array}{c}
      \theta^+\\
      0\\
   \end{array}\right]
   - \left[
   \begin{array}{c}
      \theta^-\\
      0\\
   \end{array}\right]\\
   &=& \left[
   \begin{array}{ccc}
      (\alpha^*)^2   & \alpha^*\frac{\alpha\alpha^* + 1}{\psi} & \left(
      \frac{\alpha\alpha^* + 1}{\psi}\right)^2\\
      -2\alpha^*\psi & -2\alpha\alpha^* - 1                    & -2\alpha
      \frac{\alpha\alpha^* + 1}{\psi}\\
      \psi^2         & \psi\alpha                              & \alpha^2
   \end{array}\right]
   \frac{1}{\psi}\left[
   \begin{array}{c}
      \deebar\sigma + \frac{\deebar\alpha}{\psi}t\\
      0\\
      0
   \end{array}\right]\\
   && - \frac{1}{\psi}\left[
   \begin{array}{c}
      \deebar\left(\frac{\alpha\alpha^* + 1}{\psi}\alpha^*t + (\alpha^*)^2\sigma\right)
      + \frac{\deebar\alpha^*}{\psi}t\\
      0\\
      0
   \end{array}\right]\\
   &=& \left[
   \begin{array}{c}
      \frac{(\alpha^*)^2}{\psi}\left(\deebar\sigma + \frac{\deebar\alpha}{\psi}
      t\right) - \frac{1}{\psi}\left\{\deebar\left(\frac{\alpha\alpha^* + 1}{\psi}
      \alpha^* t + (\alpha^*)^2\sigma\right) + \frac{\deebar\alpha^*}{\psi}t\right\}\\
      -2\alpha^*\left(\deebar\sigma + \frac{\deebar\alpha}{\psi}t\right)\\
      \psi\left(\deebar\sigma + \frac{\deebar\alpha}{\psi}t\right)
   \end{array}\right]\\
   &=& \left[
   \begin{array}{ccc}
      \deebar & \frac{\deebar\alpha^*}{\psi} & 0\\
      0 & \deebar & 2\frac{\deebar\alpha^*}{\psi}\\
      0 & 0 & \deebar
   \end{array}\right]
   \left[
   \begin{array}{c}
      0\\
      -2\left(\alpha^*\sigma + \frac{\alpha\alpha^* + 1}{\psi}t\right)\\
      \psi\sigma + \alpha t
   \end{array}\right]
\end{eqnarray*}
We thus have a decomposition of $\gamma$ as a sum of two terms $\gamma =
\gamma_a + \gamma_b$ where
\begin{eqnarray}\label{gamma_formula}
   \gamma_a &=& \left[
   \begin{array}{c}
      0\\
      -2\alpha^*\sigma\\
      \psi\sigma + \alpha t
   \end{array}\right]\notag\\
   \gamma_b &=& \left[
   \begin{array}{c}
   0\\
   -2\frac{\alpha\alpha^* + 1}{\psi}t\\
   0
   \end{array}\right]
\end{eqnarray}
Note that $\gamma_a$ is compactly supported but $\gamma_b$ is not.
\par
Now let $\phi = P(i\delta(s))$ where
\begin{eqnarray*}
   P : H^1(\T, S^2 \widetilde E(-2)) \to \left\{ \phi \in \Omega^0(
   \R^3, S^2 E)\quad | \quad \Delta_{(A, \Phi)}\phi = 0\right\}
\end{eqnarray*}
is the Penrose transform (see eg: \cite{MR1054377}) and
\begin{eqnarray*}
   \Delta_{(A, \Phi)} = \nabla_A^*\nabla_A + \Phi^*\Phi
\end{eqnarray*}
Fix $O \in \R^3$ and let $r$ be the radial distance in $\R^3$ to $O$.
We claim that
\begin{eqnarray*}
   \|\phi\| = O(r^{-1})
\end{eqnarray*}
To see this, let $x \in \R^3$ and let
$P_x \subset \T$ be the corresponding twistor line. Let $l^+_x$ be the
oriented line $xO$ and $l^-_x$ be the oriented line $Ox$.
Now, just as in \cite{MR709461}, since $\theta^\pm$
are compactly supported in a neighbourhood of $S$, we can find
disjoint neighbourhoods $V^\pm_x$ of $l^\pm_x$ in $P_x$ such that
$\supp\theta^\pm|_{P_x} \subset V^+_x \cup V^-_x$.
Let $v_1, v_2$ be an $SU(2)$ basis
for $E_x$ and let $f_1, f_2 \in H^0(P_x, \widetilde E)$ be the corresponding
sections. Then, using $S^2 E \simeq \End_0(E)$ and using
the recipe for the Penrose transform we have
\begin{eqnarray}\label{Euc_phi_formula}
   \phi_{ij}(x) &=& \langle v_i, \phi(x)v_j\rangle\notag\\
    &=& \int_{P_x} \langle f_i, \theta^+ f_j\rangle
    = \int_{V^+_x} \langle f_i, \theta^+ f_j\rangle +
   \int_{V^-_x} \langle f_i, \theta^+ f_j \rangle\notag\\
   &=& \int_{V^+_x} \langle f_i, \theta^+ f_j\rangle + \int_{V^-_x} \langle f_i,
   \theta^- f_j \rangle + \int_{V^-_x} \langle f_i, \deebar \gamma f_j\rangle   
\end{eqnarray}
where $\langle,\rangle$ is the skew form. Now of the three terms above, Hitchin's
arguments in \cite{MR709461} show that
the first two $\int_{V^\pm_x} \langle f_i, \theta^\pm f_j\rangle$ have exponential
decay in $x$ since $\theta^\pm$ takes values in $(L^\pm)^2(-2)$. Thus we need only
show that the last term decays as $O(\|x\|^{-1})$. But if $H^-_x$ is the hemisphere
in $P_x$ containing $V^-_x$ with boundary consisting of those lines perpendicular
to the line joining $O$ and $x$ then using the fact that $\deebar\gamma|_{P_x}$ is
supported in the disjoint subsets $V^\pm_x$ we have
\begin{eqnarray*}
   \int_{V^-_x} \langle f_i, \deebar \gamma f_j\rangle &=& \int_{H^-_x} \langle
   f_i, \deebar\gamma f_j\rangle = \int_{H^-_x} \deebar \langle f_i, \gamma f_j
   \rangle = \int_{\partial H^-_x} \langle f_i, \gamma f_j\rangle = 
   \int_{\partial H^-_x} \langle f_i, \gamma_b f_j\rangle
\end{eqnarray*}
since $\gamma_a$ is supported in $V^\pm_x$ which doesn't meet $\partial H^-_x$.
Now under our identification $S^2 \widetilde E \simeq \End_0(\widetilde E)$
and using $\widetilde E \simeq E^- = L^*(-k) \oplus L(k)$ we have from
\eqref{gamma_formula}
\begin{eqnarray*}
   \gamma_b = \left[
   \begin{array}{cc}
      -2t/\psi & 0\\
      0        & 2t/\psi
   \end{array}\right]
\end{eqnarray*}
on $\partial H^-_x$ since $\alpha$ vanishes there.
\par
Consider now $2t/\psi$. Introduce the usual coordinates $\zeta, \eta$ on $\T$
with $O$ corresponding to the zero section $\eta = 0$. Then it is elementary to
check that on $\partial H^-_x$ we have
\begin{eqnarray*}
   \|x\| = \frac{|\eta|}{1 + |\zeta|^2}
\end{eqnarray*}
Now $t \in H^0(\T, \cO(2k - 2))$ and so let
\begin{eqnarray*}
   t = a_0 \eta^{k-1} + a_1(\zeta) \eta^{k-2} + \cdots + a_{k-1}(\zeta)
\end{eqnarray*}
where $a_i$ has degree $2i$ in $\zeta$. Thus on $\partial H^-_x$
\begin{eqnarray*}
   |t| &\le& |a_0|(1 + |\zeta |^2)^{k-1}\|x\|^{k-1} + |a_1(\zeta)|
   (1 + |\zeta |^2)^{k-2}\|x\|^{k-2} + \cdots + |a_{k-1}(\zeta)|\\
   &=& Q^t_\zeta(\|x\|)
\end{eqnarray*}
for a polynomial $Q^t_\zeta$. Thus provided $\|x\| \ge 1$ we have
\begin{align*}
   \left|\frac{Q^t_\zeta(\|x\|) - |a_0|(1 + |\zeta|^2)^{k-1}\|x\|^{k-1}}{(1 + |\zeta|)^{k-1}
   \|x\|^{k-1}}\right| &= \frac{1}{\|x\|}\left(\frac{a_1(\zeta)}{1 + |\zeta|^2} +
   \cdots + \frac{1}{\|x\|^{k-2}}\frac{|a_{k-1}(\zeta)|}{(1 + |\zeta|^2)^{k-1}}\right)\\
   &\le \frac{1}{\|x\|}g(\zeta)
\end{align*}
where
\begin{eqnarray*}
   g(\zeta) = \frac{|a_1(\zeta)|}{1 + |\zeta|^2} + \cdots + \frac{|a_{k-1}(\zeta)|}
   {(1 + |\zeta|^2)^{k-1}}
\end{eqnarray*}
Now $g$ is bounded since it is a sum of ratios of polynomials of the same degree.
Thus if we choose $R^t$ large enough then if $\|x\| > R^t$ we have
\begin{eqnarray*}
   |t| &\le& \left|Q^t_\zeta(\|x\|) - |a_0|(1 + |\zeta|^2)^{k-1}\|x\|^{k-1}\right| +
   |a_0|(1 + |\zeta|^2)^{k-1}\|x\|^{k-1}\\
   &\le& (1 + |a_0|)\|x\|^{k-1}(1 + |\zeta|^2)^{k-1}
\end{eqnarray*}
on $\partial H^-_x$. By a similar argument, we can choose $R^\psi$ large enough such
that if $\|x\| > R^\psi$ then
\begin{eqnarray*}
   |\psi| \ge \frac{1}{2}(1 + |\zeta|^2)^k \|x\|^k
\end{eqnarray*}
on $\partial H^-_x$. Thus there exists $R$ such that if $\|x\| > R$ then
\begin{eqnarray*}
   \left|\frac{t}{\psi}\right| &\le& \frac{(1 + |a_0|)(1 + |\zeta|^2)^{k-1}\|x\|^{k-1}}
   {\frac{1}{2}(1 + |\zeta|^2)^k\|x\|^k}\\
   &=& \frac{2(1 + |a_0|)}{(1 + |\zeta|^2)\|x\|} \le \frac{2(1 + |a_0|)}{\|x\|}
\end{eqnarray*}
on $\partial H^-_x$.
\par
In view of this, if we can show that the function $\langle f_i, \left[\begin{array}{cc}
1 & 0\\0 & -1\end{array}\right]f_j\rangle$ is bounded for sufficiently large $\|x\|$
(by a constant independent of $x$) on $\partial H^-_x$ then it will follow that $\phi$
has the required decay. However the required boundedness is exactly what is proved
in lemma \ref{vanishing_lemma}. We thus have $\|\phi\| = O(r^{-1})$ as claimed.
\par
Finally note that since for any section $\phi$
\begin{eqnarray*}
   \frac{1}{2}\Delta\|\phi\|^2 = \|\nabla\phi\|^2 - (\nabla^*\nabla\phi, \phi)
\end{eqnarray*}
(where $\Delta$ is just the Laplacian for functions on $\R^3$) and since for our
section $\phi$
\begin{eqnarray*}
   \left(\nabla^*\nabla + \Phi^*\Phi\right)\phi = 0
\end{eqnarray*}
we have
\begin{eqnarray}\label{Euc_maxprin_ineq}
   \frac{1}{2}\Delta\|\phi\|^2 = \|\nabla\phi\|^2 + \|\Phi(\phi)\|^2 \ge 0
\end{eqnarray}
So $\|\phi\|^2$ is subharmonic. Since we also have $\|\phi\|^2 = O(r^{-2})$ the
maximum principle for subharmonic functions on $\R^3$ (see eg \cite{MR1814364})
requires that we in fact have $\phi = 0$ as required.
\myqed
\par
Having established the necessary vanishing theorem, we are
ready to proceed with the deformation theory. Firstly we have
\begin{corollary}\label{vanish_cor_dims}
   In the notation of corollary \ref{Euc_norm_bnd_id} we have
   \begin{enumerate}
      \item
      $h^1(S, \hat N) = 0$
      \item
      $h^0(S, \hat N) = 4k$
      \item\label{vanish_cor_dim_3}
      $h^0(S, \hat N(-1)) = 2k$
   \end{enumerate}
   (where we have identified $\hat S$ with $S$)
\end{corollary}
\noindent {\bf Proof}
Firstly, note that in view of corollary \ref{Euc_norm_bnd_id} we can restate
theorem \ref{Euc_vanishing_theorem} as
\begin{eqnarray*}
   H^0(S, \hat N(-2)) = 0
\end{eqnarray*}
Next note that in view of lemma \ref{NormExtLemma} we have the short exact sequence
\begin{eqnarray}\label{Euc_nrm_bnd_ext}
   0 \to \cO \to \hat N \to \cO(2k) \to 0
\end{eqnarray}
We thus have $\wedge^2 \hat N \simeq \cO(2k)$ and so
\begin{eqnarray*}
   \hat N^*(2k) \simeq \hat N
\end{eqnarray*}
Furthermore, this sequence also allows us to calculate the Chern
class of $\hat N$. Indeed, recalling that $S$ is a k-fold (branched) covering
of $\CP$, we have
\begin{eqnarray*}
   c_1(\hat N)[S] = 2k^2
\end{eqnarray*}
Finally recall that since the spectral curve is a divisor of $\cO(2k)$
on $\T$ its canonical bundle is (by adjunction) $\cO(2k-4)$.
\begin{enumerate}
   \item
   By Serre duality and our observations above
   \begin{eqnarray*}
      h^1(S, \hat N) = h^0(S, \hat N^*(2k-4)) = h^0(S, \hat N(-4)) = 0
   \end{eqnarray*}
   since $h^0(S, \hat N(-2)) = 0$ (and $h^0(S, \cO(2)) > 0$).
   \item
   By the Riemann-Roch formula
   \begin{eqnarray*}
      h^0(S, \hat N) - h^1(S, \hat N) &=& (2 + c_1(\hat N))(1 + \frac{1}{2}
      c_1(TS))[S]\\
      &=& 4k
   \end{eqnarray*}
   and $h^1(S, \hat N) = 0$ by part (i).
   \item
   By the Riemann-Roch formula
   \begin{eqnarray*}
      h^0(S, \hat N(-1)) - h^1(S, \hat N(-1)) = 2k
   \end{eqnarray*}
   We thus need only show $h^1(S, \hat N(-1)) = 0$ but this follows
   immediately since
   \begin{eqnarray*}
      h^1(S, \hat N(-1)) = h^0(S, \hat N^*(2k-3)) = h^0(S, \hat N(-3)) = 0
   \end{eqnarray*}
   since $h^0(S, \hat N(-2)) = 0$ (and $h^0(S, \cO(1)) > 0$). \myqed
\end{enumerate}
\par
The results of this corollary are quite significant. Not only have we
established the necessary vanishing of $H^1(S, \hat N)$ so that we can apply
Kodaira's theorem but, as claimed in section \ref{DefSpecOverview}, we find
that the space of deformations has dimension $4k$ without appealing to the
analytical results of Taubes \cite{MR723549}.
\begin{remark}
   \par\noindent
   \begin{enumerate}
      \item
      In the proof of the vanishing theorem, it was clear a priori that we could
      obtain (compactly supported) representatives $\theta^{\pm}$ taking values
      in $\widetilde E L^{\pm}(-2)$. It was important that we did better and found
      (compactly supported) representatives taking values in $(L^{\pm})^2(-2)$.
      Otherwise after the Penrose transform we will only obtain a bounded
      (rather than decaying) solution of the covariant Laplace equation
      $\Delta_{(A,\Phi)}\phi = 0$ on $\R^3$. Since the Higgs field $\Phi$ is
      bounded and (after identifying $\End_0(E) \simeq S^2 E$) satisfies
      $\Delta_{(A, \Phi)}\Phi = 0$ it is clear that boundedness is not enough
      to deduce vanishing. (For another consequence of the observation
      $\Delta_{(A, \Phi)}\Phi = 0$ see subsection \ref{atiyah_higgs_subsect}.)
      \item
      Another point about the vanishing theorem is that, unlike Hitchin's proof
      that $H^0(S, L^z(k-2)) = 0$ for $z \in (0,2)$, our proof does not generalise
      from the case $H^0(S, \widetilde EL(k-2)) = 0$ to a possible vanishing result
      for $H^0(S, \widetilde EL^z(k-2))$. This is because we used the injection
      \begin{eqnarray*}
         H^0(S, \widetilde EL(k-2)) \hookrightarrow H^0(S,\cO(2k-2))
      \end{eqnarray*}
      to obtain a class in $H^0(S, \cO(2k-2))$ and then noted that the restriction
      map $H^0(\T, \cO(2k-2)) \to H^0(S, \cO(2k-2))$ is an isomorphism. If we replace
      $L$ with $L^z$, we get an injection
      \begin{eqnarray*}
         H^0(S, \widetilde EL^z(k-2)) \hookrightarrow H^0(S,L^{z-1}(2k-2))
      \end{eqnarray*}
      but we can no longer extend a class in $H^0(S,L^{z-1}(2k-2))$ to all of 
      $\T$ unless $z=1$. Indeed $H^0(\T,L^{z-1}(2k-2))$ vanishes unless $z=1$.
   \end{enumerate}
\end{remark}
\par
Now that we have established the natural isomorphism \eqref{Euc_tang_spc_model}
we have a good model for the tangent space to the moduli space $M_k$, ie: a
natural isomorphism
\begin{eqnarray*}
   T_{\hat S}M_k \otimes_{\R}\C \simeq H^0(S, \hat N)
\end{eqnarray*}
Using this, we can understand the geometry of the moduli space. We begin by
noting that we have a natural isomorphism
\begin{eqnarray*}
   H^0(S, \cO(1)) \otimes H^0(S, \hat N(-1)) \simeq H^0(S, \hat N)
\end{eqnarray*}
To see this, consider the extension on $\CP$
\begin{eqnarray*}
   0 \to \cO(-1) \to \C^2 \to \cO(1) \to 0
\end{eqnarray*}
Pull this back to $S$ and tensor it with $\hat N(-1)$. Taking the long exact
sequence of cohomology groups and using \ref{Euc_vanishing_theorem} to see
$H^0(S, \hat N(-2)) =  H^1(S, \hat N(-2)) = 0$ we obtain the required
isomorphism.
\par
We thus have a natural isomorphism
\begin{eqnarray}\label{EucModTanDecom}
   T_{\hat S}M_k \otimes_{\R} \C \simeq H^0(S, \hat N(-1))
   \otimes_{\C} H^0(S, \cO(1))
\end{eqnarray}
This is important because of the following
\begin{proposition}\label{HK_EH_prop}
   Let $M$ be a smooth manifold of dimension $4k$. Let $E$ be a rank $2k$ smooth
   complex vector bundle on $M$ carrying a quaternionic structure (ie: an anti-
   linear endomorphism whose square is $-1$) and compatible symplectic structure.
   Let $H$ be a trivial rank $2$ smooth complex vector bundle on $M$ carrying a
   quaternionic structure. Then an isomorphism
   \begin{eqnarray*}
      TM \otimes_{R}\C \simeq E\otimes_{\C}H
   \end{eqnarray*}
   respecting real structures naturally induces a reduction of the structure group
   of $TM$ to $Sp(k)$.
\end{proposition}
\noindent {\bf Proof}
   See \cite{MR664330}. \myqed
\par
Immediately we notice that since $H^0(S, \cO(1)) \simeq H^0(\CP, \cO(1))$, the
rank $2$ bundle with fibre $H^0(S, \cO(1))$ naturally trivialises and carries a
natural quaternionic structure. Furthermore we have
\begin{lemma}\label{EucModSkewForm}
   $H^0(S, \hat N(-1))$ carries a natural quaternionic structure and
   a compatible non-degenerate complex skew-symmetric pairing.
\end{lemma}
\noindent {\bf Proof}
   The real structure induces the quaternionic structure.
   To define the skew-symmetric pairing $\omega$, say, we begin by noting
   that in view of the natural exact sequences
   \begin{eqnarray*}
      0 \to T\hat S \to TL^2|_{\hat S} \to \hat N \to 0
   \end{eqnarray*}
   and
   \begin{eqnarray*}
      0 \to \cO \to T(L^2\setminus 0) \to T\T \to 0
   \end{eqnarray*}
   we have a natural isomorphism
   \begin{eqnarray*}
      \wedge^2\hat N \simeq K_S\otimes K_{\T}^*|_{S}
   \end{eqnarray*}
   where $K_S$ and $K_{\T}$ are the canonical bundles of $S$ and $\T$ respectively.
   Note also that we have a natural isomorphism $K_{\T} \simeq \cO(-4)$.
   Now let $\lambda \in H^1(S, \cO(-2)) \simeq H^1(\CP, \cO(-2))$ be the 
   tautological element. Using the above, we define $\omega$ via the following
   bilinear map 
   \begin{eqnarray*}
      H^0(S, \hat N(-1))\times H^0(S, \hat N(-1)) &\to& H^0(S, K_S\otimes K_T^*(-2))\\
      &\simeq& H^0(S, K_S(2))\\
      &\overset{\lambda}\to& H^1(S, K_S)\\
      &\simeq& \C
   \end{eqnarray*}
   To see that $\omega$ is non-degenerate, note that $\lambda$ is the extension
   class of
   \begin{eqnarray}\label{Euc_lambda_ext_def}
      0 \to \cO(-1) \to \C^2 \to \cO(1) \to 0
   \end{eqnarray}
   Now $\omega$ defines a map
   \begin{eqnarray*}
      \hat\omega : H^0(S, \hat N(-1)) &\to& H^0(S, \hat N(-1))^*\\
      &\simeq& H^1(S, K_S\otimes \hat N^*(1))\\
      &\simeq& H^1(S, \hat N\otimes K_{\T}(1))\\
      &\simeq& H^1(S, \hat N(-3))
   \end{eqnarray*}
   and the statement that $\omega$ is non-degenerate is the same as saying
   $\Ker\hat\omega = 0$. However this map can also
   be recognised as the connecting homomorphism of the long exact sequence
   of cohomology groups associated to the short exact sequence of vector bundles
   \begin{eqnarray*}
      0 \to \hat N(-3) \to \C^2\otimes \hat N(-2) \to \hat N(-1) \to 0
   \end{eqnarray*}
   formed by tensoring the sequence \eqref{Euc_lambda_ext_def} with $\hat N(-2)$.
   Thus
   \begin{eqnarray*}
      \Ker \hat\omega = H^0(S, \C^2\otimes \hat N(-2))
   \end{eqnarray*}
   which vanishes by theorem \ref{Euc_vanishing_theorem} and corollary
   \ref{Euc_norm_bnd_id}. \myqed
\par
We thus have a natural reduction of the structure group of $TM_k$ to $Sp(k)$.
To complete the story from this point of view we need to show that this in
fact defines a hyperk\"ahler structure. One way to accomplish this would be
to demonstrate the existence of a torsion free connection compatible with the
reduction. To do this, a theory of connections of Kodaira deformation spaces
would be necessary and indeed some steps have been taken in this direction
(see eg: \cite{MR1327027}). Although this would be a satisfactory resolution
of the issue, it is easier to resort to the more elementary approach detailed below.
Furthermore, from the coordinate expressions below we will be able to see that we
have recovered the usual hyperk\"ahler structure on the monopole moduli space.
\par
First, let us recall some general facts about hyperk\"ahler geometry.
A hyperk\"ahler structure on a manifold $M$ includes a hypercomplex structure.
Indeed if we have $TM\otimes_{\R}\C \simeq E\otimes_{\C}H$ as in proposition
\ref{HK_EH_prop} then the complex structures on $M$ that make up the hypercomplex
structure are naturally parameterised by the projectivisation $\mathbb{P}(H)$. If we
pick a non-zero vector $h\in H$ then the image of $E$ under the injection
\begin{eqnarray*}
   E &\hookrightarrow& TM\otimes_{\R}\C\\
   e &\mapsto& e\otimes h
\end{eqnarray*}
depends only on the point $h$ defines in $\mathbb{P}(H)$ and is the bundle of $(1,0)$
vectors for a complex structure on $M$.
\par
In our case, $H = H^0(\CP, \cO(1))$ and the fibre of the bundle $E$ at a point $S$ is
$E_S = H^0(S, \hat N(-1))$. Thus the lines $[v] \in H^0(\CP, \cO(1))$ parameterise
almost complex structures on the monopole moduli space $M_k$. A first step towards
proving that our reduction of structure group is integrable is to prove that these
complex structures are integrable. To do this we need to show that the
bundle $T^{1,0}_{[v]}$ of $(1,0)$ forms for this almost complex structure is Frobenius
integrable. Since the space of real spectral curves $M_k$ is a real submanifold of the
complex space $M_k^\C$ of all complex deformations spectral curves, we naturally have
\begin{eqnarray*}
   T^{1,0}_{[v]} \subset TM_k^\C
\end{eqnarray*}
Thus to show that this bundle is Frobenius integrable, it
is enough to exhibit a complex submanifold of $M_k^\C$ with $T^{1,0}_{[v]}$ as tangent
space. Now $v$ vanishes at a point $p \in \CP$. Pick a point $S \in M_k$ and let
\begin{eqnarray*}
D = \{p_1, \ldots, p_k\}
\end{eqnarray*}
be the divisor consisting of the $k$ points in $S$ lying
over $p$. For simplicity we assume $p$ is not a branch point of the map $S \to \CP$.
The section $v$ on $\CP$ pulls back to a section of $\cO(1)$ on $S$ with divisor $D$.
As we will see, the fibre of $T^{1,0}_{[v]}$ at $S$ consists of those sections
in $H^0(S, \hat N) \simeq TM_k^\C$ vanishing at all the $p_i$. Now $p_i$
are the points where $S$ intersects the fibre of $L^2\setminus 0$ over $p \in \CP$ and
the sections vanishing at $p_i$ correspond to tangents to the submanifold $M_k^\C(D)$
of those deformations of $S$ that pass
through the same $k$ points $p_1, \ldots, p_k$ in the fibre of $L^2\setminus 0$ over $p$.
Indeed, consider the fibre $(L^2\setminus 0)_p \simeq \C\times\C^*$ of $L^2\setminus 0$
over $p$ and the map
\begin{eqnarray}\label{spec_fib_int}
   M_k^\C &\to& s^k(L^2\setminus 0)_p\\
   S &\mapsto& S \cap (L^2\setminus 0)_p\notag
\end{eqnarray}
($s^k$ is the symmetric product operation).
We will work in a neighbourhood of $S$ in $M_k^\C$ so that
the image of \eqref{spec_fib_int} avoids the singularities in the symmetric product
(since $p$ is not a branch point of $S \to \CP$).
The fibre of this map containing $S$ is our manifold $M_k^\C(D)$.
Furthermore the derivative of this map at $S$ is just the restriction map
\begin{eqnarray}\label{spec_fib_int_deriv}
   H^0(S, \hat N) \to H^0(D, \hat N)
\end{eqnarray}
Thinking about this map in terms of the short exact sequence of sheaves on $S$
\begin{eqnarray*}
   0 \to \hat N(-1) \overset{v}\to \hat N \to \hat N_D \to 0
\end{eqnarray*}
and recalling (corollary \ref{vanish_cor_dims}) that we have $h^0(S, \hat N(-1)) = 2k$,
$h^0(S, \hat N) = 4k$ and $h^0(S, \hat N_D) = 2k$ we see that \eqref{spec_fib_int_deriv}
in fact fits into the exact sequence
\begin{eqnarray}
   0 \to H^0(S, \hat N(-1)) \overset{v}\to H^0(S, \hat N) \to H^0(D, \hat N) \to 0
\end{eqnarray}
The map \eqref{spec_fib_int_deriv} is thus surjective which tells us that a generic
set of $k$ points in $(L^2\setminus 0)_p$ is the intersection with a curve in $M_k^\C$
and that the fibre  $M_k^\C(D)$ is a manifold of the correct dimension and so we have
integrability of the complex structure defined by $[v]$.
\par
Now, up to a factor of $\C^*$, $M_k^\C(D)$ carries a natural non-degenerate skew form
$\omega_D$ defined by using the isomorphism $T_S M_k^\C(D) = T^{1,0}_{[v]} \simeq
H^0(S, \hat N(-1))$ above and then using the skew form we defined on $H^0(S, \hat N(-1))$
in lemma \ref{EucModSkewForm}. Unraveling our definition of the skew form on
$H^0(S, \hat N(-1))$, we see that the recipe for $\omega_D$ is as follows. Let
\begin{eqnarray*}
   \rho \in H^0(L^2\setminus 0, K_{L^2\setminus 0}(4))
\end{eqnarray*}
be the natural non-vanishing section, where $K_{L^2\setminus 0}$ denotes the canonical
bundle of $L^2\setminus 0$. (From one point of view, $\rho$ exists because
$L^2\setminus 0$ is the twistor space of $\R^3\times S^1$ which carries a hyperk\"ahler
metric.) Now pick two vectors $X_1$, $X_2$ in the tangent space $T^{1,0}_{[v]}$ of
$M_k^\C(D)$. These are sections of $\hat N$ on $S$ vanishing at the points of $D$.
From these we obtain sections $X_i/v$ of $\hat N(-1)$ on $S$. If we pair these with
$\rho$ we obtain
\begin{eqnarray*}
   \rho(X_1/v, X_2/v) \in H^0(S, K_S(2))
\end{eqnarray*}
Finally we multiply this section with the pull back of the canonical element of
$H^1(\CP, \cO(-2))$ and apply Serre duality (integrate) to get a complex number.
This is the value of $\omega_D(X_1, X_2)$. To complete our discussion of integrability
we now show that $\omega_D$ is a holomorphic symplectic form on $M_k^\C(D)$.
\par
We can obtain an even more explicit description of $\omega_D$ (from which it will be
obvious that it is holomorphic) if we introduce some natural local coordinates on
$M_k^\C(D)$. We parameterise curves in $M_k^\C(D)$ by their points of intersection
with the fibre of $L^2\setminus 0$ over $0 \in \CP$ (which for simplicity we assume
is not a branch point of $S \to \CP$, we also assume $p \ne 0$). Thus consider the map
\begin{eqnarray}\label{symm_prod_param}
   M_k^\C(D) &\to& s^k(L^2\setminus 0)_0\\
   S &\mapsto& S\cap (L^2\setminus 0)_0\notag
\end{eqnarray}
Identifying $T^{1,0}_{[v]}$ with $H^0(S, \hat N(-1))$ as usual, the derivative of this
map at $S$ is the restriction map
\begin{eqnarray}\label{symm_prod_param_deriv}
   H^0(S, \hat N(-1)) \to H^0(S \cap (L^2\setminus 0)_0, \hat N(-1))
\end{eqnarray}
This fits into the long exact cohomology sequence associated to the following short
exact sequence of sheaves on $S$
\begin{eqnarray*}
   0 \to \hat N(-2) \overset{v_0}\to \hat N(-1) \to \hat N(-1)_{S\cap (L^2\setminus 0)_0}
   \to 0
\end{eqnarray*}
where $v_0$ is a non-zero section of $\cO(1)$ vanishing at $0$.
Since $H^0(S, \hat N(-2)) = H^1(S, \hat N(-2)) = 0$,
the map \eqref{symm_prod_param_deriv} is an isomorphism and so \ref{symm_prod_param} is a
biholomorphism near $S$. This provides us with some useful coordinates for $M_k^\C(D)$,
namely the $k$ points of intersection of a point $S$ with the fibre of $L^2\setminus 0$
over $0$. To be completely explicit, we introduce a patching description of
$L^2\setminus 0$. Thus we view $L^2\setminus 0$ as two copies of $\C\times\C\times\C^*$
glued together appropriately. More precisely let the first copy have coordinates
$\zeta, \eta, u$ and the second copy have coordinates $\tilde\zeta, \tilde\eta,
\tilde u$. Then the identification is defined by
\begin{eqnarray*}
   \tilde\zeta &=& \zeta^{-1}\\
   \tilde\eta &=& \zeta^{-2}\eta\\
   \tilde u &=& e^{\eta/\zeta}u
\end{eqnarray*}
Also let $s$ and $\tilde s$ be non-vanishing local sections of $\cO(1)$ on the relevant
patches defining local trivialisations (so that $\tilde s = \zeta s$). In these coordinates
our section $\rho$ appears as
\begin{eqnarray}\label{rho_coord_form}
   \rho = d\zeta\wedge d\eta\wedge (du/u) s^4
\end{eqnarray}
and similarly on the other coordinate patch. We also need to note that we have a natural
isomorphism
\begin{eqnarray*}
   \hat N \simeq T_F|_S
\end{eqnarray*}
where $T_F$ is the bundle of tangents to the fibres of the map $L^2\setminus 0 \to \CP$. In
terms of our coordinates $\zeta, \eta, u$ on $L^2\setminus 0$, $T_F$ is spanned by
$\frac{\partial}{\partial \eta}, \frac{\partial}{\partial u}$. In a neighbourhood of
the $k$ points of $S$ above $0 \in \CP$ we can represent the $k$ sheets of $S \to \CP$ as the
graphs of functions $\zeta \mapsto (\eta_i(\zeta), u_i(\zeta))$, $i=1, \ldots, k$ and so
an infinitesmal deformation $X$ of $S$ is represented on the $i^{\rm th}$ sheet by
\begin{eqnarray*}
   \eta'_i\frac{\partial}{\partial \eta} + u'_i\frac{\partial}{\partial u}
\end{eqnarray*}
for some functions $\eta'_i, u'_i$, $i=1,\ldots, k$. We can summarise our observations
above by saying that the map
\begin{eqnarray*}
   X \mapsto (\eta'_1(0), u'_1(0), \ldots, \eta'_k(0), u'_k(0))
\end{eqnarray*}
is an isomorphism for vectors $X$ such that $\eta'_i(\zeta_0) = u'_i(\zeta_0) = 0$
where $\zeta_0$ is the value of $\zeta$ corresponding to the point $p\in\CP$ defining $D$.
We will express $\omega_D$ in terms of these coordinates $\eta'_i(0), u'_i(0)$.
\par
Thus let $X_j$, $j=1,2$ be two tangent vectors to $M_k^\C(D)$ at $S$ represented on
the $i^{\rm th}$ sheet by
\begin{eqnarray*}
   \eta'_{i,j}\frac{\partial}{\partial \eta} + u'_{i,j}\frac{\partial}{\partial u}
\end{eqnarray*}
for $j=1,2$. To calculate $\omega_D(X_1, X_2)$, we first take a section of $\cO(1)$
vanishing at $p$. We use $v = (\zeta\zeta_0^{-1} - 1)s$. Dividing $X_j$ by $v$ we obtain
the expression
\begin{eqnarray*}
   \frac{1}{(\zeta\zeta_0^{-1} - 1)s}\left(\eta'_{i,j}\frac{\partial}{\partial\eta} 
   + u'_{i,j}\frac{\partial}{\partial u}\right)
\end{eqnarray*}
which we plug into our formula \eqref{rho_coord_form} to obtain
\begin{eqnarray*}
   \frac{\eta'_{i,1}u'_{i,2} - \eta'_{i,2}u'_{i,1}}{(\zeta\zeta_0^{-1} - 1)^2 u_i}d\zeta s^2
\end{eqnarray*}
This is an expression for the section of $H^0(S, K_S(2))$ on the $i^{\rm th}$ sheet.
\par
Finally, the canonical element of $H^1(\CP, \cO(-2))$ is represented as a \v{C}ech
cohomology class by $1/(\zeta s^2)$. Pulling this back to $S$, pairing it with the
above expression for our section of $K_S(2)$ and carrying out the contour integral
to implement Serre duality we obtain
\begin{eqnarray*}
   \omega_D(X_1, X_2) &=& \sum_{i=1}^k\int_{|\zeta| = 1}\frac{\eta'_{i,1} u'_{i,2}
   - \eta'_{i,2} u'_{i,1}} {(\zeta\zeta_0^{-1} - 1)^2 u_i}\frac{d\zeta}{\zeta}\\
   &=& \sum_{i=1}^k\frac{\eta'_{i,1}(0) u'_{i,2}(0) - \eta'_{i,2}
  (0) u'_{i,1}(0)}{u_i(0)}
\end{eqnarray*}
We have thus expressed $\omega_D$ in terms of the desired local coordinates on
$s^k(L^2\setminus 0) \simeq s^k(\C\times\C^*)$. We see that it is the usual
holomorphic product symplectic structure on $s^k(\C\times\C^*)$ (away from the
singularities). This establishes integrability and shows that we have recovered
the usual hyperk\"ahler structure on the monopole moduli space, see \cite{MR934202}.
\subsection{The Higgs field and the Atiyah class}\label{atiyah_higgs_subsect}
Consider the Penrose transform for $\End_0(E)$
\begin{eqnarray*}
   P : H^1(\T, \End_0(\widetilde E)(-2)) \to \left\{ \phi \in \Omega^0(
   \R^3, \End_0(E))\quad | \quad \Delta_{(A, \Phi)}\phi = 0 \right\}
\end{eqnarray*}
where
\begin{eqnarray*}
   \Delta_{(A, \Phi)}\phi = [\nabla^*, [\nabla, \phi]] + [\Phi^*, [\Phi, \phi]]
\end{eqnarray*}
Since (by the Bogomolny equations and the Bianchi identity) we have
\begin{eqnarray*}
\Delta_{(A, \Phi)}\Phi = 0
\end{eqnarray*}
there must be a natural class in $H^1(\T, \End_0(\widetilde E)(-2))$. The
question arises of identifying it in twistor terms. Before
we state the theorem which answers this question, let us first quickly
recall some elementary facts about holomorphic vector bundles.
\par
Thus let $V$ be a rank $n$ holomorphic vector bundle on a complex manifold
$X$. Let $\pi : F \to X$ be the principal $GL(n, \C)$ bundle on $X$
associated to $V$. We have the following natural exact sequence of vector
bundles on $F$
\begin{eqnarray*}
   0 \to \mathfrak{gl}(n, \C) \to TF \to \pi^*TX \to 0
\end{eqnarray*}
Everything in the sequence is acted on by $GL(n, \C)$ and the actions are
compatible with the maps so that we can take quotients to get the following
exact sequence of vector bundles on $X$
\begin{eqnarray*}
   0 \to \End(V) \to TF/GL(n, \C) \to TX \to 0
\end{eqnarray*}
The Atiyah class
\begin{eqnarray*}
\Lambda \in H^1(X, T^*X\otimes \End(V))
\end{eqnarray*}
of $V$ is defined to be the extension class of this sequence\footnote{Another
way to define the Atiyah class is to say it is the extension class of the
first jet bundle sequence of $V$
\begin{eqnarray*}
   0 \to T^*X\otimes V \to J_1(V) \to V \to 0
\end{eqnarray*}}.
All we need to know is that if
$A$ is a smooth connection on $V$ compatible with the holomorphic structure
on $V$ (ie: such that $F_A^{0,2} = 0$ where $F_A$ is the curvature of $A$)
then the $(1,1)$ part of the curvature
\begin{eqnarray*}
   F_A^{1,1} \in \Omega^{1,1}(X, \End(V)) = \Omega^{0,1}(X, T^*X\otimes\End(V))
\end{eqnarray*}
is a Dolbeault representative for $\Lambda$, see \cite{MR0098194}. (That 
$\deebar F_A^{1,1} = 0$ follows from the Bianchi identity and the fact that
$F_A^{0,2} = 0$.) We are ready to state
\begin{theorem}\label{atiyah_higgs_Euc_thm}
   Let $\widetilde E$ be the holomorphic vector bundle on $\T$
   corresponding to a solution $(A, \Phi)$ of the $SU(2)$ Bogomolny equations
   on $E \to \R^3$. Let
   \begin{eqnarray*}
      \Lambda \in H^1(\T, T^*\T\otimes \End(\widetilde E))
   \end{eqnarray*}
   be the Atiyah class of $\widetilde E$ and let $\Gamma : T^*\T \to \cO(-2)$
   be the dual of the natural inclusion $\cO(2) \hookrightarrow T\T$. Then
   \begin{eqnarray*}
      P(\Gamma(\Lambda)) = 4\pi\Phi
   \end{eqnarray*}
   where $P$ is the Penrose transform.
\end{theorem}
\noindent {\bf Proof} 
Suppose that (by abuse of notation) $\Lambda \in \Omega^{1,1}(\T, \End(
\widetilde E))$ is a Dolbeault representative for the Atiyah class. Fix
a point $O$ in $\R^3$ and let $j : Z \hookrightarrow \T$ denote the
inclusion of the corresponding twistor line $Z$ into twistor space.
Let $t : \widetilde E|_{Z} \to E_O$ be the natural trivialisation. Then
\begin{eqnarray*}
   tj^*\Gamma(\Lambda) \in \Omega^{1,1}(Z, E_O)
\end{eqnarray*}
is a $(1,1)$ form on $Z \simeq \CP$ taking values in the vector space $E_O$
and the value of the Penrose transform of $\Gamma(\Lambda)$ at $O$ is
\begin{eqnarray*}
   P(\Gamma(\Lambda))(O) = \int_Z t j^* \Gamma(\Lambda)
\end{eqnarray*}
\par
Now recall from chapter \ref{SingHypChap} that if $f : \T \to
\R^3$ is the function that takes a geodesic $\gamma \in \T$ to its closest
point to $O$ then we can define $\widetilde E = f^*E$ with $\deebar$-operator
\begin{eqnarray*}
   \deebar = (\widetilde\nabla - i\tilde\Phi\otimes\theta)^{0,1}
\end{eqnarray*}
where $\widetilde\nabla = f^*\nabla$, $\tilde\Phi = f^*\Phi$ and $\theta$ is
the real 1-form on $\T$ such that $\theta^{0,1} = \frac{2\eta}{(1 +
\zeta\bar\zeta)^2}d\bar\zeta$ in the usual coordinates $(\eta, \zeta)$ for
$\T \simeq T\CP$ in which $Z$ now appears as the zero section $\{\eta = 0\}$.
This means that a Dolbeault representative for the Atiyah class is the $(1,1)$
part of the curvature $F_{\hat\nabla}$ of the connection $\hat\nabla =
\widetilde\nabla - i\tilde\Phi\otimes\theta$.
\par
Also note that with this definition of $\widetilde E$, since $f(\gamma) = O$
for $\gamma \in Z$, the trivialisation $t$ is just the natural identification
\begin{eqnarray*}
   t : \widetilde E_\gamma = E_{f(\gamma)} = E_O
\end{eqnarray*}
Thus since the curvature of the connection $\hat\nabla$ is given by
\begin{eqnarray*}
   F_{\hat\nabla} = f^*F_\nabla + i\theta\wedge f^*(\nabla\Phi) - if^* \Phi
   \otimes d\theta
\end{eqnarray*}
we have
\begin{eqnarray*}
   P(\Gamma(\Lambda))(O) = \int_Z j^* \Gamma((f^* F_\nabla)^{1,1}
   + (i\theta \wedge f^*(\nabla\Phi))^{1,1} - if^*\Phi\otimes d\theta)
\end{eqnarray*}
We deal with the three terms of the integrand separately. Firstly
\begin{eqnarray*}
   \int_Z j^*\Gamma(f^*\Phi\otimes d\theta) &=& 
   \int_Z (f\circ j)^*\Phi\otimes j^*\Gamma(d\theta)\\
   &=& \Phi(O) \otimes \int_Z j^*\Gamma(d\theta)
\end{eqnarray*}
But in the coordinates $(\eta, \zeta)$ on $\T$ we have
\begin{eqnarray*}
   \theta = \frac{2}{(1 + |\zeta|^2)^2}\left(\eta d\overline\zeta +
   \overline \eta d\zeta\right)
\end{eqnarray*}
and
\begin{eqnarray*}
   \Gamma : \left\{
   \begin{aligned}
      d\zeta &\mapsto& 0\\
      d\eta &\mapsto& d\zeta
   \end{aligned}
   \right.
\end{eqnarray*}
Thus
\begin{eqnarray*}
   \int_Z j^*\Gamma(d\theta) = \int_{\C}\frac{2}{(1 + |\zeta|^2)^2}
   d\zeta\wedge d\overline\zeta = 4\pi i
\end{eqnarray*}
\par
Next note that
\begin{eqnarray*}
   j^*\Gamma \left(\theta\wedge f^*(\nabla\Phi)\right)^{1,1} &=& j^*
   \left((\Gamma\theta^{1,0})\wedge (f^*(\nabla\Phi))^{0,1} + \theta^
   {0,1}\wedge\Gamma(f^*(\nabla\Phi))^{1,0}\right) = 0
\end{eqnarray*}
since $\theta$ vanishes along $Z$.
\par
Evidently, to prove the theorem we need only show that the final term
$j^*\Gamma (f^*F_\nabla)^{1,1}$  in the integrand contributes nothing.
To show this it is sufficient to prove that the composition
\begin{eqnarray*}
   \wedge^2 T^*_O\R^3 \overset{f^*}\longrightarrow\wedge^2 T^*_\gamma\T
   \longrightarrow \wedge^{1,1}T^*_\gamma\T \overset{\Gamma}\longrightarrow
   \cO(-2)_\gamma\otimes \overline T^*_\gamma\T \overset{j^*}\longrightarrow
   \wedge^{1,1}T^*_\gamma Z
\end{eqnarray*}
is zero for every $\gamma \in Z$.
\par
To see this, suppose that $\omega \in \wedge^2 T^*_O\R^3$ and let
\begin{eqnarray*}
   f^*\omega = a_1d\zeta\wedge d\eta + a_2 d\zeta\wedge d\overline\eta
   + a_3 d\eta\wedge d\overline\zeta\\
   + a_4 d\zeta\wedge d\overline\zeta + a_5 d\eta\wedge
   d\overline\eta + a_6 d\overline\zeta\wedge d\overline\eta
\end{eqnarray*}
Using the formula for $\Gamma$ we find
\begin{eqnarray*}
   j^*\Gamma(f^*\omega)^{1,1} = a_3d\zeta\wedge d\overline\zeta
\end{eqnarray*}
We thus need to show $a_3 = 0$. To do this we shall use the formula
for $f$ as a function of our coordinates $(\zeta, \eta)$. This is
\begin{eqnarray*}
   f(\eta, \zeta) = \frac{1}{(1 + |\zeta|^2)^2}\Re\left\{\overline\eta
   (1 - \zeta^2, i(1 + \zeta^2), 2\zeta)\right\}
\end{eqnarray*}
Using the above formula for $f$ we find
$\frac{\partial f}{\partial\overline\zeta}(\zeta, 0) = 0$ from which it follows that
$a_3 = 0$ as required. \myqed
\section{The hyperbolic case}
\subsection{Geometry on the moduli space}
Hyperbolic monopoles have a description in terms of spectral curves that
has much in common with the spectral curve description of Euclidean
monopoles. A hyperbolic monopole of mass $m$ and charge $k$ has a
spectral curve $S \subset \QQ$ which is a compact algebraic curve in the
linear system $|\cO(k,k)|$ such that
\begin{itemize}
   \item
   $S$ is real
   \item
   $L^{2m+k}$ is trivial on $S$
\end{itemize}
As usual, $L$ is the line bundle $\cO(1,-1)$.
\par
Although all of the conditions for a compact algebraic curve in the linear
system $|\cO(k,k)|$ to be the spectral curve of a hyperbolic monopole are
not known, the two conditions above are the only ones that we shall need
(just as in the Euclidean case).
\par
We are interested in the moduli space $M_k$ of gauged hyperbolic monopoles
which naturally corresponds to the space of pairs consisting of a spectral
curve $S$ together with a trivialisation of $L^{2m+k}$ (satisfying an
appropriate reality condition) over $S$. We therefore
use the trivialisation to lift the curve $S$ to its image $\hat S$ in
$L^{2m+k}\setminus 0$ and denote the normal bundle
by $\hat N$. We shall again show $H^1(\hat S, \hat N) = 0$ and apply Kodaira's
theorem to obtain the space of deformations of $\hat S$ in $L^{2m+k}\setminus 0$.
Again the space of deformations carries a real structure for which the real
points correspond to curves $S$ that are real. We thus again find
\begin{eqnarray*}
T_{\hat S}M_k \otimes{_\R}\C \simeq H^0(\hat S, \hat N)
\end{eqnarray*}
\par
As in the Euclidean case, we must first identify the normal bundle of $\hat S$
in $L^{2m+k}$. Appealing to lemma \ref{NormIdLemma} we have the following
\begin{proposition}\label{hyp_nrm_bnd_id}
   Let $\hat N$ be the normal bundle of $\hat S \subset L^{2m+k}$. Then
   (identifying $\hat S$ and $S$)
   \begin{eqnarray*}
      \hat N \simeq \widetilde EL^m(k,0)|_{S}
   \end{eqnarray*}
   (Note that $L^m(k,0)|_{S} \simeq L^{-m}(0,k)|_{S}$.)
\end{proposition}
\noindent {\bf Proof} This follows from lemma \ref{NormIdLemma} together with
the construction of $\widetilde E$ from $S$ (see \cite{MR1399482}). \myqed
\par
The next result we need is the cohomology vanishing
theorem which (amongst other things) ensures $H^1(\hat S, \hat N) = 0$.
\begin{theorem}\label{hyp_vanishing}
   Let $\widetilde E$ be the holomorphic vector bundle on $\Q$
   corresponding to an $SU(2)$ monopole of charge $k$ on $\Hyp$. Let
   $S$ be the spectral curve of the monopole. Then
   \begin{eqnarray*}
      H^0(S, \widetilde EL^m(k-1, -1)) = 0
   \end{eqnarray*}
\end{theorem}
\noindent {\bf Proof} The proof is similar to the corresponding
result, theorem \ref{Euc_vanishing_theorem}, in the Euclidean case. Reusing
the notation where possible from the proof of \ref{Euc_vanishing_theorem},
we sketch the details as follows.
\par
Using the exact sequences (see \cite{MR1399482})
\begin{eqnarray*}
   0 \to L^m(0, -k) \to \widetilde E \to L^{-m}(0, k) \to 0
\end{eqnarray*}
and
\begin{eqnarray*}
   0 \to L^{-m}(-k, 0) \to \widetilde E \to L^m(k, 0) \to 0
\end{eqnarray*}
together with the result (see \cite{MR1961371}) that $H^0(\Q, L^s(n, 0)) = 0$
for any $s \in (0,\infty)$ and any non-negative $n \in \Z$ we obtain injections
\begin{eqnarray*}
   H^0(S, \widetilde EL^m(k-1, -1)) \overset{\delta}\hookrightarrow H^1(\Q,
   \widetilde EL^m(-1, -k-1)) \overset{i}\hookrightarrow H^1(\Q, S^2\widetilde E(-1, -1))
\end{eqnarray*}
We also have an injection
\begin{align}\label{hyp_t_inject}
   H^0(S, \widetilde EL^m(k-1, -1)) \hookrightarrow H^0(S, \cO(k-1, k-1))
   \simeq H^0(\Q, \cO(k-1, k-1))
\end{align}
Finally we have the Penrose transform (see subsection \ref{Hyp_Higgs_Atiyah} for a
brief discussion of the Penrose transform in hyperbolic case)
\begin{eqnarray*}
   P : H^1(\Q, S^2\widetilde E(-1, -1)) \to \left\{ \phi \in \Omega^0(\Hyp, S^2 E)
   \quad |\quad \Delta_{(A, \Phi)}\phi = \phi\right\}
\end{eqnarray*}
where as usual
\begin{eqnarray*}
   \Delta_{(A, \Phi)} = \nabla^*\nabla + \Phi^*\Phi
\end{eqnarray*}
Thus let $s \in H^0(S, \widetilde EL^m(k-1, -1))$ and let
\begin{eqnarray*}
   t \in H^0(\Q, \cO(k-1, k-1))
\end{eqnarray*}
be the corresponding element under the injection \eqref{hyp_t_inject}.
Also let $\phi = P(i\delta s)$. Fix a point $O \in \Hyp$. Then, just as in
the Euclidean case, if $x \in \Hyp$ and $v_1, v_2$ are an $SU(2)$ basis for
$E_x$ and $f_1, f_2$ are the corresponding sections of $\widetilde E$ on the
twistor line $P_x \subset \Q$ corresponding to $x$, we have the following formula
for $\phi$ (cf. equation \eqref{Euc_phi_formula})
\begin{eqnarray}\label{hyp_phi_formula}
   \phi_{ij}(x) = \int_{V^+}\langle f_i, \theta^+ f_j\rangle + \int_{V^-}\langle
   f_i, \theta^- f_j\rangle + \int_{\partial H^-_x}\langle f_i, \gamma_b f_j\rangle
\end{eqnarray}
where $V^{\pm}$ are appropriate disjoint neighbourhoods of the points in $P_x$
corresponding to the oriented lines joining $O$ and $x$, $H^-_x$ is the hemisphere
in $P_x$ containing $V^-_x$ and whose boundary consists of those points
corresponding to the geodesics perpendicular to the line joining $x$ and $O$.
\par
Now on $\Q\setminus S$, we have the holomorphic splitting of $\widetilde E$
\begin{eqnarray}\label{E_splitting}
   \widetilde E = L^{-m}(-k, 0) \oplus L^m(k, 0)
\end{eqnarray}
and (cf equation \eqref{gamma_formula}) in terms of this splitting we have the
following formula for $\gamma_b \in \Omega^0(\Q, \End_0(\widetilde E)(-2))$ on
$\partial H^-_x$
\begin{eqnarray*}
   \gamma_b = \left[
   \begin{array}{cc}
      -2t/\psi & 0\\
      0 & 2t/\psi
   \end{array}\right]
\end{eqnarray*}
where $\psi \in H^0(\Q,\cO(k,k))$ is a section defining the spectral curve.
\par
Let $\rho$ be the hyperbolic distance from $x$ to $O$. Just as in
the Euclidean case (see also \cite{MR1961371}) the first two terms in
\eqref{hyp_phi_formula} have exponential decay in $\rho$, in fact they decay
like $e^{-(2m+1)\rho}$. We claim that the third term decays like $e^{-\rho}$.
Using the splitting \eqref{E_splitting} this term is
\begin{eqnarray*}
   \int_{\partial H^-_x}\langle f_i, \gamma_b f_j\rangle = 
   \int_{\partial H^-_x}2t/\psi
   \langle f_i, \left[
   \begin{array}{cc}
      -1 & 0\\
      0 & 1
   \end{array}\right] f_j\rangle
\end{eqnarray*}
Using the boundary conditions given in part (ii) of definition \ref{def_sngmpole},
exactly the same method of proof as that used in the proof of lemma
\ref{vanishing_lemma} establishes that $\langle f_i, \left[\begin{array}{cc} -1
& 0\\0 & 1\end{array} \right] f_j\rangle$ is bounded as $r \to \infty$ and so we
need only show that $t/\psi$ has exponential decay on $\partial H^-$ as $r\to\infty$.
To see this, note that we may choose our coordinates $(z,w)$ on $\Q \simeq \QQ$ such
that the oriented line $Ox$ has coordinates $(0, 0)$. A quick calculation (best
done in the upper half space model of $\Hyp$) reveals that in these coordinates
\begin{eqnarray*}
   P_x = \left\{ (z, w) \in \Q\quad |\quad z = e^{-2\rho}w\right\}
\end{eqnarray*}
Recalling that the circle $\partial H^-_x \subset P_x$ consists of those lines through
$x$ perpendicular to the line $Ox$, we find
\begin{eqnarray*}
   \partial H^-_x = \left\{(z, \overline z^{-1}) \in \Q \quad |\quad |z| = e^{-\rho}\right\}
\end{eqnarray*}
Let
\begin{eqnarray*}
   t(z, w) = p_1(z) + p_2(z)w + \cdots + p_{k-1}(z)w^{k-1}
\end{eqnarray*}
where $p_i$ is a polynomial of degree $k-1$ for $i=1,\ldots,k-1$. Thus, on
$\partial H^-_x$ we have
\begin{eqnarray*}
   |t(z,w)| &\le& |p_1(z)| + |p_2(z)||w| + \cdots + |p_{k-1}(z)||w|^{k-1}\\
            &=& |p_1(z)| + |p_2(z)|e^\rho + \cdots + |p_{k-1}(z)|e^{(k-1)\rho}\\
            &\le& C_t e^{(k-1)\rho}
\end{eqnarray*}
for some constant $C_t$ and large enough $\rho$. Next we claim that there exists
a constant $C_\psi > 0$ such that on $\partial H^-_x$
\begin{eqnarray*}
   |\psi(z,w)| \ge C_\psi e^{k\rho}
\end{eqnarray*}
for large enough $r$. Thus, note that we have
\begin{eqnarray*}
   \psi(z,w) = q_1(z) + q_2(z)w + \cdots + q_k(z)w^k
\end{eqnarray*}
where $q_i$ is a polynomial of degree $k$ for $i=1,\ldots,k$. A simple polynomial
estimate as above will give the required inequality for $\psi$ provided the coefficient
of the $w^k$ term in the expression for $\psi$ is non-zero. This is the same as saying
$\psi(\infty, 0) \ne 0$. However this follows since we know that the spectral curve $S$
does not meet the anti-diagonal $\overline \Delta$. We thus have that there exist
constants $C$ and $R$ such that for $\rho > R$ we have
\begin{eqnarray*}
   \left|\frac{t}{\psi}\right| \le Ce^{-\rho}
\end{eqnarray*}
on $\partial H^-_x$. Putting all of our estimates together, we find
\begin{eqnarray*}
   \|\phi\| = O(e^{-\rho})
\end{eqnarray*}
\par
Finally we show that this decay of $\phi$ is in fact enough to ensure that it
vanishes. Again, as in the Euclidean case, we use a maximum principle for an appropriate
elliptic operator. However we cannot simply use the hyperbolic Laplacian since in place
of \eqref{Euc_maxprin_ineq} we only obtain
\begin{eqnarray*}
   \frac{1}{2}\Delta\|\phi\|^2 = \|\nabla\phi\|^2 + \|\Phi(\phi)\|^2 - \|\phi\|^2
\end{eqnarray*}
which is not obviously non-negative. Thus introduce upper half space coordinates $(x,y,z)$
on $\Hyp$ and use the elliptic operator
\begin{eqnarray}\label{Euc_Laplacian}
   \Delta_E = \frac{\partial^2}{\partial x^2} + \frac{\partial^2}{\partial y^2}
   + \frac{\partial^2}{\partial z^2} + z^{-1}\frac{\partial}{\partial z}
\end{eqnarray}
This is the Laplacian for $S^1$-invariant functions on $\Hyp\times S^1$ but with the
metric conformally rescaled by a factor of $z^2$ to the flat Euclidean metric. Now a
quick calculation shows that
\begin{eqnarray*}
   z^3\nabla^{*_E}\nabla(z^{-1}\phi') = (\nabla^*\nabla - 1)\phi'
\end{eqnarray*}
where $*_E$ indicates that we're taking the adjoint of $\nabla$ with respect to the
Euclidean metric on $\Hyp\times S^1$ and $\phi'$ is any section. We thus have
\begin{eqnarray*}
   (\nabla^{*_E}\nabla + z^{-2}\Phi^*\Phi) (z^{-1}\phi) = 0
\end{eqnarray*}
and so, still working in the Euclidean metric, we have
\begin{eqnarray*}
   \frac{1}{2}\Delta_E\|z^{-1}\phi\|^2 = \|\nabla(z^{-1}\phi)\|^2_E + \|z^{-1}
   \Phi(z^{-1}\phi)\|^2 \ge 0
\end{eqnarray*}
   Thus, by the maximum principle for the differential operator $\eqref{Euc_Laplacian}$,
   $z^{-2}\|\phi\|^2$ satisfies the maximum principle on $\Hyp$. To see that this is
   enough to ensure vanishing we note
\begin{eqnarray*}
   \|\phi\| &=& O\left(e^{-\rho}\right)\\
   \Rightarrow \|\phi\| &=& O\left(\frac{1}{\cosh\rho}\right)\\
   \Rightarrow \|\phi\| &<& \frac{C}{\cosh\rho} \mbox{\quad outside a compact
   set in $\Hyp$}
\end{eqnarray*}
for an appropriate $C > 0$. Recalling that provided we choose our coordinate
functions $(x,y,z)$ such that $O$ has coordinates $(0, 0, 1)$
\begin{eqnarray*}
   \cosh\rho = \frac{x^2 + y^2 + z^2 + 1}{2z}
\end{eqnarray*}
we find that we have
\begin{eqnarray*}
   z^{-2}\|\phi\|^2 < \frac{4C^2}{(x^2 + y^2 + z^2 + 1)^2}
\end{eqnarray*}
outside a compact set. Since this is decaying and satisfies the maximum principle,
we must have $\phi = 0$ as required. \myqed
\begin{corollary} In the notation of corollary \ref{hyp_nrm_bnd_id} we have
   \begin{enumerate}
      \item
      $h^1(\hat S, \hat N) = 0$
      \item
      $h^0(\hat S, \hat N) = 4k$
      \item
      $h^0(\hat S, \hat N(-1, 0)) = 2k$
      \item
      $h^0(\hat S, \hat N(0, -1)) = 2k$
   \end{enumerate}
\end{corollary}
\noindent {\bf Proof} Firstly, note that in view of corollary \ref{hyp_nrm_bnd_id}
 we can restate theorem \ref{hyp_vanishing} as
\begin{eqnarray*}
   H^0(S, \hat N(-1, -1)) = 0
\end{eqnarray*}
Next note that in view of lemma \ref{NormExtLemma} we have the short exact sequence
\begin{eqnarray*}
   0 \to \cO \to \hat N \to \cO(k, k) \to 0
\end{eqnarray*}
We thus have $\wedge^2\hat N \simeq \cO(k, k)$ and so
\begin{eqnarray*}
   \hat N^*(k,k) \simeq \hat N
\end{eqnarray*}
Furthermore, this sequence also allows us to calculate the Chern class of $\hat N$.
Indeed, recalling that $S$ is a (branched) k-fold covering of the left and right
$\CP$ in $\Q\simeq \QQ$, we have
\begin{eqnarray*}
   c_1(\hat N)[S] = 2k^2
\end{eqnarray*}
Finally recall that since the spectral curve is a divisor of $\cO(k,k)$ on $\Q$
its canonical bundle is (by adjunction) $\cO(k-2, k-2)$.
\begin{enumerate}
   \item
   By Serre duality and our observations above
   \begin{eqnarray*}
      h^1(S, \hat N) = h^0(S, \hat N^*(k-2, k-2)) = h^0(S, \hat N(-2, -2)) = 0
   \end{eqnarray*}
   since $h^0(S, \hat N(-1, -1)) = 0$ (and $h^0(S, \cO(1, 1)) > 0$).
   \item
   By the Riemann-Roch formula
   \begin{eqnarray*}
      h^0(S, \hat N) - h^1(S, \hat N) &=& (2 + c_1(\hat N))(1 + \frac{1}{2})
      c_1(TS))[S]\\
      &=& 4k
   \end{eqnarray*}
   and $h^1(S, \hat N) = 0$ by part (i).
   \item
   By the Riemann-Roch formula
   \begin{eqnarray*}
      h^0(S, \hat N(-1, 0)) - h^1(S, \hat N(-1, 0)) = 2k
   \end{eqnarray*}
   We thus need only show $h^1(S, \hat N(-1, 0)) = 0$ but this follows immediately
   since
   \begin{eqnarray*}
      h^1(S, \hat N(-1, 0)) = h^0(S, \hat N^*(k-1, k-2)) = h^0(S, \hat N(-1, -2)) = 0
   \end{eqnarray*}
   \item
   This follows just as in (iii).
\end{enumerate}
\myqed
\begin{remark}
   Note that in \cite{MR1961371} it is shown that
   \begin{eqnarray*}
      H^0(S, L^z(k-1, -1)) = 0\mbox{\quad for $z \in [0, 2m]$}
   \end{eqnarray*}
   Our proof that $H^0(S, \widetilde EL^m(k-1,-1))=0$
   does not generalise to a possible vanishing for $H^0(S, \widetilde EL^z(k-1, -1))$
   and instead works only for the midpoint value of $z$ (which is of course,
   as in the Euclidean case, precisely the value we need).
\end{remark}
\par
We thus see that, as in the Euclidean case, the obstruction class of Kodaira's
deformation theory for $\hat S \subset L^{2m+k}\setminus 0$ vanishes and that
we obtain a family of deformations of the correct dimension $4k$.
\par
It is here that the differences between the hyperbolic and Euclidean cases
begin to emerge. Firstly, in view of the above results we see that we have
{\it two} natural tensor product decompositions of the complexified tangent space
to the hyperbolic monopole moduli space:
\begin{eqnarray}\label{hyp_mod_space_decomp}
   \begin{aligned}
      T_{\hat S}M_k \otimes_{\R}\C &\simeq& H^0(S, \hat N(-1, 0))\otimes_{\C}
      H^0(S, \cO(1, 0))\\
      T_{\hat S}M_k \otimes_{\R}\C &\simeq& H^0(S, \hat N(0, -1))\otimes_{\C}
      H^0(S, \cO(0, 1))
   \end{aligned}
\end{eqnarray}
These decompositions are of course the analogue of the single decomposition
\eqref{EucModTanDecom}. Crucially, there is also an analogue of lemma
\ref{EucModSkewForm}.
\begin{lemma}
   $H^0(S, \hat N(-1, 0))$ and $H^0(S, \hat N(0, -1))$ are naturally anti-isomorphic
   and carry natural non-degenerate complex skew-symmetric pairings.
\end{lemma}
\noindent {\bf Proof} The real structure induces the anti-isomorphism between
   $H^0(S, \hat N(-1, 0))$ and $H^0(S, \hat N(0, -1))$. We shall
   demonstrate the existence of the skew pairing for $H^0(S, \hat N(-1, 0))$
   since the other case is completely analogous. We begin by noting that
   we have a natural isomorphism
   \begin{eqnarray*}
      \wedge^2\hat N \simeq K_S\otimes K^*_{\Q}|_S
   \end{eqnarray*}
   where $K_S$ and $K_{\Q}$ are the canonical bundles of $S$ and $\Q$ respectively.
   Note also that we have a natural isomorphism $K_{\Q} \simeq \cO(-2, -2)$. Now
   let $\lambda \in H^1(S, \cO(0, -2)) \simeq H^1(\CP, \cO(-2))$ be the tautological
   element. Using the above, we define the skew form $\omega$ via the following
   bilinear map
   \begin{eqnarray*}
      H^0(S, \hat N(-1, 0)) \times H^0(S, \hat N(-1, 0)) &\to& H^0(S, K_S\otimes
      K_{\Q}^*(-2, 0))\\
      &\simeq& H^0(S, K_S(0, 2))\\
      &\overset{\lambda}\to& H^1(S, K_S)\\
      &\simeq& \C
   \end{eqnarray*}
   To see that $\omega$ is non-degenerate, note that $\lambda$ is the extension
   class of
   \begin{eqnarray*}\label{hyp_lambda_ext}
      0 \to \cO(0, -1) \to \C^2 \to \cO(0, 1) \to 0
   \end{eqnarray*}
   Now $\omega$ defines a map
   \begin{eqnarray*}
      \hat\omega : H^0(S, \hat N(-1, 0)) &\to& H^0(S, \hat N(-1, 0))^*\\
      &\simeq& H^1(S, K_S\otimes\hat N^*(1, 0))\\
      &\simeq& H^1(S, \hat N\otimes K_{\Q}(1, 0))\\
      &\simeq& H^1(S, \hat N(-1, -2))
   \end{eqnarray*}
   However this map can also be recognised as the
   connecting homomorphism of the long exact sequence of cohomology groups
   associated to the short exact sequence of vector bundles
   \begin{eqnarray*}
      0 \to \hat N(-1, -2) \to \C^2\otimes \hat N(-1, -1) \to \hat N(-1, 0) \to 0
   \end{eqnarray*}
   obtained by tensoring the sequence \eqref{hyp_lambda_ext} with $\hat N(-1, -1)$.
   Thus
   \begin{eqnarray*}
      \Ker\hat\omega = H^0(S, \C^2\otimes \hat N(-1, -1))
   \end{eqnarray*}
   which vanishes by theorem \ref{hyp_vanishing} and corollary
   \ref{hyp_nrm_bnd_id}. \myqed
\par
In the Euclidean case we saw that the decomposition \eqref{EucModTanDecom}
together with the skew form on $H^0(S, \hat N(-1))$ reflected the underlying
hyperk\"ahler structure of the moduli space. In the hyperbolic case it is as yet
unclear what type of geometry our analogous decompositions
\eqref{hyp_mod_space_decomp} together with the skew forms we have identified reflect.
\par
The geometry we have identified on the hyperbolic monopole moduli space appears
to be a real geometry whose complexification is very like the
complexification of hyperk\"ahler geometry but which is subtly different.
Furthermore whatever the geometry, it should converge to hyperk\"ahler geometry
in the limit as the mass of the monopoles tends to infinity.
\subsection{The Higgs field and the Atiyah class}\label{Hyp_Higgs_Atiyah}
Not surprisingly there is an analogue of theorem \ref{atiyah_higgs_Euc_thm}
in the hyperbolic case. Corresponding to the Euclidean Penrose transform
\begin{eqnarray*}
   H^1(\T, \cO(-2)) \simeq \left\{ \phi \in C^\infty(\R^3)\quad | \quad
   \Delta\phi = 0\right\}
\end{eqnarray*}
we have the hyperbolic version
\begin{eqnarray*}
   H^1(\Q, \cO(-1, -1)) \simeq \left\{ \phi \in C^\infty(\Hyp)\quad | \quad
   \Delta\phi = \phi\right\}
\end{eqnarray*}
The hyperbolic Penrose transform gives 1-eigenfunctions of the Laplacian
rather than harmonic functions because the Penrose transform for
the Laplacian on $\Hyp$ is really the $S^1$-invariant version of the Penrose
transform for the conformal Laplacian on $\Hyp\times S^1$. The conformal Laplacian
on an oriented 4-manifold with metric $g$ is $\Delta_g + R/6$ where $\Delta_g$
is the metric Laplacian and $R$ is the scalar curvature. $\Hyp\times S^1$ has
scalar curvature $R = -6$.
\par
The above isomorphisms generalise to
\begin{eqnarray*}
   H^1(\T, L^s(-2)) \simeq \left\{ \phi \in C^\infty(\R^3)\quad | \quad
   \Delta\phi = -s^2\phi\right\}
\end{eqnarray*}
where $L$ is the line bundle with transition function $e^{-\eta/\zeta}$ on
$\T$ and \cite{MR1961371}
\begin{eqnarray*}
   H^1(\Q, L^s(-1, -1)) \simeq \left\{ \phi \in C^\infty(\Hyp)\quad | \quad
   \Delta\phi = (1 - s^2)\phi\right\}  
\end{eqnarray*}
where $L = \cO(1, -1)$.
\par
In particular, taking $s=\pm 1$ in the hyperbolic case, we have
\begin{eqnarray*}
   H^1(\Q, \cO(-2, 0)) \simeq H^1(\Q, \cO(0, -2)) \simeq \left\{ \phi \in
   C^\infty(\Hyp)\quad | \quad \Delta\phi = 0\right\}
\end{eqnarray*}
Finally note that as usual we may couple the Penrose transform to a
solution of the self-duality equations (in this case the Bogomolny
equations) to obtain a covariant version. Thus if $\widetilde E$ is the
holomorphic bundle on $\Q$ corresponding to a solution $(A, \Phi)$
of the $SU(2)$ Bogomolny equations on $E \to \Hyp$, the Penrose transform
in the adjoint representation provides natural isomorphisms
\begin{eqnarray*}
   P_L : H^1(\Q, \End(\widetilde E)(-2, 0)) &\simeq& \left\{ \phi \in \Omega^0
   (\Hyp, \End(E))\quad | \quad \Delta_{(A, \Phi)}\phi = 0\right\}
\end{eqnarray*}
and
\begin{eqnarray*}
   P_R : H^1(\Q, \End(\widetilde E)(0, -2)) &\simeq& \left\{ \phi \in \Omega^0
   (\Hyp, \End(E))\quad | \quad \Delta_{(A, \Phi)}\phi = 0\right\}
\end{eqnarray*}
where $\Delta_{(A, \Phi)}\phi = \left[\nabla_A^*,\left[\nabla_A,\phi\right]
\right] + \left[\Phi^*,\left[\Phi, \phi\right]\right]$. This is important
because as in the Euclidean case we have $\Delta_{(A, \Phi)}\Phi = 0$ and
so the question arises of identifying the natural classes in $H^1(\Q, \End
\widetilde (E)(0, -2))$ and $H^1(\Q, \End\widetilde (E)(-2, 0))$ which correspond
to $\Phi$ under the Penrose transform.
\par
We answer this question with the following
\begin{theorem}
   Let $\widetilde E$ be the holomorphic vector bundle on $\Q$
   corresponding to a solution $(A, \Phi)$ of the $SU(2)$ Bogomolny equations
   on $E \to \Hyp$. Let
   \begin{eqnarray*}
      \Lambda \in H^1(\Q, T^*\Q\otimes\End(\widetilde E))
   \end{eqnarray*}
   be the Atiyah class of $\widetilde E$ and let $\Gamma_L : T^*\Q = \cO(-2, 0)
   \oplus \cO(0, -2) \to \cO(-2, 0)$ be the natural projection. Then
   \begin{eqnarray*}
      P_L(\Gamma_L(\Lambda)) = 4\pi\Phi
   \end{eqnarray*}
   (Similarly the corresponding result for $\Gamma_R$ and $P_R$ holds.)
\end{theorem}
\noindent {\bf Proof}
   As in the proof of theorem \ref{atiyah_higgs_Euc_thm}, if $O \in \Hyp$
   we have
   \begin{eqnarray*}
      P_L(\Gamma_L(\Lambda))(O) = \int_\Delta j^*\Gamma_L((f^*F_A)^{1,1}
      + i(\theta\wedge f^*(\nabla\Phi))^{1,1} - if^*\Phi\otimes d\theta)
   \end{eqnarray*}
   where $\Delta \subset \Q \simeq \QQ$ is the diagonal (the twistor
   line of the point $O \in \Hyp$), $j : \Delta \to \Q$ is the inclusion,
   $f : \Q \to \Hyp$ is the map taking a geodesic $\gamma \in \Q$ to its
   closest point to $O$ and $\theta$ is the real 1-form on $\Q$ such that
   \begin{eqnarray*}
      \theta_{\Q}^{0,1} = (z-w)\left(\frac{d\bar z}{(1+z\bar z)(1+\bar zw)} +
      \frac{d\bar w}{(1+w\bar w)(1+z\bar w)}\right)
   \end{eqnarray*}
   in the coordinates $([z,1],[w,1])$ on (an open set of) $\Q$.
   \par
   As before, we deal with the three terms of the integrand separately.
   Firstly
   \begin{eqnarray*}
      \int_\Delta j^*\Gamma_L(f^*\Phi\otimes d\theta) &=& \int_\Delta
      (f\circ j)^*\Phi \otimes j^*\Gamma_L(d\theta)\\
      &=&\Phi(O)\otimes \int_\Delta j^*\Gamma_L(d\theta)
   \end{eqnarray*}
   But in the coordinates $(z, w)$ on $\Q$, we have
   \begin{eqnarray*}
      \Gamma_L : \left\{
      \begin{aligned}
         dz &\mapsto& dz\\
         dw &\mapsto& 0
      \end{aligned}
      \right.
   \end{eqnarray*}
   Using our above formula for $\theta$ we find
   \begin{eqnarray*}
      \int_\Delta j^*\Gamma_L(d\theta) = \int_\C \frac{2}{(1 + |z|^2)^2}
      dz\wedge d\bar z = 4\pi i
   \end{eqnarray*}
   Next note that
   \begin{eqnarray*}
      j^*\Gamma_L \left(\theta\wedge f^*(\nabla\Phi)\right)^{1,1} &=& j^*
      \left((\Gamma_L\theta^{1,0})\wedge (f^*(\nabla\Phi))^{0,1} + \theta^
      {0,1}\wedge\Gamma_L(f^*(\nabla\Phi))^{1,0}\right) = 0
   \end{eqnarray*}
   since $\theta$ vanishes along $\Delta$.
   \par
   Evidently, to finish the proof of the theorem we need only show that the
   final term $j^*\Gamma_L (f^*F_A)^{1,1}$ in the integrand contributes
   nothing. To show this it is sufficient to prove that the composition
   \begin{eqnarray*}
      \wedge^2 T^*_O\Hyp \overset{f^*}\longrightarrow\wedge^2 T^*_\gamma\Q
      \longrightarrow \wedge^{1,1}T^*_\gamma\Q \overset{\Gamma_L}\longrightarrow
      \cO(-2,0)_\gamma\otimes \overline T^*_\gamma\Q \overset{j^*}\longrightarrow
      \wedge^{1,1}T^*_\gamma\Delta
   \end{eqnarray*}
   is zero for every $\gamma \in \Delta$.
   \par
   To see this, suppose that $\omega \in \wedge^2 T^*_O\Hyp$ and let
   \begin{eqnarray*}
      f^*\omega &=& a_1 dz\wedge dw + a_2 dz\wedge d\overline w + a_3
      dw\wedge d\overline z\\
      && a_4 dz\wedge d\overline z + a_5 dw\wedge d\overline w + a_6
      d\overline z\wedge d\overline w
   \end{eqnarray*}
   Using the formula for $\Gamma_L$ we find
   \begin{eqnarray*}
      j^*\Gamma_L (f^*F_A)^{1,1} = (a_2 + a_4) dz\wedge
      d\overline z
   \end{eqnarray*}
   We thus need to show $a_2 + a_4 = 0$. To do this we
   shall need the formula for $f$ as a function of our coordinates 
   $(z,w)$. Using the upper half-space model of $\Hyp$ with coordinates
   $(x + iy, u)$ and metric $\frac{dx^2 + dy^2 + du^2}{u^2}$, the formula for
   $f$ is
   \begin{eqnarray}\label{f_formula_hyp}
      f(z,w) = \mu(z,w)\left((1 + |w|^2)z - (1 + |z|^2)w,
      \sqrt{(1 + |z|^2)(1 + |w|^2)}|1 + z\overline w|\right)
   \end{eqnarray}
   where
   \begin{eqnarray*}
      \mu(z,w) = \frac{1}{1 + 2|w|^2 + |zw|^2}
   \end{eqnarray*}
   We need to know about the derivatives of $f$ on $\Delta$.
   Thus, let $f = (f_x + if_y, f_u)$. Using the above formula for $f$ we find
   \begin{eqnarray*}
      \left[
      \begin{array}{ccc}
         \frac{\partial f_x}{\partial z} & \frac{\partial f_x}{\partial \overline z}
         & \frac{\partial f_x}{\partial \overline w}\\
         &&\\
         \frac{\partial f_y}{\partial z} & \frac{\partial f_y}{\partial \overline z}
         & \frac{\partial f_y}{\partial \overline w}\\
         &&\\
         \frac{\partial f_u}{\partial z} & \frac{\partial f_u}{\partial \overline z}
         & \frac{\partial f_u}{\partial \overline w}
      \end{array}
      \right](z,z) = \frac{1}{2(1 + |z|^2)^2}\left[
      \begin{array}{ccc}
         1 - \overline z^2 & 1 - z^2 & -1 + z^2\\
         &&\\
         -i(1 + \overline z^2) & i(1 + z^2) & -i(1 + z^2)\\
         &&\\
         2\overline z & 2z & -2z
      \end{array}
      \right]
   \end{eqnarray*}
   (We have omitted the formulae for the derivatives of $f$ with respect to $w$ since
   we do not need them.)
   Using these formulae it is straightforward to verify that we do indeed have
   $a_2 + a_4 = 0$ on $\Delta$ as required. \myqed

\chapter{Further properties and open issues}\label{OpenIssueChap}
\section{Instantons and the hypercomplex quotient}\label{HCInstQuotSect}
One of the ways to see that the Euclidean monopole moduli spaces
carry a natural hyperk\"ahler structure is to use the hyperk\"ahler
quotient construction \cite{MR877637}. Indeed, as shown in
\cite{MR934202}, the Euclidean Bogomolny equations may be regarded 
(at least formally) as an infinite dimensional hyperk\"ahler moment map
for the action of the group of $SU(2)$ gauge transformations on $\R^3$
on the infinite dimensional hyperk\"ahler manifold that is the
space of pairs consisting of an $SU(2)$ connection and Higgs field on
$\R^3$.
\par
It is interesting to wonder whether the \emph{hyperbolic} Bogomolny equations
may be regarded as the moment map for the analogous action of the group
of gauge transformations in hyperbolic space provided we endow the space
of connections and Higgs fields with an appropriate geometry. This remains
an open question.
\par
Although there are many examples\footnote{Other than the Euclidean monopole
moduli spaces we also have Kronheimer and Nakajima's beautiful generalisation
\cite{MR1075769} of the ADHM construction \cite{MR598562} in mind as well
as the moduli spaces of instantons on a compact hyperk\"ahler 4-manifold
\cite{MR935967}.} of situations in which the self-duality equations
can be regarded as a moment map, they all occur in the presence
of hyperk\"ahler geometry. It would be interesting to
have an example in which the self-duality equations can be
regarded as a moment map in the presence of a different type of geometry.
In particular it would be interesting to have such an example in the
presence of a non-metric geometry since the geometry of the hyperbolic monopole
moduli space appears to be non-metric.
To this end we show that the moduli space of instantons on a compact
\emph{hypercomplex} 4-manifold may be obtained as an infinite dimensional
hypercomplex quotient in the sense of Joyce \cite{MR1109637}.
\par
In addition, as we shall see later, there is
another reason why this result about instantons on a hypercomplex 4-manifold
is relevant to the geometry of hyperbolic monopoles. This is because the
quotient of certain Hopf surfaces by an appropriate $S^1$-action is a
hyperbolic solid torus.
\par
Let us briefly recall the details of both the hyperk\"ahler and hypercomplex
quotient constructions. Thus suppose that a Lie group $G$ acts on a
hyperk\"ahler manifold $(M, \omega_1, \omega_2, \omega_3)$ by diffeomorphisms
preserving the hyperk\"ahler structure. Recall \cite{MR877637} that a
hyperk\"ahler moment map for this action is an equivariant map
\begin{eqnarray*}
   \mu = (\mu_1, \mu_2, \mu_3) : M \to \g^*\otimes\R^3
\end{eqnarray*}
such that if $\lambda \in \g$ and $X_\lambda$ is the vector field on $M$
generated by $\lambda$ via the $G$-action then
\begin{eqnarray*}
   \langle\lambda, d\mu_i\rangle = \omega_i(X_\lambda, \cdot) \qquad\mbox{for $i=1,2,3$}
\end{eqnarray*}
When a moment map exists, the quotient manifold $N = \mu^{-1}(\{\xi\})/G$
admits a natural hyperk\"ahler structure and we say $N$ has been obtained
from $M$ by a hyperk\"ahler quotient.
\par
Now suppose that a Lie group $G$ acts on a hypercomplex manifold
$(M, I_1, I_2, I_3)$ by diffeomorphisms preserving the hypercomplex structure.
In this case, Joyce \cite{MR1109637} defines a hypercomplex moment map to be
an  equivariant map
\begin{eqnarray*}
   \mu = (\mu_1, \mu_2, \mu_3) : M \to \g^*\otimes\R^3
\end{eqnarray*}
satisfying the \lq\lq Cauchy-Riemann equations\rq\rq
\begin{eqnarray}\label{HCQeqns}
   I_1 d\mu_1 = I_2 d\mu_2 = I_3 d\mu_3
\end{eqnarray}
and the \lq\lq transversality condition\rq\rq which states that the function
\begin{eqnarray}\label{HCtrans}
   \langle\lambda, (I_1 d\mu_1)(X_\lambda)\rangle : M \to \R
\end{eqnarray}
does not vanish on $M$ for any non-zero $\lambda \in \g$. When a moment map
exists, the quotient manifold $N = \mu^{-1}(\{\xi\})/G$ admits a natural
hypercomplex structure. We say $N$ has been obtained as a hypercomplex quotient
of $M$ by $G$.
\par
Before we give a proof of the main result in this section, we need to prove
an elementary lemma.
\begin{lemma}\label{HCGaudLem}
   Let $(M, I_1, I_2, I_3, g)$ be a hypercomplex 4-manifold together
   with a metric in the conformal class defined by the hypercomplex structure.
   Let $\omega_1 = g(I_1\cdot, \cdot)$. Then
   \begin{eqnarray*}
      d\omega_1 = \alpha\wedge\omega_1
   \end{eqnarray*}
   for a certain 1-form $\alpha$ (the Lee form, see eg \cite{MR557668})
   independent of $I_1$ satisfying
   \begin{eqnarray}\label{d_s_alpha_eqn}
      d*\alpha = -I_1 dd^c_{I_1} \omega_1
   \end{eqnarray}
   (where $d^c_{I_1} = I_1^{-1}dI_1$).
   \par
   In particular if $g$ is Gauduchon\footnote{Recall that if $(M, J, g)$ is a
   4-manifold with complex structure and Hermitian metric, then we say $g$ is
   Gauduchon (or is a Gauduchon metric) iff $dd^c\omega = 0$ where $\omega =
   g(J\cdot,\cdot)$ and $d^c = J^{-1}dJ$.} with respect to one of the complex
   structures in the hypercomplex 2-sphere, it is Gauduchon with respect to all
   of them.
\end{lemma}
\noindent {\bf Proof} Let $D$ be the Obata connection \cite{MR0095290} on $M$.
As noted in \cite{MR1199072}, $D$ is a Weyl connection and so $Dg = \alpha
\otimes g$ for some 1-form $\alpha$ on $M$. Thus $d\omega_1 = \alpha\wedge\omega_1$
since $D$ is torsion free.
\par
Next, according to \cite{MR2217300} we have
\begin{eqnarray}\label{HCGaudLemma_Step1}
   d*\alpha = \omega_1\wedge d^c_{I_1}\alpha + d^c_{I_1}\omega_1\wedge\alpha
\end{eqnarray}
 for any 1-form $\alpha$ on $M$. Also note
\begin{eqnarray*}
   dd^c_{I_1}\omega_1 &=& -d I_1 d I_1 \omega_1\\
   &=& -d I_1 d\omega_1 = -d I_1 (\alpha\wedge\omega_1)\\
   &=& -d ((I_1\alpha)\wedge\omega_1)\\
   &=& -d(I_1\alpha)\wedge\omega_1 + (I_1\alpha)\wedge\alpha\wedge\omega_1\\
   \Rightarrow d(I_1\alpha)\wedge\omega_1 &=& (I_1\alpha)\wedge\alpha\wedge\omega_1
   - dd^c_{I_1}\omega_1
\end{eqnarray*}
Thus
\begin{eqnarray}\label{HCGaudLemma_Step2}
   \omega_1\wedge d^c_{I_1}\alpha &=& \omega_1\wedge (I_1dI_1\alpha) =
   I_1(d(I_1\alpha)\wedge\omega_1)\notag\\
   &=& I_1((I_1\alpha)\wedge\alpha\wedge\omega_1 - dd^c_{I_1}\omega_1)\notag\\
   \Rightarrow \omega_1\wedge d^c_{I_1}\alpha &=& -\alpha\wedge (I_1\alpha)
   \wedge\omega_1 - I_1dd^c_{I_1}\omega_1
\end{eqnarray}
But 
\begin{eqnarray}\label{HCGaudLemma_Step3}
   d^c_{I_1}\omega_1\wedge\alpha &=& -(I_1d\omega_1)\wedge\alpha\notag\\
   &=& -(I_1(\alpha\wedge\omega_1))\wedge\alpha = -(I_1\alpha)\wedge\omega_1
   \wedge\alpha\notag\\
   \Rightarrow d^c_{I_1}\omega_1\wedge\alpha &=& \alpha\wedge(I_1\alpha)\wedge\omega_1
\end{eqnarray}
Equation \eqref{d_s_alpha_eqn} now follows after substituting \eqref{HCGaudLemma_Step3}
and \eqref{HCGaudLemma_Step2} in \eqref{HCGaudLemma_Step1}. \myqed
\par
Recall \cite{MR742896} that it is always possible to find a Gauduchon metric
in a given conformal class on a complex surface. With this in mind, we state
our main result of this section.
\begin{theorem}\label{HCInstQuot}
   Let $(M, I_1, I_2, I_3)$ be a compact hypercomplex 4-manifold.
   Let $P$ be a principal $SU(2)$ bundle on $M$. Then the moduli
   space of irreducible instantons on $P$ can be formally obtained as an infinite
   dimensional hypercomplex quotient of the set of irreducible connections on $P$ by
   the group of gauge transformations of $P$.
\end{theorem}
\noindent {\bf Proof}
Let $adP$ be the Lie algebra bundle on $M$ associated to $P$ by the adjoint
action of $SU(2)$ on $\mathfrak{su(2)}$. Let $\cG$ be the group of gauge
transformations (ie: automorphisms) of $P$. The Lie algebra of $\cG$ is naturally
$\Omega^0(M, adP)$ and using the pairing
\begin{eqnarray*}
   \Omega^0(M, adP) \otimes \Omega^4(M, adP) &\to& \C\\
   \zeta \otimes \xi &\mapsto& \int_M\limits \langle\zeta, \xi\rangle
\end{eqnarray*}
(the pairing $\langle\cdot,\cdot\rangle$ is defined using the Killing form
of $\mathfrak{su(2)}$) we can naturally identify the dual of the Lie algebra
of $\cG$ with $\Omega^4(M, adP)$.
\par
Let $\cA^*$ be the set of irreducible connections on $P$.  Since $\cA^*$ is an
affine space on $\Omega^1(M, adP)$ we have a natural trivialisation
\begin{eqnarray*}
   T\cA^* \simeq \cA^* \times \Omega^1(M, adP)
\end{eqnarray*}
For $i=1,2,3$, define an almost complex structure $\hat I_i$ on $\cA^*$ by
\begin{eqnarray*}
   \hat I_i : \cA^* \times \Omega^1(M, adP) &\to& \cA^* \times \Omega^1(M, adP)\\
   (A, a) &\to& (A, I_ia)
\end{eqnarray*}
These complex structures obviously satisfy the quaternionic relations and 
make $\cA^*$ into an infinite dimensional hypercomplex manifold. Furthermore
the natural action of $\cG$ on $\cA^*$ preserves this hypercomplex structure.
\par
Now let $g$ be a Gauduchon metric in the conformal class on $M$ defined by the
hypercomplex structure and for $i=1,2,3$ define
\begin{eqnarray*}
   \omega_i = g(I_i\cdot, \cdot) \in \Omega^2(M)
\end{eqnarray*}
We claim
\begin{enumerate}
   \item
   For $i=1,2,3$ the $\cG$-equivariant maps 
   \begin{eqnarray}
      \mu_i : \cA^* &\to& \Omega^4(M, adP)\notag\\
      A &\mapsto& -\omega_i \wedge F_A\label{HCInstMom}
   \end{eqnarray}
   satisfy the \lq\lq Cauchy-Riemann\rq\rq equations (\ref{HCQeqns})
   (where $F_A\in\Omega^2(M, adP)$ is the curvature of a connection $A$)
   \item
   $\mu_i$ also satisfy the \lq\lq transversality condition\rq\rq (\ref{HCtrans}).
\end{enumerate}
To see the first part of the claim, note that for $i=1,2,3$ we have
\begin{eqnarray}\label{HCInstMomDer}
   d\mu_i : \cA^* \times \Omega^1(M, adP) &\to& \Omega^4(M, adP)\notag\\
   (A, a) &\mapsto& \omega_i \wedge d_A a
\end{eqnarray}
where $d_A : \Omega^p(M, adP) \to \Omega^{p+1}(M, adP)$ is the natural extension
of the covariant derivative associated to the connection $A$ to $p$-forms.
Thus
\begin{eqnarray}\label{HCInstHoriz}
   \hat I_i d\mu_i (A, a) = -\omega_i \wedge d_A I_i a
\end{eqnarray}
but since $I_i \eta = - *(\omega_i \wedge \eta)$ for a 1-form $\eta$ on $M$ where
$*$ is the Hodge star operator associated to $g$
\begin{eqnarray*}
   \hat I_i d\mu_i (A, a) = -\omega_i \wedge d_A(*(\omega_i \wedge a))
\end{eqnarray*}
And so if $\alpha$ is the 1-form of lemma \ref{HCGaudLem} we have
\begin{eqnarray*}
   \hat I_i d\mu_i (A, a) &=& -\omega_i \wedge d_A*(\omega_i \wedge a)\\
   &=& -d_A(\omega_i\wedge *(\omega_i\wedge a)) + d\omega_i\wedge *(\omega_i\wedge a)\\
   &=& d_A *a + d\omega_i\wedge *(\omega_i\wedge a)\\
   &=& d_A *a + \alpha\wedge\omega_i\wedge *(\omega_i\wedge a)\\
   &=& d_A *a - \alpha\wedge *a
\end{eqnarray*}
(where we have used the identity $\omega_i\wedge *(\omega_i\wedge \eta) = -*\eta$
twice above). This shows that $\mu_i$ satisfy the equations (\ref{HCQeqns}).
\par
For the second part of the claim, let $\zeta \in \Omega^0(M, adP)$ be non-zero.
The vector field $X_{\zeta}$ on $\cA^*$ generated by $\zeta$ using $\cG$ is
\begin{eqnarray*}
   X_{\zeta}(A) = (A, d_A\zeta) \in \cA^*\times \Omega^1(M, adP)
\end{eqnarray*}
We must show that the function
\begin{eqnarray*}
   f_\zeta : \cA^* &\to& \C\\
   A &\mapsto& \int_M\limits tr(\zeta(d_A(*d_A\zeta) - \alpha\wedge *d_A\zeta))
\end{eqnarray*}
is non-vanishing. But
\begin{eqnarray*}
   dtr(\zeta *d_A\zeta) &=& tr(d_A\zeta\wedge *d_A\zeta) + tr(\zeta d_A *d_A\zeta)\\
   \Rightarrow f_\zeta (A) &=& -\int_M\limits tr(d_A\zeta\wedge *d_A\zeta)
   -\int_M\limits\alpha \wedge tr(\zeta *d_A\zeta)\\
   &=& \| d_A\zeta\|^2 - \int_M\limits \alpha\wedge *tr(\zeta\wedge d_A\zeta)\\
   &=& \| d_A\zeta\|^2 - \frac{1}{2}\int_M\limits \alpha\wedge *dtr\zeta^2\\
   &=& \| d_A\zeta\|^2 - \frac{1}{2}\int_M\limits |\zeta |^2d*\alpha
\end{eqnarray*}
The result now follows since $A$ is irreducible, $\zeta \ne 0$ and $d*\alpha = 0$
by lemma \ref{HCGaudLem} since $g$ is Gauduchon.
\par
Finally note that a two form $\eta$ is anti-self dual if and only if $\omega_i\wedge
\eta = 0$ for $i=1,2,3$. It follows that the hypercomplex quotient $\mu^{-1}(0)/
\cG$ of $\cA^*$ by $\cG$ using the moment maps \eqref{HCInstMom} is the moduli space
of instantons on $P$ and acquires a natural hypercomplex structure. \myqed
\pagebreak
\begin{remark}
\par\noindent
\begin{enumerate}
   \item
   After this proof was completed it was brought to the author's attention that
   Joyce remarks in \cite{MR1163458} that the moduli space of instantons on a
   compact hypercomplex 4-manifold may be obtained as a hypercomplex quotient.
   However no details are provided.
   \item
   It is interesting to look at what model the hypercomplex quotient construction gives
   us for the tangent space to the moduli space of instantons. In \cite{MR1109637}, Joyce
   points out that if $N = \mu^{-1}(\{\xi\})/G$ is a hypercomplex quotient of $M$ by $G$
   using moment map $\mu=(\mu_1, \mu_2, \mu_3)$ then we have a canonical isomorphism
   \begin{eqnarray*}
      T_{\left[ p\right]}N \simeq \left\{ v\in T_pM \quad |\quad d\mu_1(v) = d\mu_2(v)
      = d\mu_3(v) = (I_1d\mu_1)(v) = 0\right\}
   \end{eqnarray*}
   In our case, if $\mathcal{M}^*_P$ is the moduli space of irreducible instantons on
   $P$ then a quick glance at equations \eqref{HCInstMomDer} and \eqref{HCInstHoriz}
   reveals that we have a canonical isomorphism
   \begin{eqnarray*}
      T_{\left[A\right]}\mathcal{M}^*_P \simeq \left\{ a\in \Omega^1(M, adP) \quad |\quad
      d^+_A a = 0 \mbox{ and } \omega_1\wedge d_{I_1 A}^c a = 0 \right\}
   \end{eqnarray*}
   We thus see that the horizontality condition on the tangent space is the same as the
   one used by L\"ubke and Teleman in \cite{MR1370660} (see page 169) and more
   recently by Hitchin in \cite{MR2217300}. Note that this model for the tangent
   space is independent of the choice of Gauduchon metric in the conformal class and
   so therefore is the hypercomplex structure we obtain on the instanton moduli space.
   \item
   Note that the only compact four-dimensional hypercomplex manifolds that are not
   hyperk\"ahler are the Hopf surfaces \cite{MR915736}. Thus it is only for Hopf surfaces
   that theorem \ref{HCInstQuot} gives new information.
   \end{enumerate}
\end{remark}
   Although theorem \ref{HCInstQuot} concerns the geometry of instanton moduli spaces,
   it is related to the geometry of hyperbolic monopole moduli spaces. To see why, we
   follow \cite{MR1013728} and consider the Hopf surface 
   \begin{eqnarray*}
      M = \C^2\setminus \{0\}/\Z
   \end{eqnarray*}
   where the free $\Z$ action is
   \begin{eqnarray*}
      n\cdot(z,w) = (\lambda^n z, \lambda^n w)
   \end{eqnarray*}
   for some fixed $\lambda \in \C^*$, $|\lambda| < 1$. The $S^1$ action on $\C^2$
   defined by
   \begin{eqnarray*}
      e^{i\theta}\cdot(z,w) = (z, e^{i\theta}w)
   \end{eqnarray*}
   commutes with this $\Z$ action and so descends to $M$. On $\C^2$, this $S^1$ action is
   free except on the fixed point set $\C\times\{0\}$. On $M$, the action is free except
   on the quotient
   \begin{eqnarray*}
      B = \C^*\times\{0\}/\Z \simeq S^1\times S^1
   \end{eqnarray*}
   Now it is well known \cite{MR893593} that for appropriate conformal scaling of the
   Euclidean metric on $\C^2$, the $S^1$ action is isometric and we have an isometry
   \begin{eqnarray*}
      (\C^2\setminus \C\times\{0\})/S^1 \simeq H^3
   \end{eqnarray*}
   and so we find
   \begin{eqnarray*}
       (M\setminus B)/S^1 \simeq H^3/\Z \simeq D^2\times S^1
   \end{eqnarray*}
   (where $D^2$ is the unit disc). Thus there is an $S^1$ action on our Hopf surface $M$
   whose quotient away from the fixed point set is a hyperbolic solid torus, the boundary
   of which appears as the fixed point set $B$.
   \par
   This means that $S^1$--invariant instantons on the Hopf surface $M$ correspond to
   monopoles on the hyperbolic three manifold $H^3/\Z$ (just as \cite{MR893593}
   $S^1$--invariant instantons on $S^4$ correspond to monopoles on $H^3$). Since $\partial
   (\Hyp/\Z)$ has one connected component, monopoles on $\Hyp/\Z$ have one charge $k \in \N$
   and (see \cite{MR1010167}) their moduli space\footnote{Note that we
   are not considering a space of gauged monopoles on $\Hyp/\Z$ as we do when we study
   monopoles on $\Hyp$. The space of unbased instantons on the Hopf surface $M$ carries the
   interesting geometry and correspondingly we are interested in the space of monopoles
   with no gauging condition.} has dimension\footnote{Braam's formula \cite{MR1010167}
   states that for a wide class of hyperbolic 3-manifolds, the moduli space of monopoles
   with positive masses and charges $k_i$ has dimension $4\sum k_i - \chi$ (where $\chi$ is
   the Euler characteristic of the hyperbolic manifold).} $4k$.
   \par
   The moduli space of monopoles on $H^3/\Z$ thus appears as the fixed point set under an $S^1$
   action on the moduli space of instantons on $M$, which we have shown carries a natural
   hypercomplex structure. This is analogous to the fact that the moduli space of
   gauged monopoles on $\Hyp$ appears as the fixed point set under an $S^1$ action on
   the moduli space of based instantons on $S^4$, which carries a hyperk\"ahler structure.
   Using this point of view to look at the moduli spaces of monopoles on $\Hyp/\Z$ and $\Hyp$
   we learn that they carry a natural complex structure and K\"ahler metric (see section
   \ref{HypMonop_KahMetric_sect}) respectively.
   Although the result for monopoles on $\Hyp/\Z$ is weaker (merely a complex structure and
   connection) we avoid the issue that gauging a monopole on $\Hyp$ seems to single out a
   point on $\partial\Hyp$ (again, see section \ref{HypMonop_KahMetric_sect}). Furthermore,
   in view of \cite{MR1919716} this is all the geometry we expect to find on the fixed point
   set of an $S^1$ action on a hypercomplex manifold.
\section{Hyperbolic monopoles and K\"ahler metrics}\label{HypMonop_KahMetric_sect}
In this section we shall see that, given a choice of horospherical height function on
$\Hyp$, then, at least in the integral mass case, we can define a natural K\"ahler metric
on the moduli space of hyperbolic monopoles.
\par
This result should be compared
with the results of chapter \ref{SingHypChap} where we showed that the moduli space
of charge $1$ singular hyperbolic monopoles naturally carries a $2$-sphere of
(self-dual) K\"ahler metrics where the $2$-sphere parameterisation naturally factors
through horospherical height functions. In view of the link with horospherical height
functions we expect that the K\"ahler metrics defined here coincide with those of
chapter \ref{SingHypChap} in charge $1$. For higher charge, there appears to be a
difference related to the choice of base point when gauging the hyperbolic monopoles.
\par
Our starting point is the observation that $\Hyp\times S^1$ is just the simplest
case of the space $M$ investigated in section \ref{HWSingU1Sect}. It corresponds to
the case when we have no singular points. In view of this, we know that a point $O
\in \Hyp$ determines a natural $2$-sphere of K\"ahler structures on $\Hyp\times S^1$,
the metrics of which all lie in the same self-dual conformal class. Indeed if $q$ is a
horospherical height function on $\Hyp$ such that the level set $q=1$ contains $O$
then we define a K\"ahler structure on $\Hyp\times S^1$ via the metric and $2$-form
\begin{eqnarray*}
   g &=& q^2(h + d\theta^2)\\
   \Omega &=& \frac{1}{2}\left(*d(q^2) + d(q^2)\wedge d\theta\right)
\end{eqnarray*}
where $h$ is the hyperbolic metric, $*$ is the hyperbolic Hodge $*$-operator
and $\theta$ is the $S^1$ coordinate. As we showed in section \ref{HWSingU1Sect},
such functions $q$ are in 1-1 correspondence with the points of $\partial\Hyp$.
Choosing appropriate coordinates $(x,y,q)$ for $\Hyp$ so that it has metric
$\frac{dx^2 + dy^2 + dq^2}{q^2}$, these equations read
\begin{eqnarray*}
   g &=& dx^2 + dy^2 + dq^2 + q^2d\theta^2\\
   \Omega &=& dx\wedge dy + qdq\wedge d\theta
\end{eqnarray*}
from which we can see that we have merely recovered the usual conformal equivalence
\begin{eqnarray}\label{HypS1ConfEquiv}
   \Hyp\times S^1 \simeq \C\times\C^*
\end{eqnarray}
(see \cite{MR893593}). The $S^1$ action on $\C\times\C^*$ is $e^{i\theta}\cdot(z_1,z_2) =
(z_1,e^{i\theta}z_2)$. The point we wish
to emphasize is that this conformal equivalence singles out a horospherical height
function on $\Hyp$ and so also a point on $\partial\Hyp$.
\par
As is well known, $S^4 = \C^2\cup \{\infty\}$ is an $S^1$-equivariant conformal
compactification of $\C\times\C^*$ and integral mass hyperbolic monopoles extend
to $S^1$-invariant instantons on $S^4$. The $S^1$ action on $S^4$ is free except on
the fixed point set $S = (\C\times\{0\})\cup\{\infty\} \subset S^4$ which is naturally
identified with $\partial\Hyp$. The point we have singled out on $\partial\Hyp$ in making
the conformal identification \eqref{HypS1ConfEquiv} corresponds to $\infty \in S
\subset S^4$.
\par
Consider now the moduli space $X$ of anti-self-dual $SU(2)$ instantons on $E \to S^4$
based at $\infty$. This is a principal $SO(3)$ bundle over the space of unbased instantons
on $S^4$ and is naturally identified with the space of instantons on $\C^2 = S^4 \setminus
\{\infty\}$ framed at $\infty$. As such (see \cite{MR723549} as well as \cite{MR1091573}),
the tangent space at a point $[A] \in X$ can be naturally identified
\begin{eqnarray*}
   T_{[A]}X \simeq L^2_{\C^2} \cap \Ker(d_A^* \oplus d_A^+)
\end{eqnarray*}
where $d_A$ and $d_A^+$ are the usual differentials in the instanton deformation complex
\begin{eqnarray*}
   \Omega^0(\C^2, \End_0(E)) \overset{d_A}\longrightarrow \Omega^1(\C^2, \End_0(E))
   \overset{d_A^+}\longrightarrow \Omega^2_+(\C^2, \End_0(E))
\end{eqnarray*}
and the adjoint $d_A^*$ is defined using the Euclidean metric on $\C^2$. In view of this
$X$ carries a natural $L^2$ metric (which is in fact hyperk\"ahler).
\par
We wish to study hyperbolic monopoles as $S^1$-invariant instantons on $S^4$. We begin
by making sure that we have chosen an $S^1$-equivariant basing of $E$ at $\infty$ (this
cuts down the $SU(2)$ choice of basings to an $S^1$). The space $M$ of $S^1$-invariant
based instantons is the space of gauged hyperbolic monopoles. The $S^1$
action on $X$ preserves the metric and so we have a metric on the monopole moduli space
$M$. If $a \in \Omega^1(\C^2, \End_0(E))^{S^1}$ is an infinitesimal deformation of a
hyperbolic monopole, the length of $a$ is the $L^2$ metric
\begin{eqnarray}\label{HypL2metric}
   \|a\|^2 = \int_{\C^2} (a,a)_E vol_E
\end{eqnarray}
where $(\cdot,\cdot)_E$ uses the Euclidean inner product on $1$-forms and $vol_E$ is the
Euclidean volume form. Since $a$ is $S^1$-invariant, we can express this as an integral
on $\Hyp$. To see this, note that we have an isomorphism
\begin{eqnarray*}
   \Omega^1(\Hyp, \End_0(E/S^1))\oplus \Omega^0(\Hyp, \End_0(E/S^1)) &\simeq&
   \Omega^1(\Hyp\times S^1, \End_0(E))^{S^1}\\
   (a',\psi) &\mapsto& a' + \psi d\theta
\end{eqnarray*}
Thus since the Euclidean metric is a conformal scaling of
the metric on $\Hyp\times S^1$ by a factor of $q^2$ and $a$ is $S^1$-invariant, a quick
calculation reveals that we may also express \eqref{HypL2metric} as
\begin{eqnarray*}
   \|a\|^2 = 2\pi\int_{\Hyp} q^2\left((a',a')_H + (\psi,\psi)\right) vol_H
\end{eqnarray*}
where $(\cdot,\cdot)_H$ uses the hyperbolic inner product on $1$-forms and $vol_H$ is the
hyperbolic volume form. Here $a'$ is an infinitesimal deformation of the monopole
connection on $\Hyp$ and $\psi$ is an infinitesimal deformation of the Higgs field.
\par
We have thus seen that a choice of horospherical height function $q$ defines a metric
on the hyperbolic monopole moduli space. Furthermore, since the $S^1$ action preserves
the complex structure on $\C^2$, $q$ also defines a natural complex structure on the
moduli space so that we in fact have a K\"ahler metric.
\par
In chapter \ref{SingHypChap} we saw that for each horospherical height function $q$
we obtained a (self-dual) K\"ahler metric on the charge $1$ singular monopole moduli
space. After fixing a point $O\in\Hyp$ we used this observation to define a family
of natural K\"ahler metrics on the same moduli space parameterised by $\partial\Hyp$.
In view of the above, it is tempting to think that the same phenomenon occurs here
for moduli spaces of monopoles of arbitrary charge. However
the metrics we have defined here depend on a choice of base point in
$S\simeq\partial\Hyp$ and so these metrics exist on different spaces. There does not
appear to be a natural way to identify spaces of $S^1$-invariant instantons that are
based at different points in $S$. It is thus possible that the results of chapter
\ref{SingHypChap} represent a coincidence that occurs only in charge $1$.
\par
Certainly there is some subtlety when gauging hyperbolic monopoles. On the one hand,
as we have said, when thinking in terms of $S^1$-invariant instantons it appears that
gauging requires us to choose a base point for the instanton in $S$ and so the gauging
depends on a choice of base point in $\partial\Hyp$. However from the point of view of
the spectral curve there is a natural gauging requiring no choice of point on
$\partial\Hyp$: we just choose a trivialisation of the line bundle $L^{2m+k}|_S$
satisfying the appropriate reality condition.
\section{Constructing the twistor space}\label{TwistSpace_Section}
\subsection{Overview}
In chapter \ref{DefSpecChap} we saw how to obtain the moduli space
of charge $k$ monopoles as a space of deformations of genus $(k-1)^2$
curves (spectral curves) in a complex $3$-manifold (the total space
of a certain holomorphic $\C^*$ bundle on a complex surface).
\par
Of course in the Euclidean case it is also possible to obtain the
moduli space as a space of deformations of genus $0$ curves (twistor
lines) in a more complicated space: the twistor space of the moduli
space. This compromise (deforming simple curves in a complicated space
instead of complicated curves in a simple space) is useful because
it enables us to make contact with mainstream twistor theory and
because of the classification of holomorphic vector bundles on $\CP$.
\par
In this section we will show how our approach (deforming the spectral
curve) relates to the usual twistor theoretic approach in the Euclidean
case. The spectral curves appear as $k$-fold branched covers of the
corresponding twistor lines. We will make a conjecture that relates the
normal bundle of the twistor lines to that of the spectral curves in the
appropriate ambient spaces.
\par
We will also indicate what happens in the hyperbolic case. We will find
that it is possible to introduce twistor spaces of the moduli space of
hyperbolic monopoles but that there are complications because it is
difficult to understand the reality conditions on the \lq\lq twistor
lines\rq\rq. Usually a twistor space carries a real structure and the
real structure on the space of deformations of the twistor lines is
derived from this. However this is not what happens in the hyperbolic
case. Thus the approach taken in chapter \ref{DefSpecChap} has a
significant advantage in that the real structure on the space of
deformations does arise in the usual way. Indeed this was one of the
reasons that approach was introduced.
\subsection{The Euclidean case}
We review the construction of the twistor space of the moduli space of
Euclidean monopoles given in \cite{MR934202}. This will motivate our
constructions in the hyperbolic case to follow as well as establishing
the necessary notation so that we can show how the approach introduced
in chapter \ref{DefSpecChap} relates to this more conventional twistor
theoretic approach.
\par
Thus, as in \cite{MR934202}, let
\begin{eqnarray*}
   Y_k = \cO(2) \oplus \cO(4) \oplus \cdots \oplus \cO(2k)
\end{eqnarray*}
(so $Y_k$ is just the fibrewise symmetric product of the line bundle
$T\CP \to \CP$) and let
\begin{eqnarray*}
   D_k = \left\{ (\eta, \eta_1, \ldots, \eta_k) \in T\CP \oplus Y_k \quad|\quad
   \eta^k + \eta_1\eta^{k-1} + \cdots + \eta_k = 0 \right\}
\end{eqnarray*}
(so a point in $D_k$ is a set of $k$ unordered points (possibly with repetition)
in some fibre of the line bundle $T\CP \to \CP$ such that one of the points is
marked). We have the following commutative diagram
\begin{eqnarray}\label{EucTwstDiag}
   \begin{diagram}
      \node[2]{D_k} \arrow{sw,t}{p_1} \arrow{se,t}{p_2}\\
      \node{Y_k} \arrow{se} \node[2]{T\CP} \arrow{sw}\\
      \node[2]{\CP}
   \end{diagram}
\end{eqnarray}
We define the vector bundle $V_k$ on $Y_k$ by
\begin{eqnarray*}
   V_k = p_{1*} p_2^*L^2    
\end{eqnarray*}
and the open set $Z_k \subset V_k$ by
\begin{eqnarray*}
   Z_k = p_{1*} p_2^*(L^2\setminus 0)
\end{eqnarray*}
As shown in \cite{MR934202}, $Z_k$ is the twistor space of the moduli
space of Euclidean monopoles of charge $k$.
\par
Note that the fibre of $Z_k$ over a point $(\eta_1, \ldots, \eta_k)
\in Y_k$ is naturally isomorphic to the direct sum of the fibres
$L^2_{\beta_i} \setminus 0$ for $i=1,2, \ldots, k$ where $\beta_i$ are
the roots of $\beta^k + \eta_1\beta^{k-1} + \cdots + \eta_k = 0$, ie
\begin{eqnarray}\label{Zk_fibre}
   {Z_k}_{(\eta_1, \ldots, \eta_k)} \simeq L^2_{\beta_1}\setminus 0 \oplus
   \cdots \oplus L^2_{\beta_k}\setminus 0
\end{eqnarray}
\par
Now $Z_k$ is a holomorphic fibre bundle over $\CP$. Let us denote the
map by
\begin{eqnarray*}
   p : Z_k \to \CP
\end{eqnarray*}
Since $Z_k$ is the twistor space of a hyperk\"ahler manifold, there
is a family of holomorphic sections of $p$. It is easy to see how these
correspond to curves $S \subset T\CP$ in the linear system $|\cO(2k)|$
together with a trivialisation of $L^2|_S$ as follows. First note that
a section $s$ of $p$ includes a section $s'$ of $Y_k \to \CP$. Let
$S' \subset Y_k$ be the image of $s'$. Now if $\eta$ is the
tautological section of $\cO(2)$ on $T\CP$ and
\begin{eqnarray*}
   s' = (a_1, \ldots, a_k)
\end{eqnarray*}
$a_i \in H^0(\CP, \cO(2i))$, then pulling back $a_i$ to $T\CP$ we let
\begin{eqnarray*}
   \psi = a_1\eta^{k-1} + \cdots + a_k \in H^0(T\CP, \cO(2k))
\end{eqnarray*}
The divisor $S$ of $\psi$ is the curve corresponding $s$. To get the
trivialisation of $L^2$ on $S$, we simply use \eqref{Zk_fibre} and the
lifting of $S' \subset Y_k$ to $\tilde S = im s \subset Z_k$. It is useful
to summarise this correspondence in the following diagram
\begin{eqnarray}\label{Euc_corresp_diag}
   \begin{array}{ccccccccc}
      \tilde S &\subset & Z_k        &                      &D_k &                      &L^2\setminus 0  &\supset & \hat S\\
               &        & \downarrow &\overset{p_1}\swarrow &    &\overset{p_2}\searrow &\downarrow      &        &\\
      S'       &\subset & Y_k        &                      &    &                      &T\CP            &\supset &S\\
               &        &            &\searrow              &    &\swarrow              &                &\\
               &        &            &                      &\CP &                      &                &
\end{array}
\end{eqnarray}
where $\hat S$ is the lifting of $S \subset T\CP$ to $L^2\setminus 0$
provided by the trivialisation. Note also that we have
\begin{eqnarray*}
   S = p_2(p_1^{-1}(S'))
\end{eqnarray*}
that $p_2 : p_1^{-1}(S') \to S$ is a biholomorphism and that $p_1 :
p_1^{-1}(S') \to S'$ is a $k$-fold branched covering. Identifying $\tilde S$
with $S'$ and $\hat S$ with $S$ we see that the curves $\hat S$ we were
deforming in chapter \ref{DefSpecChap} are naturally $k$-fold branched covers
of the lines $\tilde S$ in $Z_k$.
\par
An important property of the twistor space $Z_k$ is that it carries a
real structure
\begin{eqnarray*}
   \sigma : Z_k \to Z_k
\end{eqnarray*}
defined as follows. Let $x \in Z_k$ and using \eqref{Zk_fibre} let $x = (x_1,
\ldots, x_k)$ where $x_i \in L^2_{\beta_i} \setminus 0$. $L^2\setminus 0$
carries a real structure $\tau$ say (indeed it is the twistor space of the
real 4-manifold $\R^3\times S^1$) and we define the real structure $\sigma$
on $Z_k$ by
\begin{eqnarray*}
   \sigma(x_1, \ldots, x_k) = \left(\tau(x_1), \ldots, \tau(x_k)\right)
\end{eqnarray*}
Note that $\sigma$ covers the natural real structure on $Y_k$ and $p$ commutes
with $\sigma$ and the usual real structure on $\CP$.
\par
Now, in order to show that $Z_k$ is the twistor space of a hyperk\"ahler manifold,
it is necessary to prove that the normal bundle of a section of $p$ is
isomorphic to $\C^{2k} \otimes \cO(1)$ and to exhibit a holomorphic section
$\omega$ of $\wedge^2T_F^*\otimes \cO(2)$ defining a symplectic form on each
fibre of $p$ (where $T_F = \Ker p_*$ is the tangent bundle to the fibres of
$p$). In \cite{MR934202}, these two goals are achieved by passing to a more
explicit patching description of $Z_k$ (ie showing how to construct $Z_k$
by gluing two copies of $\C^{2k+1}$ together appropriately). We shall not need
these explicit details here and instead point out that it should be possible to
establish these results intrinsically using the ideas below (in particular
conjecture \ref{NrmBnd_conj_push_frwd}).
\par
We have already noted that the curves $\hat S$ we deform in $L^2\setminus 0$ are
naturally $k$-fold branched covers of the lines $\tilde S$ that are deformed in $Z_k$
and so we have begun to understand the connection between these two spaces of
deformations. In order to complete the picture, we need to relate the normal
bundles of each of these curves in their ambient spaces. We begin by noting the
following two results.
\begin{proposition}\label{HyeonLemma}
   Let $S \subset T\CP$ be the spectral curve of a charge $k$ Euclidean monopole.
   Let $\pi_S : S \to \CP$ be the natural degree $k$ branched covering.
   Then there is a natural isomorphism
   \begin{eqnarray*}
       \pi_{S*}\cO \simeq \cO \oplus \cO(-2) \oplus \cdots \oplus \cO(2-2k)
   \end{eqnarray*}
\end{proposition}
\noindent{\bf Proof} See \cite{MR1841397}. \myqed
\par
(Note that this lemma is a useful way to encode the content of proposition 3.1,
proposition 4.5 and lemma 5.2 of \cite{MR709461}.)
\begin{proposition}\label{Norm_push_frwd}
   Let $S \subset T\CP$ be the spectral curve of a charge $k$ Euclidean monopole.
   Let $\pi_S : S \to \CP$ be the natural degree $k$ branched covering.
   Let $\hat N$ be the normal bundle of a lift of $S$ to $L^2\setminus 0$.
   Then
   \begin{eqnarray*}
      \pi_{S*}\hat N \simeq \C^{2k}\otimes \cO(1)
   \end{eqnarray*}
\end{proposition}
\noindent {\bf Proof} We must show that $\pi_{S*}\hat N(-1)$ is trivial.
It follows easily from the classification of vector bundles on $\CP$ that
a vector bundle $V$ on $\CP$ is trivial if and only if
\begin{eqnarray*}
   h^0(\CP, V) &=& \rank V\\
   h^0(\CP, V(-1)) &=& 0
\end{eqnarray*}
We thus need to show
\begin{eqnarray*}
   h^0(\CP, \pi_{S*}\hat N(-1)) &=& 2k\\
   h^0(\CP, \pi_{S*}\hat N(-2)) &=& 0
\end{eqnarray*}
ie
\begin{eqnarray*}
   h^0(S, \hat N(-1)) &=& 2k\\
   h^0(S, \hat N(-2)) &=& 0
\end{eqnarray*}
and we have already noted on several occasions that both of these results
are implied by theorem \ref{Euc_vanishing_theorem}. \myqed
\par
Now, referring to diagram \eqref{Euc_corresp_diag} let $N$, $\hat N$, $N'$
and $\tilde N$ be the normal bundles of $S$, $\hat S$, $S'$ and $\tilde S$
in their ambient spaces. Then we have a natural isomorphism
\begin{eqnarray*}
   N' \simeq p_{1*} p_2^* N
\end{eqnarray*}
To see why we need only note
\begin{itemize}
   \item
      $p_1 : p_1^{-1}(S^{'}) \to S^{'}$ is naturally isomorphic to the
      $k$-fold branch cover $\pi_S : S \to \CP$ obtained from the projection
      $\pi : T\CP \to \CP$
   \item
      $N' \simeq Y_k = \cO(2) \oplus \cdots \oplus \cO(2k)$ naturally since
      $S'$ is the image of a section of $Y_k$
   \item
      $N \simeq \cO(2k)$ naturally since $S$ is the divisor of a section of
      $\cO(2k)$
   \item
      By proposition \ref{HyeonLemma} $\pi_{S*} \cO(2k) \simeq \cO(2) \oplus
      \cdots \oplus \cO(2k)$ naturally.
\end{itemize}
Now note that (in view of lemma \ref{NormExtLemma}) we have the following
natural exact sequence on $S'$
\begin{eqnarray}\label{NrmBndExt3}
   0 \to V_k|_{S'} \to \tilde N \to N' \to 0
\end{eqnarray}
However pulling back using $p_2$ and applying the direct image under $p_1$
to
\begin{eqnarray*}
   0 \to L^2|_S \to \hat N \to N \to 0
\end{eqnarray*}
we also have the following natural exact sequence of bundles on $S'$
\begin{eqnarray}\label{NrmBndExt_push_frwd}
   0 \to p_{1*}p_2^* L^2|_S \to p_{1*}p_2^* \hat N \to p_{1*}p_2^* N \to 0
\end{eqnarray}
(the final $0$ follows by counting ranks). Since we have just noted that we
have a natural isomorphism $N' \simeq p_{1*}p_2^*N$ and also by definition
$V_k = p_{1*}p_2^* L^2$ we make the following
\begin{conjecture}\label{NrmBnd_conj_push_frwd}
   The sequences \eqref{NrmBndExt3} and \eqref{NrmBndExt_push_frwd} are
   naturally isomorphic.
\end{conjecture}
Note that since $\tilde N$ is the normal bundle of a twistor line, we know
that $\tilde N \simeq \C^{2k}\otimes \cO(1)$ and so proposition
\ref{Norm_push_frwd} says that $\tilde N$ and $p_{1*}p_2^* \hat N$ are
isomorphic bundles.
\par
Note also that a natural isomorphism $p_{1*}p_2^*\hat N \simeq \tilde N$
would induce a natural isomorphism of sections
\begin{eqnarray*}
   H^0(\tilde S, \tilde N) \simeq H^0(\hat S, \hat N)
\end{eqnarray*}
which we expect since they are both the complexified tangent spaces to the
monopole moduli space and also a natural isomorphism
\begin{eqnarray}\label{Euc_skewformspace_iso}
   H^0(\tilde S, \tilde N(-1)) \simeq H^0(\hat S, \hat N(-1))
\end{eqnarray}
and this isomorphism should carry the symplectic structure on $H^0(\tilde S,
\tilde N(-1))$ that exists because $Z_k$ is the twistor space of a
hyperk\"ahler manifold to the one we defined on $H^0(\hat S, \hat N(-1))$.
\subsection{The hyperbolic case}
In order to construct the analogue of the twistor space $Z_k$ for the hyperbolic
monopole moduli space, we first need to define analogues of the spaces $Y_k$ and
$D_k$. In the Euclidean case, $Y_k$ was obtained by taking the fibrewise symmetric
product of $T\CP \to \CP$. $T\CP$ is of course replaced with $Q = \QQ$ in the
hyperbolic case. However $\QQ$ has two maps to $\CP$, $\pi_L$ and $\pi_R$ (the
projections onto the left and right factors of $\QQ$ respectively). We cannot
favour either factor and so we have a version of the constructions in the Euclidean
case for each factor. We thus define
\begin{eqnarray*}
   Y_k^L = \CP\times s^k\CP \setminus \Delta_{Y^L}
\end{eqnarray*}
where
\begin{eqnarray*}
   \Delta_{Y^L} = \left\{ (p, \{q_1, \ldots, q_k\}) \in \CP\times s^k\CP
   \quad|\quad \tau(p) \in \{q_1, \ldots, q_k\}\right\}
\end{eqnarray*}
Note that we're representing points in the symmetric product $s^k\CP$ as multisets
$\{q_1, \ldots, q_k\}$.
\par
Next define
\begin{eqnarray*}
   D_k^L = \left\{(p, \{q_1, \ldots, q_k\}, r) \in Y_k^L\times\CP \quad|\quad r \in
   \{q_1, \ldots, q_k\}\right\}
\end{eqnarray*}
There are obvious similar definitions of the spaces $Y_k^R$ and $D_k^R$. The diagram
corresponding to \eqref{EucTwstDiag} is
\begin{eqnarray}\label{TwstCorrDiag}
   \begin{diagram}
      \node[2]{D_k^L}\arrow{sw,t}{p_1^L}\arrow{se,t}{p_2^L} \node[2]{D_k^R}\arrow{se,t}
      {p_1^R}\arrow{sw,t}{p_2^R}\\
      \node{Y_k^L}\arrow{se} \node[2]{\Q}\arrow{se}\arrow{sw} \node[2]{Y_k^R}\arrow{sw}\\
      \node[2]{\CP} \node[2]{\CP}
   \end{diagram}
\end{eqnarray}
Finally we define the vector bundle
\begin{eqnarray*}
   V_k^L = {(p_1^L)}_* {(p_2^L)}^*L^{2m+k}
\end{eqnarray*}
and the open set $Z_k^L \subset V_k^L$
\begin{eqnarray*}
   Z_k^L = {(p_1^L)}_* {(p_2^L)}^*(L^{2m+k}\setminus 0)
\end{eqnarray*}
Note that $Z_k^L$ is a holomorphic fibre bundle over $\CP$. An
appropriate section $s$ of $Z_k^L \to \CP$ yields a curve $S \subset
\QQ$ in the linear system $|\cO(k, k)|$ together with a trivialisation of
$L^{2m+k}|_S$ as follows. A section $s$ of $Z_k^L$ includes a section
of $Y_k^L$ which is really just a map
\begin{eqnarray*}
s' : \CP \to s^k\CP
\end{eqnarray*}
If we are to obtain a curve in the linear system $|\cO(k, k)|$, we must
only consider those sections $s$ for which $s'$ is a degree $k$ map.
Now use the natural isomorphism $s^k\CP \simeq \mathbb{CP}^k$ to represent
$s'$ as
\begin{eqnarray*}
   s' = [s'_0, s'_1, \ldots, s'_k]
\end{eqnarray*}
where $s'_i$ (defined up to a factor of $\C^*$, the same for each $i$) is
a polynomial of degree $k$. We define the section $\psi$ of $\cO(k, k)$ on
$\QQ$ by
\begin{eqnarray*}
   \psi(\zeta_L, \zeta_R) = s'_0(\zeta_L) + s'_1(\zeta_L)\zeta_R + \cdots +
   s'_k(\zeta_L) \zeta_R^k
\end{eqnarray*}
where $\zeta_L$, $\zeta_R$ are coordinates on the left and right factors
in $\QQ$. $S \subset \QQ$ is defined to be the divisor of $\psi$. We
obtain the trivialisation of $L^{2m+k}|_S$ as in the Euclidean case from
the natural isomorphism that identifies the fibre of $Z_k^L$ over a point
$(p, \{q_1, \ldots, q_k\}) \in Y_k^L$
\begin{eqnarray}\label{hyp_Zk_fibre}
   {(Z_k^L)}_{(p, \{q_1, \ldots, q_k\})} \simeq L^{2m+k}_{(p,q_1)}\setminus 0
   \oplus \cdots \oplus L^{2m+k}_{(p,q_k)}\setminus 0
\end{eqnarray}
and from the fact that the section of $Y_k^L$ defined by $s'$ lifts to a
section of $Z_k^L$.
\par
Note also that if $S' \subset Y_k^L$ is the image of the section of $Y_k^L$
defined by $s'$ then we have
\begin{eqnarray*}
   S = p_2^L(p_1^L)^{-1}(S')
\end{eqnarray*}
so that identifying $S'$ with the line above it in $Z_k^L$ and $S$ with its
lift to $L^{2m+k}$, the spectral curve we deformed in $L^{2m+k}$ in chapter
\ref{DefSpecChap} is naturally a branched cover of the twistor line in $Z_k^L$
(just as in the Euclidean case).
\par
Again, there are similar constructions of $V_k^R$ and $Z_k^R$.
\par
Now in the Euclidean case, the twistor space came with a real structure. Here
we meet a fundamental difference between the Euclidean and hyperbolic cases. Instead
of a real structure $\sigma : Z_k \to Z_k$, we have an anti-holomorphic bijection
\begin{eqnarray*}
   \sigma : Z_k^L \to Z_k^R
\end{eqnarray*}
defined as follows. Use the isomorphism \eqref{hyp_Zk_fibre} for $Z_k^L$ to
represent a point in $Z_k^L$ as $(x_1, \ldots, x_k)$ with $x_i \in L^{2m+k}
\setminus 0$. Then, using the version of \eqref{hyp_Zk_fibre} for $Z_k^R$, define
\begin{eqnarray*}
   \sigma(x_1, \ldots, x_k) = (\tau(x_1), \ldots, \tau(x_l))
\end{eqnarray*}
where $\tau$ is the real structure on $L^{2m+k}\setminus 0$. Note that $\sigma$ covers
the natural anti-holomorphic bijection $\sigma_Y : Y_k^R \to Y_k^L$ defined by
\begin{eqnarray*}
   \left(p, \{q_1, \ldots, q_k\}\right) \mapsto \left(\{\tau(q_1), \ldots, \tau(q_k)\},
   \tau(p)\right)
\end{eqnarray*}
So instead of the pair $(Z_k, \sigma)$ we have the triple $(Z_k^L, Z_k^R, \sigma)$. This
is an important difference. In the Euclidean case, $\sigma : Z_k \to Z_k$
was not just any anti-holomorphic bijection, it was a real structure, ie: it
was an involution. Of course this condition does not make sense in the case $\sigma : Z_k^L
\to Z_k^R$ and it is not clear what it should be replaced with.
\par
It is worth pointing out that this phenomenon with real structures has been seen before in
twistor theory. The simplest family of examples is the twistor spaces of oriented conformal
$(4n+2)$-manifolds. The conformal $6$-sphere $S^6$ is probably the simplest example. As in
our situation, rather than obtaining a twistor space with a real structure one obtains two
twistor spaces with an anti-holomorphic bijection between them. (See \cite{MR1397410} for a
nice account of this.) Even in this case, it is not clear exactly what other data we
must supply with these conjugate twistor spaces in order to be able to recover the original
real $S^6$.
\par
The point is that given a twistor space with an embedded compact submanifold that admits
deformations, we obtain a complex manifold of deformations. If we are to recover a real manifold
we need this family of deformations to carry a real structure. If the underlying twistor space
itself carried a real structure then the family of deformations will too. As we shall see (just
like in the case of $S^6$) we find ourselves in the case where we can describe a real structure
on the family of deformations but we cannot see how it arises from the more fundamental twistor
spaces.
\par
Of course we can tell which sections of $Z_k^L \to \CP$ we would like to call real
because we know they correspond to curves in $L^{2m+k}\setminus 0$ which does have
a real structure. ie if $s^L$ is a section of $Z_k^L \to \CP$ corresponding to
$\hat S \subset L^{2m+k}\setminus 0$ then we say $s^L$ is real iff $\hat S$ is.
Also, consider the following commutative diagram
\begin{eqnarray*}
   \begin{CD}
      Z_k^L @>\sigma>> Z_k^R\\
      @VVV @VVV\\
      \CP @>\tau>> \CP
   \end{CD}
\end{eqnarray*}
From this, $s^L$ defines a section $s^R = \sigma \circ s^L \circ \tau$ of
$Z_k^R \to \CP$. However, $s^L$ also defines a section of $Z_k^R \to \CP$ by
using $\hat S$. $\hat S$ is real iff these two sections coincide.
\par
Finally, we say something about the symplectic forms carried by the hyperbolic
monopole moduli space. Indeed we make the following
\begin{conjecture}
   $Z_k^L$ carries a natural holomorphic section
   \begin{eqnarray*}
      \omega^L \in H^0(Z_k^L, \wedge^2 T_F^*\otimes \cO(2))
   \end{eqnarray*}
   defining a symplectic structure on the fibres of $p^L : Z_k^L \to \CP$
   (where $T_F = \Ker p^L_*$ is the bundle of tangent vectors to the fibres of
   $p^L$). Similarly for $Z_k^R$.   
\end{conjecture}
As evidence for this conjecture let $\tilde S^L \subset Z_k^L$ be the image
of a section of $p^L : Z_k^L \to \CP$ corresponding to the spectral curve
$\hat S \subset L^{2m+k}\setminus 0$. Let the normal bundles of these curves in
their ambient spaces be $\tilde N^L$ and $\hat N$ respectively. Then we expect
that there should be a natural isomorphism
\begin{eqnarray*}
   H^0(\tilde S^L, \tilde N^L(-1)) \simeq H^0(\hat S, \hat N(-1))
\end{eqnarray*}
analogous to \eqref{Euc_skewformspace_iso}. However $\tilde S^L$ being the
image of a section, we also have the natural isomorphism
\begin{eqnarray*}
   \tilde N^L \simeq T_F|_{\tilde S^L}
\end{eqnarray*}
Since we saw in chapter \ref{DefSpecChap} that $H^0(\hat S, \hat N(-1))$
carries a natural symplectic form, we thus expect to find one on $H^0(\tilde S^L,
\tilde N^L(-1)) \simeq H^0(\tilde S^L, T_F(-1))$ and inspired by what happens
for the twistor spaces of hyperk\"ahler manifolds we expect this symplectic form
should be induced by a section $\omega^L$ as conjectured above.

\addcontentsline{toc}{chapter}{Bibliography}

\bibliography{thesis}
\bibliographystyle{plain}
\end{document}